\documentclass[a4paper,12pt,centertags]{amsart}
\usepackage{pslatex}
\usepackage{mathrsfs}
\usepackage{dsfont}
\usepackage{hyperref}
\newtheorem{theorem}{Theorem}[section]
\newtheorem{lemma}[theorem]{Lemma}
\newtheorem{proposition}[theorem]{Proposition}
\newtheorem{corollary}[theorem]{Corollary}
\newtheorem{assumption}[theorem]{Assumption}
\theoremstyle{definition}
\newtheorem{definition}[theorem]{Definition}
\theoremstyle{remark}
\newtheorem{remark}[theorem]{Remark}
\numberwithin{equation}{section}

\DeclareMathOperator{\sgn}{sgn}
\DeclareMathOperator{\Span}{span}
\DeclareMathOperator{\comp}{Comp}
\DeclareMathOperator{\Div}{div}
\DeclareMathOperator{\Tr}{Tr}
\newcommand{\E}{\mathbb{E}}
\newcommand{\Q}{\mathcal{Q}}
\newcommand{\abP}{\mathds{P}} 
\newcommand{\navsto}{{\text{\Tiny NS}}}
\newcommand{\R}{\mathbf{R}}
\newcommand{\N}{\mathbf{N}}
\newcommand{\Z}{\mathbf{Z}}
\newcommand{\T}{\mathds{T}}    
\renewcommand{\H}{\mathcal{H}} 
\newcommand{\V}{\mathcal{V}}   
\newcommand{\W}{\mathcal{W}}   
\newcommand{\B}{\mathscr{B}}
\newcommand{\Ps}{\mathcal{P}}  
\newcommand{\Bns}{\B_\navsto}
\newcommand{\BNS}{\B^\navsto}
\newcommand{\ms}{\mathscr{C}}
\newcommand{\ns}{\ms_\navsto}
\newcommand{\rv}{\mathcal{R}^+}
\newcommand{\F}{\mathscr{F}}
\newcommand{\test}{\mathcal{D}^\infty}
\newcommand{\Om}{\Omega}
\newcommand{\Omns}{\Om_\navsto}
\newcommand{\e}{\mathtt{e}}
\newcommand{\ka}{\mathbf{k}}
\newcommand{\el}{\mathbf{l}}
\newcommand{\ma}{\mathbf{m}}
\newcommand{\set}[2]{\left\{\,#1\,|\,#2\,\right\}}
\newcommand{\ps}[2]{\langle#1,#2\rangle}
\newcommand{\sss}[1]{{\text{\tiny $#1$}}}
\newcommand{\uno}{\mathds{1}}
\newcommand{\s}{\sigma}
\newcommand{\om}{\omega}
\newcommand{\la}{\lambda}
\newcommand{\al}{\alpha}
\newcommand{\oth}{\overline{\theta}}
\newcommand{\be}{\beta}
\newcommand{\ep}{\varepsilon}
\newcommand{\oep}{{\overline{\kappa}}}
\renewcommand{\th}{\theta}
\newcommand{\loc}{{\text{\Tiny loc}}}
\newcommand{\map}[1]{\textbf{\textsf{\footnotesize [MP#1]}}}
\newcommand{\er}[1]{{\text{\Tiny ($#1$)}}}
\newcommand{\ov}[1]{\overline{#1}}
\newcommand{\wt}[1]{\tilde{#1}}
\hypersetup{pdftitle={Markov selections for the 3D stochastic Navier-Stokes equations},
pdfsubject={Markov selections for the Navier-Stokes equations},
pdfauthor={F. Flandoli, M. Romito},
pdfkeywords={stochastic Navier-Stokes equations, martingale problem, Markov property,
Markov selections, strong Feller property, well posedness},
pdfcreator={M. Romito}}
\begin{document}
\title[Markov selections for the Navier-Stokes equations]{Markov selections
for the 3D stochastic Navier-Stokes equations}
\author[F. Flandoli]{Franco Flandoli}
\address{Dipartimento di Matematica Applicata, Universit\`a di Pisa,
         via Bonanno Pisano 25/b, 56126 Pisa, Italia}
\email{flandoli@dma.unipi.it}
\author[M. Romito]{Marco Romito}
\address{Dipartimento di Matematica, Universit\`a di Firenze,
         Viale Morgagni 67/a, 50134 Firenze, Italia}
\email{romito@math.unifi.it}
\subjclass[2000]{Primary 76D05; Secondary 60H15, 35Q30, 60H30, 76M35}
\keywords{stochastic Navier-Stokes equations, martingale problem, Markov property,
Markov selections, strong Feller property, well posedness}
\begin{abstract}
We investigate the Markov property and the continuity with respect
to the initial conditions (strong Feller property) for the solutions
to the Navier-Stokes equations forced by an additive noise.

First, we prove, by means of an abstract selection principle, that there
are Markov solutions to the Navier-Stokes equations. Due to the lack of
continuity of solutions in the space of finite energy, the Markov property
holds almost everywhere in time. Then, depending on the regularity of the
noise, we prove that any Markov solution has the strong Feller property
for regular initial conditions.

We give also a few consequences of these facts, together with a new
sufficient condition for well-posedness.
\end{abstract}
\maketitle
\renewcommand{\tocsubsection}{\quad\tocsection} 
\tableofcontents
\section{Introduction}
\subsection*{General overview}
The well posedness of 3D Navier-Stokes equations (or, briefly, NSE) is still
an open problem, both in the well known deterministic case, see Fefferman
\cite{Fef}, Temam \cite{Tem3} for reviews, and in the case of stochastic
perturbations. Weak solutions (suitably defined in the stochastic case) exist
globally in time but their uniqueness is not known. Suitably regular solutions
are unique but they are proved to exist only locally in time, for regular data.
The very strong theorems of uniqueness for stochastic ordinary equations yield
hope that white noise perturbations may help, but the question is still open.

The first difference between the deterministic and stochastic case appears on
the question of continuous dependence on initial conditions, the third
property after existence and uniqueness in Hadamard definition of well
posedness. When uniqueness is open, there is no question of continuous
dependence in a strict sense, but one may ask whether there exists a
\emph{continuous selection}. The existence of a continuous selection is not
known for the 3D deterministic NSE.
One of the main results which are true for 3D \emph{stochastic} NSE is
the existence of a continuous selection when the noise is sufficiently
non degenerate. This result was first proved by Da Prato \& Debussche
\cite{DPDe}. The main aim of the present paper is to give a new insight
into this problem. We give an entirely different proof which works in
greater generality (see Theorem \ref{t:maincont}), in particular we
prove that \emph{every} Markov selection depends continuously on initial
conditions, for suitable noise.

Indeed, a fourth structural property, beside existence, uniqueness and continuous
dependence, of stochastic differential equations is the Markov property. When
uniqueness is open, Markov property has no direct meaning but a natural
question is the existence of a \emph{Markov selection}. Another main result
of this paper is the existence of a Markov selection (Theorem \ref{t:mainmarkov})
for a very large class of 3D stochastic NSE, which in fact includes also
the deterministic case.

As remarked above, the two previous questions, namely existence of Markov selection
and continuous selection, are not unrelated. The strategy of our approach is
first to prove the existence of a Markov selection in great generality on the
additive noise (even zero noise is acceptable, see Assumption \ref{a:noiseass1}),
then to prove under strong restrictions on the noise (see Assumptions \ref{a:noiseass2})
that \emph{every} Markov selection has a property of continuous dependence
on initial conditions, a property of strong Feller type. 

An obvious question is whether the previous results have consequences on the
uniqueness of weak solutions. There are many facets of this question, but
unfortunately we have not found any true result. In principle, existence of
at least \emph{one} continuous Markov selection may be a basic step, since there
are examples of differential equations in the literature where uniqueness of
solutions is proved by means of one regular flow, both in the deterministic
and stochastic case. But suitable estimates, not available at present, on
the derivative of the flow in the initial conditions seem to be necessary.

The strong Feller property for every Markov selection does not imply
uniqueness, see Stroock \& Yor \cite{StYo}.
If one replaces the sentence ``Markov selection'' by
``measurable selection'' the answer would be positive (see Flandoli
\cite{Fla2}) but the Markov property is a very demanding one. Indeed,
under irreducibility, a single solution contains information
(by disintegration) on most of the others and under strong Feller on
all the others. This rigidity is one of the obstacles in any
attempt to deduce uniqueness. 

We have only one positive example of consequence in the direction of
uniqueness, namely a conditional theorem (see Section \ref{ss:wellposed}).
It holds true for the models where we have the strong Feller Markov selections,
which are also irreducible. Roughly it states that if the problem is well posed
for one initial condition then it is well posed for all initial conditions.
This is a dichotomy with respect to the deterministic case, where well posedness
is known for sufficiently small and regular initial conditions, a result that
cannot be extended to to the stochastic case with additive non degenerate noise
due to the absence of small invariant sets.

Throughout the paper, we shall consider the Navier-Stokes equations
on the torus $[0,1]^3$,
\begin{align*}
&\dot u+(u\cdot\nabla)u+\nabla p=\nu\Delta u+\dot\eta,\\
&\Div u=0,
\end{align*}
with periodic boundary conditions, driven by a random force (details on
the equations will follow later, see in particular the main assumptions
\ref{a:noiseass1} and \ref{a:noiseass2} of the paper). Other boundary
conditions could be analysed, but we focus on this setting because
it is the simplest in this framework.
\subsection*{Some details on the main results}
There are interesting details (and restrictions) about the existence of a
Markov selection for this problem. To prove the
existence of such a selection, we need to use a definition of weak solution
which incorporates certain energy inequalities (this is not at all surprising
in the theory of 3D NSE, since one cannot prove directly that weak solutions
satisfy energy inequalities). These energy inequalities are necessary to our
approach to prove certain compactness results at due time. However, they
introduce a technical difficulty. All conditions included in the definition of
solution must be stable under the operations required by the Markov property,
namely disintegration and reconstruction.
Thus we express all the needed energy inequalities in the form of super-martingale
properties. This is a novelty on 3D stochastic NSE, to our knowledge. In a sentence,
we translate the usual well known \emph{ energy inequality} of the deterministic
case in a \emph{super-martingale property} for the stochastic case.
There is a non-trivial gain of information in the super-martingale formulation
of energy inequality, for instance it implies stopped energy inequalities
(see the proof of Theorem \ref{t:weakstrong}).

Another important detail is that the properties included in the definition should
be invariant under time translation, since disintegration and reconstruction need
to be applied at any time $s>0$. Unfortunately, even in the deterministic
case the energy inequality is known to hold only almost surely in time. More
precisely, if one writes the energy inequality between two generic times $t>s>0$,
$t$ can be arbitrary but only a.\ e.\ $s$ is allowed. This difficulty cannot be
overcome at the present state of understanding of the NSE, see, for instance
Constantin, E \& Titi \cite{CET} and Duchon \& Robert \cite{DuRo}
on this and related technical problems.

As a consequence, we can prove only an \emph{almost sure} Markov property.
We have the impression that this is not a drawback of our approach but an
intrinsic difficulty. However, under a strong assumption on the noise
\ref{a:noiseass2} that give us the strong Feller
property, we can prove that the Markov property holds true for every time.

Concerning the strong Feller property, there are several details here too that
could be highlighted. First, the topology on the initial conditions required
for the continuous dependence is not the one of the energy space (called $H$
below) but is a more regular topology related to the assumptions \ref{a:noiseass2}
on the covariance of the noise. However, we conjecture that, with suitable
improvements of some of the arguments given here, one can prove the continuous
dependence in a suitable space of regularity depending only on the Stokes
operator, independently of the assumptions on the covariance (but non
degeneracy). This will be the object of future research.
\subsection*{Comparison with the literature}
Existence of weak solutions of the martingale problem for quite general
3D stochastic NSE is classical, see for instance Flandoli \cite{Fla2}
and the references therein. Here we introduce a new definition with special
energy inequalities (Definition \ref{d:asNSsol}), so we give a few details
of proofs in the appendices.

The existence of Markov selections for certain classes of stochastic ordinary
equations is due to Krylov \cite{Kry}. We have generalised the abstract part
of this result to Hilbert spaces, following closely the presentation of Stroock
\& Varadhan \cite{StVa}.

As regards the central results of this paper, the inspiration comes undoubtedly
from the basic work of Da Prato \& Debussche \cite{DPDe}. They have proved, among
other facts, the existence of a strong Feller selection build on the Galerkin
scheme (thus rather constructive) and several properties of the associated
Markov semi-group. In a recent paper, Debussche \& Odasso \cite{DeOd} have
proved also that the selection is Markov itself.

Our approach is different.
Our Markov selection procedure is less constructive since it is based on quite
abstract notions of solutions and then on the minimisation of quite generic
functionals. But it works in great generality on the noise, up to the zero
noise case (see Assumption \ref{a:noiseass1}), while the results based
or related to strong Feller property require a suitable non degenerate noise
(see Assumption \ref{a:noiseass2}). At the end, under such stronger
assumptions, we recover a result of type \cite{DPDe}, in the general sense
that \emph{every} Markov selection with suitable noise is strong Feller. 

Other differences with respect to \cite{DPDe} (and \cite{DeOd}) are concerned
with the assumptions on the noise -- we explicitly deal with a larger class of
covariances, but restricted to space periodic case (the technical effort is
non trivial and we had to develop estimates on the nonlinear operator $B$ that
seem to be new, see Appendix \ref{s:nonlinear}). The idea that one can
shift the problem at any level of the Hilbert scale associated to the Stokes
operator $A$ is known and also discussed in private conversations with
A.~Debussche (see also a related work on ergodicity by Ferrario \cite{Fer}).
Another, more substantial, difference is the proof of the strong Feller
property. The proof given here is very short and direct, immediately based
on the Markov property, and it shows more transparently why the strong Feller
property holds. The proof in \cite{DPDe} is based on a completely different
argument, longer but at present stronger from the viewpoint of quantitative
estimates on derivatives of the Kolmogorov semi-group. In Flandoli \cite{Fla}
one can see a variation on the proof of \cite{DPDe} but applied to the Markov
selections constructed here.

Some further consequences of the theory developed here, in particular the
equivalence of transition probabilities of all Markov selections can be found
in Flandoli \& Romito \cite{FlRo3}. Other equations with lack of well-posedness
may be approached by a variant of the method developed here, see Bl\"omker,
Flandoli \& Romito \cite{BFR} for a model of surface growth.
Finally, a preliminary version of the results presented here have been given
in Flandoli \& Romito \cite{FlRo}.
\subsection*{Layout of the paper}
The paper is organised in two main parts. In the first part (from Section
\ref{s:abstract} to Section \ref{s:consequences}) we explain the main results
and we give their proofs. In the second part, constituted by the appendixes,
we give proofs of some complementary facts that have been used in the first
part. The reason behind this unusual layout is that the many technical
proofs contained in the appendixes could obfuscate the main ideas we wish
to explain.

The precise content of the paper is the following. In Section \ref{s:abstract}
we prove the abstract selection principle for Markov processes. Section
\ref{s:martingale} introduces the solutions to the martingale problem
for the NSE. In Section \ref{s:markov} we apply the
abstract selection principle to the solutions of the martingale problem.
Some of the proofs that are needed in this section are postponed to
Appendixes \ref{s:exist} and \ref{s:tecnic}. Finally, in Sections \ref{s:regular}
and \ref{s:consequences}, we prove that any Markov selection is regular, in the
sense explained above, and some consequences of this result.

Appendix \ref{s:exist} contains an existence result for those solutions
defined in Section \ref{s:martingale}. A slight modification of this
existence proof is used also in Section \ref{s:markov}. Appendix \ref{s:tecnic}
contains some technical results on the special super-martingales defined in this
paper. In Appendix \ref{s:related} we prove some useful facts on 
an auxiliary equation that we need to handle noise roughness and an equation
with truncated non-linearity. Both equations are obtained as modification of
the original NSE. Finally, Appendix \ref{s:nonlinear} contains
two different estimates of the Navier-Stokes non-linear term. We point out that,
according to our knowledge, the inequality of Lemma \ref{l:Breg} is new
for some values of the parameter.

\subsection*{Acknowledgements}
The authors wish to thank the referee of the paper for the several suggestion
that dramatically helped in improving the paper and for pointing out an error
in an earlier version. The authors wish to thank also B. Goldys
for letting them know the interesting reference \cite{StYo}.
\section{Pre-Markov families of probability measures}\label{s:abstract}
In this section we formalise the properties that a set-valued map must
have, in order to ensure the existence of a Markov selection. The content
of the following sections is an extension of the theory presented in
Chapter $12$ of Stroock \& Varadhan \cite{StVa}, with some changes, in
order to take into account the infinite dimensional setting of stochastic
partial differential equations. In particular, trajectories can have
different regularity properties in different spaces. The main novelty is
that we need to introduce the concept of \emph{almost sure markovianity},
in order to handle the problem that the energy inequality does not hold
for all times, which ultimately is related to the lack of continuity of
the trajectories in more regular spaces.

Most of the proofs of this section follow closely those of Stroock \&
Varadhan \cite{StVa}, with obvious differences whenever the extensions
stated above apply.
\subsection{Preliminaries}\label{ss:prelim}
We start by giving a few definitions and notations. Let
$$
\V\subset\H\subset\V'
$$
be a Gelfand triple of separable Hilbert spaces with continuous injections.
Set
$$
\Om=C([0,\infty);\V'),
$$
denote by $\B$ the Borel $\s$-field of $\Om$ and by $\Pr(\Om)$ the set
of all probability measures on $(\Om,\B)$. Define the canonical process
$\xi:\Om\to\V'$ as
$$
\xi_t(\om)=\om(t).
$$
\subsubsection{Preliminaries on the state space}
Define, for each $t\ge0$, the $\s$-field $\B_t=\s[\xi_s:\ 0\le s\le t]$.
Notice that one can identify this $\s$-field with the Borel $\s$-field
of $\Om_t=C([0,t];\V')$, since the set $\Om_t$ can be seen as a Borel
subset of $\Om$. Similarly, one can define $\Om^t=C([t,\infty);\V')$ and
$\B^t=\s[\xi_s:\ s\ge t]$.
Finally, define for each given $t>0$, the map $\Phi_t:\Om\to\Om^t$ as
$$
\Phi_t(\om)(s)=\om(s-t)
\qquad s\ge t.
$$
The next lemma shows that some sets, that we shall use in the sequel,
are indeed Borel sets.
Such results are well-known in a general
framework and are stated here for the sake of completeness.
\begin{lemma}\label{l:weakcont}
The set $L^\infty_\loc([0,\infty);\H)\cap\Om$ is a Borel set in $\Om$.
Moreover,
$$
L^\infty_\loc([0,\infty);\H)\cap\Om = C([0,\infty);\H_\s)\cap\Om,
$$
where $\H_\s$ denotes the space $\H$ endowed with the weak
topology. Finally, the set $L^2_\loc([0,\infty);\V)\cap\Om$
is Borel in $\Om$ as well.
\end{lemma}
\begin{proof}
We start by proving the equality. First, we easily have, by standard
arguments, that $C([0,\infty);\H_\s)\cap\Om$ is in
$L^\infty_\loc([0,\infty);\H)\cap\Om$. The other inclusion follows from
Lemma $1.4$, \S$3$ of Temam \cite{Tem}.

In order to prove measurability, notice that the map
$(t,\om)\mapsto\|\om(t)\|_\H\in[0,\infty]$ (where
the map takes the value $+\infty$ whenever $\om(t)\not\in\H$),
with $t\in[0,\infty)$ and $\om\in\Om$, is lower semi-continuous,
hence it is measurable. Let $D\subset[0,+\infty)$ be a countable
dense set. Observe that, by semi-continuity, for each $T>0$
$\sup_{t\in[0,T]}\|\om(t)\|_\H=\sup_{t\in D\cap[0,T]}\|\om(t)\|_\H$,
and so
$$
L^\infty_\loc([0,\infty);\H)\cap\Om
=\bigcap_{T=1}^\infty\bigcup_{R=1}^\infty\bigcap_{t\in D\cap[0,T]}
    \{\om\in\Om\ |\ \|\om(t)\|_\H\le R\}
$$
is measurable. Similarly, $(t,\om)\to\|\om(t)\|_\V$ is also
lower semi-continuous, so for each $T>0$ the map
$\om\to\int_0^T\|\om(t)\|_\V^2\,dt$ is lower semi-continuous as well.
Hence,
$$
L^2_\loc([0,\infty);\V)\cap\Om
=\bigcap_{T=1}^\infty\bigcup_{R=1}^\infty\bigl\{\om\in\Om\ |\ \int_0^T\|\om(t)\|_\V\,dt\le R\bigr\}
$$
is measurable.
\end{proof}
A straightforward consequence of the above result is given in the
following lemma.
\begin{lemma}
Let $P\in\Pr(\Om)$ be such that
$$
P[C([0,\infty);\H_\s)\cap\Om]=1.
$$
Then, for any given $t\ge0$, the mapping $\om\mapsto\om(t)$ has a
$P$-modification on $\B_t$ which is $\B_t$-measurable with values
in $(\H,\B(\H))$, where $\B(\H)$ is the Borel $\s$-field of $\H$.
\end{lemma}
\subsubsection{Preliminaries on disintegration and reconstruction
of probabilities}\label{sss:prelim}
Prior to the analysis of the Markov property in its different flavours,
we need some additional definitions and notations.

Given $P\in\Pr(\Om)$ and $t>0$, we will denote by
$\om\mapsto P|_{\B_t}^\om:\Om\to\Pr(\Om^t)$ a regular conditional
probability distribution of $P$ on $\B_t$. Since $\Om$ is a \emph{Polish}
space and every $\s$-field $\B_t$ is finitely generated, such a function
exists and is unique, up to $P$-null sets. In particular,
$$
P|_{\B_t}^\om[\xi_t=\om(t)]=1
$$
for all $\om\in\Om$, and, if $A\in\B_t$ and $B\in\B^t$,
$$
P(A\cap B)=\int_A P|_{\B_t}^\om(B)\,P(d\om).
$$
As conditional probabilities correspond to \emph{disintegration} with
respect to a $\s$-field, we define below the \emph{reconstruction},
which is, in a way, a sort of inverse procedure to disintegration.
\begin{definition}\label{d:raccordo}
Consider a probability $P\in\Pr(\Om)$, a time $t>0$ and a $\B_t$-measurable
map $Q:\Om\to\Pr(\Om^t)$ such that
$$
Q_\om[\xi_t=\om(t)]=1,
\qquad\text{for all }\om\in\Om.
$$
Then denote by $P\otimes_t Q$ the unique probability measure on $\Om$ such
that
\begin{enumerate}
\item $P\otimes_t Q$ and $P$ agree on $\B_t$,
\item $(Q_\om)_{\om\in\Om}$ is a regular conditional probability
      distribution of $P\otimes_t Q$ on $\B_t$.
\end{enumerate}
\end{definition}
Details on the measure whose existence is claimed above can be found in
Lemma $6.1.1$ and Theorem $6.1.2$ of Stroock \& Varadhan \cite{StVa}.
\subsection{The Markov property}
Given a family $(P_x)_{x\in\H}$ of probability measures,
the \emph{Markov property} can be stated as
$$
P_x|_{\B_t}^\om=\Phi_t P_{\om(t)},
\qquad\text{for }P_x-\text{a.e.\ }\om\in\Om.
$$
for each $x\in\H$ and for all $t\ge0$.
In view of application of the results of this section to the Navier-Stokes
equation, this definition is too strong. We give a slightly weaker definition,
where conditions on time are relaxed and a few exceptional time instants
are allowed.
\begin{definition}[almost sure Markov property]\label{d:asmarkov}
Let $x\mapsto P_x$ be a measurable map defined on $\H$ with values
in $\Pr(\Om)$ such that $P_x[C([0,\infty);\H_\s)\cap\Om]=1$ for
all $x\in\H$.

The family $(P_x)_{x\in\H}$ has the \emph{almost sure Markov property}
if for each $x\in\H$ there is a set $T\subset(0,\infty)$ with null
Lebesgue measure, such that
$$
P_x|_{\B_t}^\om=\Phi_tP_{\om(t)}
\qquad\text{for }P_x-\text{a.e.\ }\om\in\Om,
$$
for all $t\not\in T$.
\end{definition}
Following Stroock \& Varadhan \cite{StVa}, we introduce an analogous
definition for set-valued maps of probability measures. Such multi-valued
maps shall satisfy the a.\ s.\ Markov property in a suitable way.
Denote by $\comp(\Pr(\Om))$ the family of all compact subsets of
$\Pr(\Om)$.
\begin{definition}[almost sure pre-Markov family]\label{d:aspremarkov}
Let a measurable map $\ms:\H\to\comp(\Pr(\Om))$ be given
such that $P[C([0,\infty);\H_\s)\cap\Om]=1$ for all $x\in\H$ and
$P\in\ms(x)$.

The family $(\ms(x))_{x\in\H}$  is \emph{almost surely pre-Markov} if
for each $x\in\H$ and $P\in\ms(x)$ there is a set $T\subset(0,\infty)$
with null Lebesgue measure, such that for all $t\not\in T$ the following
properties hold:
\begin{enumerate}
\item\emph{(disintegration)} there exists $N\in\B_t$ with $P(N)=0$ such
that for all $\om\not\in N$
\begin{center}
$\om(t)\in\H$\qquad and \qquad $P|_{\B_t}^\om\in\Phi_t\ms(\om(t))$;
\end{center}
\item\emph{(reconstruction)} for each $\B_t$-measurable map
$\om\mapsto Q_\om:\Om\to\Pr(\Om^t)$ such that there is $N\in\B_t$ with
$P(N)=0$ and for all $\om\not\in N$,
\begin{center}
$\om(t)\in\H$\qquad and \qquad $Q_\om\in\Phi_t\ms(\om(t))$;
\end{center}
then $P\otimes_t Q\in\ms(x)$.
\end{enumerate}
\end{definition}
\begin{remark}
We aim to make clear the meaning of measurability for a
$\comp(\Pr(\Om))$-valued map. The set $\Pr(\Om)$, with the weak convergence
of measures, is a Polish space, hence the set $\comp(\Pr(\Om))$, endowed with
the \emph{Hausdorff metric}, is a metric space. In this context, measurability
of the above map refers to Borel measurability with respect to the Hausdorff
metric.
\end{remark}
\begin{remark}\label{r:singleton}
If every $\ms(x)$ is a singleton, the a.\ s.\ pre-Markov family is indeed an
a.\ s.\ Markov family of probability measures, as stated in Definition
\ref{d:asmarkov}. Notice that, in the framework of NSE
we shall examine in the next sections, each $\ms(x)$ represents the set of
all solutions starting at $x$ and well-posedness of the martingale problem
follows if at least for one $x$, the set $\ms(x)$ contains a single point
(see Corollary \ref{c:lawuniq}).

In the general setting one can refer to Theorem $12.2.4$ of Stroock
\& Varadhan \cite{StVa}. Moreover Stroock \& Yor \cite{StYo} provide
several examples of non-uniqueness and comments on Markov selections.
\end{remark}
\subsection{Existence of Markov selections}\label{ss:marsel}
In this section we give the main abstract result concerning
the existence of Markov selections.
\begin{theorem}\label{t:mainabstract}
Let $(\ms(x))_{x\in\H}$ be an a.\ s.\ pre-Markov family with non-empty
convex values. Then there is a measurable map $x\mapsto P_x$ on $\H$
with values in $\Pr(\Om)$ such that $P_x\in\ms(x)$ for all $x\in\H$
and $(P_x)_{x\in\H}$ has the a.\ s.\ Markov property (as defined in
Definition \ref{d:asmarkov}).
\end{theorem}
\begin{remark}
In view of a possible extension of the previous theorem
to Feller selections (that is, Markov selections with some
kind of continuous dependence with respect to the initial
condition), we remark that so far, we cannot expect any more
regularity with respect to the initial condition from
a Markov selection, provided by the abstract principle, than
measurability. In general, continuity with respect to initial
condition cannot be gained in this abstract setting, as one can
see from the following classical example,
$$
\dot x=\sgn(x)\arctan\sqrt{|x|}.
$$
Indeed, any selection from solutions to the above equation
is not continuous at $t=0$. Hence, continuity needs a further
analysis and we will see that the noise plays a major role.
\end{remark}
Prior to the proof of the theorem above (which is postponed to
page \pageref{pr:mainabstract}), we need to give some definitions
and state some useful results.

In order to identify a unique representative in each class $\ms(x)$,
we define below a method to \emph{reduce} such classes, by means of
maximisations. To this aim, define for each $f\in C_b(\V',\R)$ and
each $\la>0$, the map $J_{\la,f}:\Pr(\Om)\to\R$ as
$$
J_{\la,f}(P)=\E^P\bigl[\int_0^\infty\e^{-\la t}f(\xi_t)\,dt\bigr],
\qquad P\in\Pr(\Om),
$$
and for each $x\in\H$,
$$
\rv_\la f(x)=\sup_{P\in\ms(x)}J_{\la,f}(P)
$$
and
$$
\ms_{\la,f}(x)=\set{P\in\ms(x)}{J_{\la,f}(P)=\rv_\la f(x)}.
$$
\begin{lemma}\label{l:max}
Let $(\ms(x))_{x\in\H}$ be an a.\ s.\ pre-Markov family with
non-empty convex values and let $\la>0$ and $f\in C_b(\V',\R)$.
Then $\rv_\la f(x)$ is well-defined and $(\ms_{\la,f}(x))_{x\in\H}$
is again an a.\ s.\ pre-Markov family with non-empty convex values.
\end{lemma}
\begin{proof}
Since each map $J_{\la,f}$ is continuous and $\ms(x)$ is compact,
$\rv_\la f(x)$ is well-defined and $\ms_{\la,f}(x)$ is non-empty.

We prove that $(\ms_{\la,f}(x))_{x\in\H}$ is an a.\ s.\ pre-Markov
family. First, each set $\ms_{\la,f}(x)$ is compact and
$x\mapsto\ms_{\la,f}(x)$ is measurable, by virtue of Lemma $12.1.7$
of Stroock \& Varadhan \cite{StVa}. Moreover, convexity follows since
$J_{\la,f}$ is linear in $P$ and
$$
P[C([0,\infty);\H_\s)\cap\Om]=1,
$$
for each $P\in\ms_{\la,f}(x)$, since $\ms_{\la,f}(x)\subset\ms(x)$.

Prior to the proof of the properties of \emph{disintegration} and
\emph{reconstruction}, we notice that, if $P\in Pr(\Om)$,
\begin{equation}\label{e:useful}
\E^{\Phi_tP}\int_t^\infty\e^{-\la s}f(\xi_s)\,ds
=\e^{-\la t}J_{\la,f}(P).
\end{equation}
Fix $x\in\H$ and $P\in\ms(x)$, let $T\subset(0,\infty)$ be the
null measure set corresponding to $P$ whose existence follows from the
a.\ s.\ pre-Markovianity of $(C(x))_{x\in\H}$, and fix $t\not\in T$.
We prove first the \emph{disintegration} property. Consider a regular
conditional probability distribution $P|_{\B_t}^\cdot$ of $P$ on
$\B_t$ and set
\begin{align*}
N&=\set{\om}{\om(t)\not\in\H\text{ or }P|_{\B_t}^\om\not\in\Phi_t\ms(\om(t))},\\
N'&=\set{\om\not\in N}{P|_{\B_t}^\om\not\in\ms_{\la,f}(\om(t))}.
\end{align*}
By the \emph{disintegration} property for $(\ms(x))_{x\in\H}$, we have that
$N\in\B_t$ and $P(N)=0$. Moreover, by Lemma $12.1.9$ of Stroock
\& Varadhan \cite{StVa}, $N'\in\B_t$. The property is true if
$P(N')=0$. To this end, let $x\mapsto R_x\in\ms_{\la,f}(x)$
be a measurable map (the existence of this map is ensured by a standard
measurable selection theorem, see e.\ g.\ Lemma $12.1.10$ of
Stroock \& Varadhan \cite{StVa}) and notice that
$\om\to\Phi_t R_{\om(t)}$ is again measurable. Define
$$
Q_\om=\begin{cases}
  P|_{\B_t}^\om     &\qquad\om\not\in N\cup N',\\
  \Phi_t R_{\om(t)} &\qquad\om\in N',\\
  \delta_0          &\qquad\om\in N,
\end{cases}
$$
where $\delta_0$ is the Dirac measure on the trajectory which is identically
zero and it is used as the fallback measure. By the \emph{reconstruction} property
for the family $(\ms(x))_{x\in\H}$, it follows that $P\otimes_t Q\in\ms(x)$ and so,
\begin{align*}
J_{\la,f}(P)
&\ge    J_{\la,f}(P\otimes_t Q)\\
&=      \E^P\int_0^t\e^{-\la s}f(\xi_s)\,ds
       +\E^P[\E^{Q_\om}\int_t^\infty\e^{-\la s}f(\xi_s)\,ds]\\
&=      J_{\la,f}(P)
       +\E^P\big[\chi_{N'}\E^{\Phi_tR_{\om(t)}}[\int_t^\infty\e^{-\la s}f(\xi_s)\,ds]\big]\\
&\quad -\E^P\big[\chi_{N'}\E^{P|_{\B_t}^\om}[\int_t^\infty\e^{-\la s}f(\xi_s)\,ds]\big]\\
&=      J_{\la,f}(P)+\e^{-\la t}\E^P[\chi_{N'}\big(\rv_\la f(\om(t))-J_{\la,f}(\Phi_t^{-1}P|_{\B_t}^\om)\big)]
\end{align*}
since $R_{\om(t)}\in\ms_{\la,f}(\om(t))$ and $P|_{\B_t}^\om\in\Phi_t\ms(\om(t))$
for $\om\in N'$, using also \eqref{e:useful}. Hence,
$\E^P[\chi_{N'}\big(\rv_\la f(\om(t))-J_{\la,f}(\Phi_t^{-1}P|_{\B_t}^\om))\big)]\le0$. On the
other hand, for each $\om\in N'$, $\Phi_t^{-1}P|_{\B_t}^\om\not\in\ms_{\la,f}(\om(t))$, so that
$J_{\la,f}(\Phi_t^{-1}P|_{\B_t}^\om)<\rv_\la f(\om(t))$ and in conclusion $P(N')=0$.

We prove the \emph{reconstruction} property. Let $\om\mapsto Q_\om$ be a
$\B_t$-measurable from $\Om$ to $\Pr(\Om^t)$. Let $N\in\B_t$ be such
that $P(N)=0$ and for all $\om\not\in N$, $\om(t)\in\H$ and
$Q_\om\in\Phi_t\ms_{\la,f}(\om(t))$. By the \emph{reconstruction} property
for $(\ms(x))_{x\in\H}$, it follows that $P\otimes_t Q\in\ms(x)$. Moreover,
by using \eqref{e:useful},
\begin{align*}
J_{\la,f}(P\otimes_t Q)
&=    \E^P\int_0^t\e^{-\la s}f(\xi_s)\,ds
     +\E^P[\E^{Q_\om}\int_t^\infty\e^{-\la s}f(\xi_s)\,ds]\\
&=    \E^P\int_0^t\e^{-\la s}f(\xi_s)\,ds
     +\e^{-\la t}\E^P[J_{\la,f}(\Phi_t^{-1}Q_\om)]\\
&\ge  \E^P\int_0^t\e^{-\la s}f(\xi_s)\,ds
     +\e^{-\la t}\E^P[J_{\la,f}(\Phi_t^{-1}P|_{\B_t}^\om)]\\
&=    \E^P\int_0^t\e^{-\la s}f(\xi_s)\,ds
     +\E^P[\E^{P|_{\B_t}^\om}\int_t^\infty\e^{-\la s}f(\xi_s)\,ds]\\
&=    J_{\la,f}(P),
\end{align*}
since $P|_{\B_t}^\om\in\Phi_t\ms(\om(t))$. In conclusion
$P\otimes_t Q\in\ms_{\la,f}(x)$.
\end{proof}
It is now possible to prove the main theorem of this section.
\begin{proof}[Proof of Theorem \ref{t:mainabstract}]
Let $(\s_n)_{n\in\N}$ be a dense set in $[0,\infty)$ and\label{pr:mainabstract}
$(\phi_n)_{n\in\N}$ be a dense set in $C_b(\V',\R)$, and consider an
enumeration $(\la_n,f_n)_{n\ge1}$ of $(\s_m,\phi_n)_{n,m\in\N}$.
For each $x\in\H$ define inductively
$$
\ms^{(0)}(x)=\ms(x),
\qquad
\ms^{(n+1)}=\ms^{(n)}_{\la_n,f_n}(x),\quad n\ge1,
$$
where $\ms^{(n)}_{\la_n,f_n}(x)$ is obtained from $\ms^{(n)}(x)$,
as in Lemma \ref{l:max}, to be the set of measures where the
maximum value of $J_{\la_n,f_n}$ is achieved.
Finally set $\ms^{(\infty)}(x)=\bigcap\ms^{(n)}(x)$. By Lemma
\ref{l:max}, each family $(\ms^{(n)}(x))_{x\in\H}$ is a.\ s.\ 
pre-Markov, and it is easy to check that
$(\ms^{(\infty)}(x))_{x\in\H}$ is a.\ s.\ pre-Markov as
well. The proof is complete if we show that each set
has exactly one element.

Let $P$, $Q\in\ms^{(\infty)}(x)$, then
$J_{\la_n,f_n}(P)=J_{\la_n,f_n}(Q)$ for all $n\in\N$,
that is, for all $m$, $n\in\N$,
$$
 \int_0^\infty\e^{-t\s_m}\E^P[\phi_n(\xi_t)]\,dt
=\int_0^\infty\e^{-t\s_m}\E^Q[\phi_n(\xi_t)]\,dt,
$$
so by the uniqueness of the Laplace transform and by the density
assumptions, we can deduce that $\E^P[g(\xi_t)]=\E^Q[g(\xi_t)]$,
$t\ge0$, holds for each bounded measurable function $g:\V'\to\R$.

In order to prove that $P=Q$, we need to show that all finite marginals
coincide, i.\ e.\ that
\begin{equation}\label{e:marginal}
 \E^P[g_1(\xi_{t_1})\ldots g_n(\xi_{t_n})]
=\E^Q[g_1(\xi_{t_1})\ldots g_n(\xi_{t_n})]
\end{equation}
for all bounded measurable $g_1,\ldots,g_n$ and all $0\le
t_1<\ldots<t_n$. Since we are going to use the a.\ s.\ Markov
property, we can prove \eqref{e:marginal} only for a set of times
which has full measure with respect to the Lebesgue measure. Namely,
if $T_P$ and $T_Q$ are the zero measure subsets of $(0,+\infty)$
corresponding respectively to $P$ and $Q$ (as given by Definition
\ref{d:aspremarkov}), then \eqref{e:marginal} will hold only for
times not belonging to $T_P\cup T_Q$. Since functions
$(t_1,\dots,t_n)\to\E^P[g_1(\xi_{t_1})\ldots g_n(\xi_{t_n})]$ are
continuous, \eqref{e:marginal} extends to all times. So, it is
sufficient to prove \eqref{e:marginal} only for times in $T_P\cup T_Q$.

We proceed by induction. For $n=1$, the claim is true
by virtue of the above considerations. Assume that the claim is true for
a integer $n$. Let $g_1,\ldots,g_{n+1}:\V'\to\R$ be bounded measurable
functions and $0\le t_1<\ldots<t_{n+1}$ be times in $(T_P\cup T_Q)^c$.

Let $\F_{t_1\ldots t_n}$ be the $\s$-field generated by $\xi_{t_1},
\ldots,\xi_{t_n}$, then
$$
 \E^P[g_1(\xi_{t_1})\ldots g_{n+1}(\xi_{t_{n+1}})]
=\E^P\bigl[g_1(\xi_{t_1})\ldots g_n(\xi_{t_n})\E^P[g_{n+1}(\xi_{t_{n+1}})|\F_{t_1\ldots t_n}]\bigr].
$$
So, it is sufficient to show that
$$
 \E^P[g_{n+1}(\xi_{t_{n+1}})|\F_{t_1\ldots t_n}]
=\E^Q[g_{n+1}(\xi_{t_{n+1}})|\F_{t_1\ldots t_n}]
\qquad P-a.\ s.,
$$
or, in different terms, that for every $\varphi\in C_b(\V',\R)$,
$$
 \E^{P|_{\F_{t_1\ldots t_n}}^\om}[\varphi(\xi_{t_{n+1}})]
=\E^{Q|_{\F_{t_1\ldots t_n}}^\om}[\varphi(\xi_{t_{n+1}})],
\qquad P-a.e.\ \om.
$$
If we had $\B_{t_n}$ in place of $\F_{t_1\ldots t_n}$, the above
equality would be just a consequence of the \emph{disintegration}
property. We are going to use this fact. By the \emph{disintegration}
property, there is $N_P\in\B_{t_n}$ such that $P(N_P)=0$, $\om(t_n)\in\H$
and $P|_{\B_{t_n}}^\om\in\Phi_{t_n}\ms^{(\infty)}(\om(t_n))$
for all $\om\not\in N_P$. Moreover, there is $A_P\in\F_{t_1\ldots t_n}$
such that $P(A_P)=0$ and for all $\om\not\in A_P$,
\begin{equation}\label{e:daBaF}
 P|_{\F_{t_1,\ldots t_n}}^\om
=\int P|_{\B_{t_n}}^{\om'}P|_{\F_{t_1,\ldots t_n}}^\om(d\om').
\end{equation}
Since $P(N_P)=0$, there is also a set $B_P\in\F_{t_1\ldots t_n}$ such that
$P|_{\F_{t_1\ldots t_n}}^\om(N_P)=0$ for $\om\not\in B_P$.

Now, set $N_P'=A_P\cup B_P\in\F_{t_1\ldots t_n}$, then by using
\eqref{e:daBaF} and by the convexity of sets $\Phi_{t_n}\ms^{(\infty)}(\om(t_n))$,
it follows that
$P|_{\F_{t_1\ldots t_n}}^\om\in\Phi_{t_n}\ms^{(\infty)}(\om(t_n))$
for $\om\not\in N_P'$. Similarly, one can find $N_Q'\in\F_{t_1\ldots t_n}$
such that $Q(N_Q')=0$ and
$Q|_{\F_{t_1\ldots t_n}}^\om\in\Phi_{t_n}\ms^{(\infty)}(\om(t_n))$
for $\om\not\in N_Q'$.

By the induction hypothesis, $P$ and $Q$ agree on $\F_{t_1\ldots t_n}$,
so that, if $\widetilde N=N_P'\cup N_Q'$, then $P(\widetilde N)=Q(\widetilde N)=0$
and
$$
 \E^{P|_{\F_{t_1\ldots t_n}}^\om}[\varphi(\xi_{t_{n+1}})]
=\E^{Q|_{\F_{t_1\ldots t_n}}^\om}[\varphi(\xi_{t_{n+1}})],
$$
for all $\om\not\in\widetilde N$.
\end{proof}
\begin{remark}
It is worth noticing that all results of this section hold if the
definitions of almost sure Markov property and almost sure
pre-Markov family are replaced by analogous definitions, where
each almost sure property is indeed sure. More precisely, it
is sufficient that in Definition \ref{d:aspremarkov} each set
$T$ of exceptional times is empty. Call it a pre-Markov family.
Then the theorem below holds and its proof is entirely similar
to that of Theorem \ref{t:mainabstract}.
\end{remark}
\begin{theorem}
Let $(\ms(x))_{x\in\H}$ be a pre-Markov family with non-empty
convex values. Then there is a measurable map $x\mapsto P_x$ defined
on $\H$ with values in $\Pr(\Om)$ such that $P_x\in\ms(x)$ for
all $x\in\H$ and $(P_x)_{x\in\H}$ has the Markov property.
\end{theorem}
\begin{remark}
In case of pre-Markov family (that is, without the almost sure),
one can even show, as in Stroock \& Varadhan \cite{StVa}, that
there are strong Markov selections. For almost sure families,
the strong Markov property seems to be technically more complicated.
\end{remark}
\section{The martingale problem for the Navier-Stokes equations}\label{s:martingale}
Let $\T=[0,1]^3$ be the three-dimensional torus and let $\test$ be the space
of infinitely differentiable divergence-free periodic vector fields on $\R^3$
with zero mean. Let $H$ be the closure of $\test$ in the norm of $L^2(\T,\R^3)$,
$V$ the closure of $\test$ in the norm of $H^1(\T,\R^3)$, and let
$D(A)=\{\,u\in H\,|\,\Delta u\in H\,\}$. Let $A:D(A)\to H$ be the Stokes
operator,
$$
Au=-\Delta u, \qquad u\in D(A),
$$
it is a positive linear self-adjoint operator on $H$ and we can define the
powers $A^\al$, $\al\in\R$, with domain $D(A^\al)$. By proper identifications
of dual spaces, in particular we have
$$
V\subset H\subset V'\subset D(A)'.
$$
The bi-linear operator $B:V\times V\to V'$ is defined as
$$
B(u,v)=\mathds{P}_{\!\!\!\text{\Tiny div}}(u\cdot\nabla)v,
$$
where $\mathds{P}_{\!\!\!\text{\Tiny div}}$ is the projection onto divergence-free
vector fields (for more details on the above definitions, a standard reference
is Temam \cite{Tem2}).
\begin{assumption}[Noise is trace class]\label{a:noiseass1}
The covariance $\Q:H\to H$ of the noise driving the equation
\eqref{e:nseq} is a symmetric non-negative trace-class operator on $H$.
\end{assumption}
Let $(e_i)_{i\in\N}$ be a complete orthonormal system of eigenvectors
and denote by $(\s_i^2)_{i\in\N}$ the eigenvalues of $\Q$. By the above
assumption, the following quantity is finite,
$$
\s^2=\sum_{i=1}^\infty\s_i^2.
$$
Further assumptions on $\Q$ will be given later in Section \ref{s:regular}
(see Assumption \ref{a:noiseass2}) to ensure the validity of results contained
in that section.

We consider the three-dimensional Navier-Stokes equations in its abstract
form,
\begin{equation}\label{e:nseq}
du+\bigl(\nu Au+B(u,u)\bigr)\,dt=\Q^\frac12\,dW
\end{equation}
where $W$ is a cylindrical Wiener process on $H$ (see Da Prato \& Zabczyk
\cite{DPZa2} for more details).

In the rest of the paper we shall assume, for simplicity, that $\nu=1$,
since the size of viscosity plays no essential role in this paper.
\subsection{Almost sure super-martingales}
Prior to the definition of solutions to the martingale problem associated
to the Navier-Stokes equations \eqref{e:nseq}, we need to introduce a
slightly different variant of definition of super-martingale.
\begin{definition}[a.\ s.\ super-martingale]\label{d:asSM}
An adapted process $(\theta_t,\B_t,P)_{t\ge0}$ is an \emph{a.\ s.\ super-martingale}
if $\E^P|\theta_t|<\infty$ for all $t\ge0$ and there is a Lebesgue
measurable set $T_\theta\subset(0,\infty)$, with null Lebesgue measure, such that
$$
\E^P[\theta_t\uno_A]\le\E^P[\theta_s\uno_A]
$$
holds for every $s\not\in T_\theta$, every $t\ge s$ and every $A\in\B_s$.

The set $T_\theta$ will be called the set of \emph{exceptional times}
of $\theta$.
\end{definition}
We shall see in Appendix \ref{s:tecnic} that some of the results
which are true for super-martingales, hold for a special class of
almost sure super-martingales.
\subsection{The solutions to the martingale problem}
In view of the results of previous section, we consider the particular
case where $\V=V$, $\H=H$ and $\V'=D(A)'$. We set
$$
\Omns = C([0,\infty);D(A)').
$$
This space will play the role of state space for the solutions to
\eqref{e:nseq}. As in Section \ref{ss:prelim}, we denote by $\Bns$ the
$\s$-field of Borel sets of $\Omns$, and, for each $t\ge0$,
by $\BNS_t=\s(\om|_{[0,t]}\,:\,\om\in\Omns)$ and
$\Bns^t=\s(\om|_{[t,\infty)}\,:\,\om\in\Omns)$ the $\s$-fields of
past and future, with respect to time $t$, events.

Next, we give the definition of solution to the martingale problem
associated to the Navier-Stokes equations \eqref{e:nseq} that will
be considered in the paper. As we shall see, the definition incorporates,
in a peculiar new form, the energy inequality. As in the deterministic case,
the energy estimate cannot be deduced directly from the equation, it can
be proved only for those solution suitably obtained by a regularisation
procedure (see Appendix \ref{s:exist} for an example of approximation).
\begin{definition}\label{d:asNSsol}
Given $\mu_0\in\Pr(H)$, a probability $P$ on $(\Omns,\Bns)$ is a solution
starting at $\mu_0$ to the martingale problem associated to the Navier-Stokes
equations \eqref{e:nseq} if
\begin{itemize}
\item[\map{1}] $P[L_\loc^\infty([0,\infty);H)\cap L^2_\loc([0,\infty);V)]=1$;
\item[\map{2}] for each $\varphi\in\test$ the process $M_t^\varphi$, defined
$P$--a.\ s.\ on $(\Omns,\Bns)$ as
$$
M_t^\varphi
= \ps{\xi_t-\xi_0}{\varphi}
 +\nu\int_0^t\ps{\xi_s}{A\varphi}\,ds
 -\int_0^t\ps{B(\xi_s,\varphi)}{\xi_s}\,ds
$$
is square integrable and $(M_t^\varphi,\BNS_t,P)$ is a continuous martingale with
quadratic variation
$$
[M^\varphi]_t=t|\Q^\frac12\varphi|^2_H;
$$
\item[\map{3}] the process $E^1_t$, defined $P$--a.\ s.\ on $(\Omns,\Bns)$ as
$$
E_t^1
= |\xi_t|^2_H
 +2\nu\int_0^t|\xi_s|_V^2\,ds
 -|\xi_0|^2_H
 -t\s^2
$$
is $P$-integrable and $(E_t^1,\BNS_t,P)$ is an a.\ s.\ super-martingale;
\item[\map{4}] for each $n\ge2$, the process $E^n_t$, defined $P$--a.\ s.\ on
$(\Omns,\Bns)$ as
$$
E_t^n
= |\xi_t|^{2n}_H
 +2n\nu\int_0^t|\xi_s|^{2n-2}_H|\xi_s|_V^2\,ds
 -|\xi_0|^{2n}_H
 -n(2n-1)\s^2\int_0^t|\xi_s|_H^{2n-2}\,ds
$$
is $P$-integrable and $(E_t^n,\BNS_t,P)$ is an a.\ s.\ super-martingale;
\item[\map{5}] $\mu_0$ is the marginal of $P$ at time $t=0$.
\end{itemize}
\end{definition}
\begin{remark}\label{r:exceptime}
Given a solution $P$ to the martingale problem associated to the Navier-Stokes
equation \eqref{e:nseq}, define the set $T_P\subset(0,\infty)$ of
\emph{exceptional times} of $P$ as the union of the sets of exceptional times
of all a.\ s.\ super-martingales $E^1_t$, $E_t^2$, \dots.
\end{remark}
\begin{remark}\label{r:dense}
Due to property \map{2}, \map{3}, \map{4}  of the above definition, property
\map{2} itself needs to be verified only on a countable subset of $\test$ which
is dense in $\test$ with respect to the norm of $D(A)$. Indeed, if
$\varphi_n\to\varphi$ in $D(A)$, then $M_t^{\varphi_n}(\om)\to M_t^\varphi(\om)$
for all $\om\in L_\loc^\infty([0,\infty);H)\cap L^2_\loc([0,\infty);V)$, and
$\E^P|M_t^{\varphi_n}-M_t^\varphi|^2\to0$.
\end{remark}
\begin{remark}
Notice that condition \map{3} is a restatement of the energy inequality for
the Navier-Stokes equations obtained (formally) by the It\^o formula, while
condition \map{4} states the energy inequality for higher moments of $|\xi|_H^2$.
As usually in the literature concerning the Navier-Stokes equations only the
energy balance for the second moment is required, we remark that the need of
condition \map{4} will be apparent in the following pages (see for example
Lemma \ref{l:exsuper}, where it ensures uniform integrability).
\end{remark}
The next theorem shows that, under natural assumptions on the initial
condition, there is at least one solution, according to Definition \ref{d:asNSsol},
to the martingale problem associated to the Navier-Stokes equations \eqref{e:nseq}.

The proof of this theorem in our framework will turn out to be just a by-product
of slightly more general results presented in the next section which are needed
for the proof of Theorem \ref{t:mainmarkov}. Such results are postponed to Appendix
\ref{ss:exist} and suitably stated in order to be used for the proof of both the next
theorem and Theorem \ref{t:mainmarkov} below.
\begin{theorem}\label{t:mainex}
Let $\mu_0\in\Pr(H)$ be a probability measure such that
$$
\E^{\mu_0}[|x|_H^{2m}]<\infty,
\qquad m\ge1.
$$
Then, there exists a solution, starting at $\mu_0$, to the martingale problem
associated to the Navier-Stokes equations \eqref{e:nseq}, as defined in
Definition \ref{d:asNSsol}.
\end{theorem}
\begin{remark}\label{r:weak}
We finally remark that the definition of solution to the martingale
problem assumed in this paper is slightly different from those available
in the literature (see for example Flandoli \& Gatarek \cite{FlGa}).
A \emph{weak martingale solution} starting at $\mu_0\in\Pr(H)$ is a
filtered probability space $(\Sigma,\F,(\F_t)_{t\ge0},\abP)$,
a cylindrical Wiener process $(W_t)_{t\ge0}$ on $H$
and a continuous $D(A)'$-valued adapted process $u$ on $(\Sigma,\F,\abP)$ such
that $u$ is in spaces $L^\infty_\loc([0,\infty);H)$ and $L^2_\loc([0,\infty);V)$
$\abP$--a.\ s., $u(0)$ has law $\mu_0$ and the equations hold in distribution:
$$
\ps{u(t)-u(0)}{\varphi}
 +\int_0^t\ps{u}{A\varphi}\,ds
 -\int_0^t\ps{B(u,\varphi)}{u}\,ds
=\ps{\varphi}{\Q^\frac12 W_t}
$$
$\abP$--a.\ s.\ for all test function $\varphi$. It can be proved (see Flandoli
\cite{Fla2}) that a probability measure $P$ on $\Omns$ is the law of a weak
martingale solution if and only if properties \map{1}, \map{2} and \map{5} hold.
\end{remark}
\section{Markov selection for the 3D Navier-Stokes equations}\label{s:markov}
The section is devoted to the proof of existence of a Markov selection
for the solutions to the stochastic Navier-Stokes equations.
Define for each $x\in H$ the set $\ns(x)\subset\Pr(\Omns)$ as
\begin{equation}\label{e:NSfamily}
\set{P\in\Pr(\Omns)}
{P\text{ is a solution starting at $\delta_x$ of the martingale problem}}.
\end{equation}
\begin{theorem}\label{t:mainmarkov}
Assume Assumption \ref{a:noiseass1}. Then there exists a family $(P_x)_{x\in H}$
of solutions to the martingale problem associated to the Navier-Stokes equations
with the a.\ s.\ Markov property, as defined in Definition \ref{d:asmarkov}.
\end{theorem}
By Theorem \ref{t:mainabstract}, the proof of the above theorem amounts
to showing that the family defined above in \eqref{e:NSfamily} is a
\emph{a.\ s.\ pre-Markov family} (as defined in \ref{d:aspremarkov}).
The proof of this claim will be developed, for the sake of clarity,
in the following lemmas.
\begin{lemma}
For each $x\in H$, the set $\ns(x)$ is non-empty, convex and for all
$P\in\ns(x)$,
$$
P[C([0,\infty);H_\sigma)]=1.
$$
\end{lemma}
\begin{proof}
Given $x\in H$, a solution to the martingale problem exists due to
Theorem \ref{t:mainex}, so that each set $\ns(x)$ is non-empty.
Moreover, by property \map{1} and Lemma \ref{l:weakcont}, it
follows that $C([0,\infty);H_\s)$ is a $P$-full set.

Finally, it is easy to check that each $\ns(x)$ is convex, since
all properties in Definition \ref{d:asNSsol} involve integration
with respect to elements of $\ns(x)$.
\end{proof}
\begin{lemma}\label{l:compmeas}
For each $x\in H$, the set $\ns(x)$ is compact and the map
$\ns:H\to\comp(\Pr(\Omns))$ is Borel measurable.
\end{lemma}
\begin{proof}
Both properties stated in the lemma follow from the following
claim:
\begin{center}
\framebox{\begin{minipage}{0.8\linewidth}
for each sequence $(x_n)_{n\in\N}$ converging in
$H$ to $x$ and for each $P_n\in\ns(x_n)$, the sequence
$(P_n)_{n\in\N}$ has a limit point $P$ in $\ns(x)$, with
respect to weak convergence in $\Pr(\Omns)$.
\end{minipage}}
\end{center}
Indeed, for compactness, one takes $x_n=x$, while measurability
follows from Lemma $12.1.8$ of Stroock \& Varadhan \cite{StVa}.
In order to prove the claim, let $x_n\to x$ in $H$ and let
$P_n\in\ns(x_n)$. We first show that $(P_n)_{n\in\N}$ is
tight on $\Omns\cap L^2_\loc([0,\infty);H)$. By \map{3},
\map{4} with $n=2$ and Corollary \ref{c:doob}, we have that
for all $T>0$,
$$
\E^{P_n}\Bigl[\sup_{[0,T]}|\xi_t|^2_H+\int_0^T|\xi_s|^2_V\,ds\Bigr]\le C(\s,T,|x_n|_H).
$$
Next, let $(\Sigma,\F,(\F_t)_{t\ge0},{\abP})$ be a filtered
probability space, $(W_t)_{t\ge0}$ a cylindrical Wiener process
on $H$, $u$ be a process on $\Sigma$ whose law is
$P_n$ and such that $u$ is a weak martingale solution to
\eqref{e:nseq} (see Remark \ref{r:weak}). In particular,
$$
u(t)=x_n-\int_0^t(Au(s)+B(u(s)))\,ds+\Q^\frac12 W_t
\qquad {\abP}-a.\ s.,
$$
in $D(A)'$. Let $J_t(u)$ be the integral term on the right-hand
side of the above formula. By Burkholder, Davis \& Gundy
inequality, for all $\al\in(0,\frac12)$,
$p>1$ and $T>0$,
$$
\E^{\abP}\|\Q^\frac12 W\|_{W^{\al,p}(0,T;H)}\le C(\al,p,T,\s).
$$
Moreover, for each $\gamma\in(\frac32,2)$ we have
$|B(u)|^2_{D(A^{-\gamma})}\le|u|^2_H|u|^2_V$ and so for each
$T>0$,
\begin{align*}
\|J(u)\|^2_{W^{1,2}(0,T;D(A^{-\gamma}))}
&\le C(T)\int_0^T(|Au(s)|^2_{D(A^{-\gamma})}+|B(u(s))|^2_{D(A^{-\gamma})})\,ds\\
&\le C(T,\gamma)\int_0^T|u(s)|^2_V(1+|u(s)|^2_H)\,ds,
\end{align*}
which is bounded in expectation by \map{3} and \map{4}
(with $n=2$). In conclusion, tightness of $(P_n)_{n\in\N}$
follows from Lemma \ref{l:excomp} (there, we take $\xi_1=J(\xi)$).
Hence, there is a sub-sequence $(P'_n)_{n\in\N}$ converging weakly
in $\Omns\cap L^2_\loc([0,\infty);H)$ to some $P\in\Pr(\Om)$.

To conclude the proof, we have to show that $P\in\ns(x)$.
First, we notice that the marginals of $P_n'$ at time $0$ 
converge weakly to the marginal at time $0$ of $P$ and,
since such marginals converge to the Dirac measure in $x$,
\map{5} is true for $P$. Moreover, \map{1}, \map{3} and
\map{4} hold for $P$ by Lemma \ref{l:exsuper}. Finally,
\map{2} can be proved in the same way as in
Lemma \ref{l:exmart} for the Galerkin approximations.
\end{proof}
\begin{lemma}
The \emph{disintegration} property of Definition \ref{d:aspremarkov}
holds for the family $(\ns(x))_{x\in H}$.
\end{lemma}
\begin{proof}
Let $x\in H$ and $P\in\ns(x)$. Let $T_P$ be the set of exceptional
times of $P$ (see Remark \ref{r:exceptime}) and fix $t\not\in T_P$.
Let $\om\to P|_{\BNS_t}^\om$ be a regular conditional probability
distribution of $P$ on $\BNS_t$ (see Section \ref{sss:prelim}),
we aim to find a $P$-null set $N\in\BNS_t$ such that $\om(t)\in H$
and $P|_{\BNS_t}^\om\in\Phi_t\ns(x)$ for all $\om\not\in N$. We
shall have
$$
N=N_1\cup N_2\cup N_3\cup N_4\cup N_5,
$$
where all sets $N_1$, \dots, $N_5$ will be specified along the proof
and correspond respectively to properties \map{1}, \dots, \map{5}
of Definition \ref{d:asNSsol}.

Set
\begin{align}\label{e:sottospazi}
S_t&=\{\,\om\in\Omns\,:\,\om|_{[0,t]}\in L^\infty(0,t;H)\cap L^2(0,t;V)\,\},\\
S^t&=\{\,\om\in\Omns\,:\,\om|_{[t,\infty]}\in L^\infty_\loc([t,\infty),H)\cap L^2_\loc([t,\infty),V)\,\},\notag
\end{align}
and notice that $S_t\in\BNS_t$, $S^t\in\Bns^t$ and that $S_t\cap S^t$
is a $P$-full set by property \map{1}. Hence,
$$
1=P[S_t\cap S^t]=\int_{S_t}P|_{\BNS_t}^\om[S^t]\,P(d\om)
$$
and thus there is a $P$-null set $N_1\in\BNS_t$ such that
$P|_{\BNS_t}^\om[S^t]=1$ for all $\om\not\in N_1$.

Let $(\varphi_n)_{n\in\N}$ be a family of test functions which is dense
in $\test$ for the $D(A)$-norm (see Remark \ref{r:dense}). Fix $n\in\N$,
then, since \map{2} holds for $P$, $(M_t^{\varphi_n},\BNS_t,P)_{t\ge0}$
is a continuous $P$-square integrable martingale with quadratic variation
$\zeta_t=t|\Q^{\sss{\frac12}}\varphi_n|_H^2$. By Proposition
\ref{p:condmart} there is a $P$-null set $N_2^n\in\BNS_t$ such that
$(M_s^{\varphi_n},\BNS_s,P|_{\BNS_t}^\om)_{s\ge t}$ is a continuous
$P|_{\BNS_t}^\om$-square integrable martingale with quadratic variation
$\zeta$ for all $\om\not\in N_2^n$. Set $N_2=\bigcup N_2^n$.

We next prove \map{3} for the conditional distributions. First we notice
that, if we set
$$
\al_t^1=|\xi_t|^2_H+2\int_0^t|\xi_s|_V^2\,ds
\qquad\text{and}\qquad
\be_t^1=|\xi_0|^2_H+t\s^2,
$$
$\al_t^1$ is left lower semi-continuous and $\be_t^1$ is non-decreasing,
and so $E_t^1=\al_t^1-\be_t^1$ is also left lower semi-continuous. Hence,
since \map{3} is true for $P$, by virtue of Proposition \ref{p:condsuper},
there is a $P$-null set $N_3\in\BNS_t$ such that
$(E_t^1,\BNS_t,P|_{\BNS_{t_0}}^\om)_{t\ge t_0}$ is an a.\ s.\ super-martingale
for all $\om\not\in N_3$.

As regarding \map{4}, one can proceed in a similar way with the
a.\ s.\ super-martingales $E_t^n$, $n\ge2$. Set
\begin{align*}
\al_t^n&=|\xi_t|^{2n}_H+2n\int_0^t|\xi_s|_H^{2n-2}|\xi_s|_V^2\,ds,\\
\be_t^n&=|\xi_0|^{2n}_H+n(2n-1)\s^2\int_0^t|\xi_s|_H^{2n-2}\,ds,
\end{align*}
again $\al_t^n$ is lower semi-continuous and $\be_t^n$ is increasing,
so that $E^n_t$ is lower semi-continuous. In order to prove integrability
of $\al^n$ and $\be^n$, one has to proceed iteratively, due to the integrals
in both terms. Indeed, integrability of $\al_t^2$, $\be_t^2$ follows
from integrability of $\al_t^1$, $\be_t^1$ and implies integrability of
$\al_t^3$, $\be_t^3$, and so on. Again, by Proposition \ref{p:condsuper},
there are $P$-null sets $N_4^n\in\BNS_t$ such that
$(E_t^n,\BNS_t,P|_{\BNS_{t_0}}^\om)_{t\ge t_0}$ is an a.\ s.\ 
super-martingale for all $\om\not\in N_4^n$. We take
$N_4=\bigcup N_4^n$.

Finally, there is a $P$-null set $N_5\in\BNS_t$ such that
$P|_{\BNS_t}^\om[\xi_t=\om(t)]=1$ for all $\om\not\in N_5$,
and this implies \map{5} for the conditional distributions.
\end{proof}
\begin{lemma}
The \emph{reconstruction} property of Definition \ref{d:aspremarkov}
holds for the family $(\ns(x))_{x\in H}$.
\end{lemma}
\begin{proof}
Let $x\in H$ and $P\in\ns(x)$. Let $T_P$ be the set of exceptional
times of $P$ (see Remark \ref{r:exceptime}) and fix $t\not\in T_P$.
Let $\om\mapsto Q_\om$ be a $\BNS_t$ measurable map defined on
$\Omns$ with values in $\Pr(\Omns^t)$ and let $N_Q\in\BNS_t$
be a $P$-null set such that $\om(t)\in H$ and $Q_\om\in\Phi_t\ns(\om(t))$
for all $\om\not\in N_Q$. We aim to show that the probability measure
$P\otimes_t Q$, defined in Definition \ref{d:raccordo}, is in $\ns(x)$.
First, observe that $(Q_\om)_{\om\in\Omns}$ is a regular conditional
probability distribution of $P\otimes_t Q$ given $\BNS_t$.

We prove \map{1}. By using the notation defined in \eqref{e:sottospazi},
$$
P\otimes_t Q[S_t\cap S^t]
=\int_{S_t}Q_\om[S^t]\,P(d\om)
=P[S_t]
=1,
$$
since $Q_\om[S^t]=1$ holds true due to \map{1} for $Q_\om\in\Phi_t\ns(\om(t))$.

In order to prove \map{2}, let $\varphi\in\test$, then
$(M_s^\varphi,\BNS_s,Q_\om)_{s\ge t}$ is a $Q_\om$-square integrable
martingale for all $\om\not\in N_Q$. By Proposition \ref{p:condmart},
$(M_s^\varphi,\BNS_s,P\otimes_t Q)_{s\ge t}$ is a $P\otimes_t Q$-square
integrable martingale. Since $P$ and $P\otimes_t Q$ agree on $\BNS_t$ and
$(M_s^\varphi,\BNS_s,P)_{0\le s\le t}$ is a martingale, it follows
that $(M_s^\varphi,\BNS_s,P\otimes_t Q)_{s\ge0}$ is a martingale as well.

One can proceed similarly in order to prove properties \map{3} and
\map{4}, by using Proposition \ref{p:condsuper} and the lower
semi-continuity of processes $E_t^n$, $n\ge1$.

Finally, $P$ and $P\otimes_t Q$ agree on $\BNS_t$ and so \map{5} is straightforward.
\end{proof}
\section{Regularity of Markov selections in the initial condition}\label{s:regular}
In this section we prove that, under sufficiently strong non-degeneracy
conditions on the noise (see Assumption \ref{a:noiseass2}), the
martingale solutions that are members of the same Markov selection depend
continuously on the initial conditions. The topology on the initial conditions
is that of $D(A^\theta)$, for a suitable $\theta$ depending on the regularity
of the noise. The space $D(A^\theta)$ will be also denoted by $\W$ in the
sequel. The topology on the solution at time $t$ is that of total
variation of the law; in this sense, it is a continuous dependence in law,
somewhat in the spirit of uniqueness in law. 
\subsection{The intuitive idea}\label{ss:intuitive}
Let $(P_x)_{x\in H}$ be an a.\ s.\ Markov process and denote by $P_{x,t}$ the
law on $H$ of $\xi_t$ under $P_x$ (in different words, the marginal of
$P_x$ at time $t$). For every $t\geq0$ and \emph{almost every} $s\geq0$ we
have
$$
\E^{P_{x,t+s}}[\varphi]=\int_H\E^{P_{y,s}}[\varphi]\,P_{x,t}(dy),
$$
for every $\varphi\in B_b(H)$. With more compact notations, we may write
$$
\E^{P_{x,t+s}}[\varphi]=\E^{P_{x,t}}\bigl[\E^{P_{\cdot,s}}[\varphi]\bigr].
$$
To understand the following idea, fix $\varphi\in B_b(H)$, $s>0$ and a
small $t$. The function $\E^{P_{\cdot,s}}[\varphi]$ is a-priori only
measurable, independently of any imposed additional regularity of $\varphi$.
This is a consequence of our (a-priori) lack of knowledge about continuous
dependence of solutions on initial conditions.

In the deterministic set-up, \emph{for small} $t$ and \emph{regular} $x$ the
function $x\mapsto P_{x,t}=\delta_{u(t;x)}$ would be continuous in the
topology of weak convergence of measures. This does imply nothing on
$\E^{P_{x,t}}\bigl[\E^{P_{\cdot,s}}[\varphi]\bigr]$ since
$\E^{P_{\cdot,s}}[\varphi]$ is only measurable.

On the contrary, one may hope that under proper assumptions on the noise,
again for small $t$ and regular $x$, the function $x\mapsto P_{x,t}$ would
be continuous in the \emph{topology of total variation}. In such a case
$x\mapsto\E^{P_{x,t}}\bigl[\E^{P_{\cdot,s}}[\varphi]\bigr]$ would be
continuous. This is the simple idea behind the main result of this section.

However, there is a difficulty. We cannot show that for sufficiently
small $t$ the function $x\mapsto P_{x,t}$ is continuous over the whole
space $D(A^\theta)$. Even worse, we cannot show that, given
$x_0\in D(A^\theta)$, there is a sufficiently small $t$ (depending on
$x_0$) such that $x\mapsto P_{x,t}$ is continuous at $x_0$ in $D(A^\theta)$.
We can only prove such a statement \emph{up to a given small error}. Namely,
given the error $\epsilon>0$ and $x_0\in D(A^\theta)$, there is a sufficiently
small $t$, depending on both $\epsilon$ and $x_0$, such that the jump of
$x\mapsto P_{x,t}$ around $x_0$ in $D(A^\theta)$ is smaller than the
error $\epsilon$. But via Lemma \ref{l:preliminare} below, this is sufficient.
The idea behind this is to use an approximation by a regularised problem,
which has itself strong Feller solutions. This same technique is usually
applied to handle locally Lipschitz non-linearities in stochastic equations.
\subsection{An abstract result on the strong Feller property}
The intuitive ideas of the previous subsection will reappear in the proof
of Theorem \ref{t:coupling} below. It indirectly formalises the previous
intuition by saying that if a Markov process coincides on a positive random
time with a strong Feller process, then it is strong Feller itself. We could
think of this as a sort of \emph{antithesis} of a coupling result, that can
be called \emph{starting-by-coupling}.

In this section we switch back to the abstract setting of Section
\ref{s:abstract}, to strengthen the idea that the following results
holds true in a more general framework than the Navier-Stokes equations.

Let $(P_x)_{x\in\H}$ be an a.\ s.\ Markov process on $(\Om,\B)$, as in
Definition \ref{d:asmarkov}, and  let $(\Ps_t)_{t\geq0}$ be the associated
transition semi-group on $B_b(\H)$, defined as
\begin{equation}\label{e:tran}
(\Ps_t\varphi)(x)=\E^{P_x}[\varphi(\xi_t)], \qquad x\in\H,\quad \varphi\in B_b(\H).
\end{equation}
The operators $\Ps_t:B_b(\H)\to B_b(\H)$ have the properties
\begin{enumerate}
\item $\Ps_0=I_{B_b(\H)}$,
\item $\|\Ps_t\|_{\mathcal{L}(B_b(\H),B_b(\H))}=1$,
\end{enumerate}
and they don't constitute a proper semi-group, since,
\begin{equation}\label{e:semi-group}
\Ps_{t+s}=\Ps_t\,\Ps_s
\qquad\text{for every }t\geq0\text{ and \emph{almost every} }s\geq0.
\end{equation}
\begin{remark}
We wish to clarify the above formula. The set of \emph{bad times} $s$
where the semi-group property does not hold depends on the point where
the semi-group is evaluated. More precisely, given $x\in\H$, there is
$T\subset(0,\infty)$ of null Lebesgue measure such that
$\Ps_{t+s}\varphi(x)=\Ps_t(\Ps_s\varphi)(x)$ for all $t\ge0$, all
$s\not\in T$ and all $\varphi\in C_b(\H)$.
\end{remark}
Finally, let $\W\subset\H$ be a Banach space with dense continuous injection.
\begin{definition}[$\W$--strong Feller semi-group]
A given semi-group $(\tilde\Ps_t)_{t\geq0}$ on $B_b(\H)$ is \emph{$\W$--strong Feller}
if for every $t>0$ and $\psi\in B_b(\H)$,
$$
\tilde\Ps_t\psi\in C_b(\W).
$$
\end{definition}
\begin{definition}[a.\ s.\ $\W$--Markov process]
A family $(\tilde P_x)_{x\in\H}$ of probability measures on $(\Om,\B)$
is an \emph{a.\ s.\ $\W$--Markov process} if
\begin{enumerate}
\item $\tilde P_x[C([0,\infty);\W)]=1$ for every $x\in\W$,
\item the mapping
$$
x\mapsto(\tilde\Ps_t\varphi)(x):=\E^{\tilde P_x}[\varphi(\xi_t)]
$$
is Borel measurable on $\W$ for every $t\geq0$ and $\varphi\in B_b(\H)$,
\item for every $t\geq0$ and \emph{almost every} $s\geq0$,
$$
\tilde\Ps_{t+s}\varphi=\tilde\Ps_t\tilde\Ps_s\varphi,
$$
for every $\varphi\in B_b(\H)$.
\end{enumerate}
\end{definition}
The next theorem contains the main result of this section. It translates
in a proper way the intuitive \emph{starting-by-coupling} idea explained above.
\begin{theorem}\label{t:coupling}
Let $(P_x)_{x\in\H}$, be an a.\ s.\ Markov process on $(\Om,\B)$ and, for each $R>0$,
let $(P^\er{R}_x)_{x\in\W}$ be an a.\ s.\ $\W$--Markov process on $(\Om,\B)$.

Assume that for every $\rho>0$ there are $R_\rho>0$ and a random time $\tau_\rho$
on $(\Om,\B)$, such that for all $x\in\W$, with $|x|_\W\le\rho$,
\begin{enumerate}
\item $\lim_{\ep\to0}P^\er{R_\rho}_{x+h}[\tau_\rho\ge\ep]=1$, uniformly in $h\in\W$, with $|h|_\W\le1$,
\item for every $t\geq0$, $\varphi\in B_b(\H)$,
$$
\E^{P^\er{R_\rho}_x}[\varphi(\xi_t)\uno_{\{\tau_\rho\geq t\}}]=\E^{P_x}[\varphi(\xi_t)\uno_{\{\tau_\rho\geq t\}}].
$$
\end{enumerate}
If $(\Ps_t^\er{R})_{t\geq0}$ is $\W$--strong Feller for every $R>0$,
then $(\Ps_t)_{t\geq0}$ is $\W$--strong Feller.
\end{theorem}
\subsubsection{Proof of Theorem \ref{t:coupling}}
We divide the proof in several lemmas and some definitions, since some of
them may have independent interest.
\begin{definition}[approximately $\W$--strong Feller]
Let $(\Ps_t)_{t\geq0}$ be a family of linear bounded operators on $B_b(\H)$
and let $\epsilon>0$, $x_0\in\W$ and $t_0>0$.

The family $(\Ps_t)_{t\geq0}$ is \emph{$\W$--strong Feller at $(t_0,x_0)$ up
to the error $\epsilon$} if for every $\eta>0$ there is
$\delta=\delta(x_0,t_0,\epsilon,\eta)>0$ such that
$$
|(\Ps_{t_0}\psi)(x_0+h)-(\Ps_{t_0}\psi)(x_0)|\leq\epsilon+\eta
$$
for every $\psi\in B_b(\H)$, with $|\psi|_\infty\leq1$, and for every
$h\in\W$, with $|h|_\W<\delta$.

If we can choose $\epsilon=0$ in the previous condition, we simply
say that $(P_t)_{t\geq0}$ is $\W$--strong Feller at $(t_0,x_0)$
(\emph{without error}).
\end{definition}
An important detail of the previous definition is the \emph{uniformity} in
$\psi\in B_b(\H)$, $|\psi|_\infty\leq1$. This is the key tool to transfer
the continuity property from small to arbitrary times. It is essentially
here that we implement the intuitive idea of subsection \ref{ss:intuitive}.
\begin{lemma}\label{l:preliminare}
Let $(\Ps_t)_{t\geq0}$ be a family of linear bounded operators on $B_b(\H)$
such that $\|\Ps_t\|_{\mathcal{L}(B_b(\H))}\leq1$ for every $t\geq0$
and such that \eqref{e:semi-group} holds true.

Given $T>0$, assume that for every $\epsilon>0$ and $x_0\in\W$ there is
$t_0\in(0,T)$ such that $(\Ps_t)_{t\geq0}$ is $\W$--strong Feller at
$(t_0,x_0)$ up to the error $\epsilon$.

Then $(\Ps_t)_{t\geq0}$ is $\W$--strong Feller at $(T,x_0)$ for all $x_0\in\W$.
\end{lemma}
\begin{proof}
Given $x_0\in\W$, $T>0$ and $\epsilon>0$, we have to find
$\delta=\delta(x_0,T,\epsilon)>0$ such that
$$
|(\Ps_T\psi)(x_0+h)-(\Ps_T\psi)(x_0)|\leq\epsilon
$$
for every $\psi\in B_b(\H)$, with $|\psi|_\infty\leq1$, and for every
$h\in\W$, with $|h|_\W<\delta$.

Corresponding to such given $x_0$ and $\epsilon$, by assumption we can choose
$t_0\in(0,T)$ such that $(\Ps_t)_{t\geq0}$ is $\W$--strong Feller at $(t_0,x_0)$
up to the error $\frac\epsilon2$. In particular, with the choice
$\eta=\frac\epsilon2$, there is $\delta^*=\delta(x_0,t_0,\frac\epsilon2,\frac\epsilon2)>0$
such that
$$
|(\Ps_{t_0}\psi)(x_0+h)-(\Ps_{t_0}\psi)(x_0)|\leq\epsilon
$$
for every $\psi\in B_b(\H)$, with $|\psi|_\infty\leq1$, and for every $h\in\W$, with
$|h|_\W<\delta^*$.

Fix $h\in\W$ with $|h|_\W<\delta^*$ and let $\Delta$ be the set of points $s\geq0$ such
that $\Ps_{t+s}=\Ps_t\,\Ps_s$ for every $t\geq0$ when evaluated in $x_0$ and $x_0+h$.
For $s\in\Delta$ and $\varphi\in B_b(\H)$, with $|\varphi|_\infty\leq1$, we have
$$
|(\Ps_{t_0+s}\varphi)(x_0+h)-(\Ps_{t_0+s}\varphi)(x_0)|
=|(\Ps_{t_0}\Ps_s\varphi)(x_0+h)-(\Ps_{t_0}\Ps_s\varphi)(x_0)|
\leq\epsilon.
$$
Let $\Delta_{T-t_0}$ be the set of all $s\in\Delta$ such that $T-t_0-s\in\Delta$.
Since $T-t_0>0$ and $\Delta^c$ has null Lebesgue measure, we see that $\Delta_{T-t_0}$
is non empty (in fact it is a subset of $[0,T-t_0]$ of full Lebesgue measure).
Therefore, for $s\in\Delta_{T-t_0}$ and $\psi\in B_b(\H)$, with $|\psi|_\infty\leq1$,
we have
\begin{multline*}
|(\Ps_T\psi)(x_0+h)-(\Ps_T\psi)(x_0)|=\\
\qquad=|(\Ps_{t_0+s}\Ps_{T-t_0-s}\psi)(x_0+h)-(\Ps_{t_0+s}\Ps_{T-t_0-s}\psi)(x_0)|
\leq\epsilon.
\end{multline*}
Therefore, $\delta^*$ is a good choice and the proof is complete.
\end{proof}
\begin{remark}
In fact, we have proved that, given any $x_0\in\W$ and $t>0$, for
every $\epsilon>0$ there is $\delta=\delta(x_0,t,\epsilon)>0$ such that
$$
|(\Ps_t\psi)(x_0+h)-(\Ps_t\psi)(x_0)|\leq\epsilon
$$
for every $\psi\in B_b(\H)$, with $|\psi|_\infty\leq1$, and for every $h\in\W$,
with $|h|_\W<\delta$. This is a stronger property, which is uniform in $\psi\in B_b(\H)$.
\end{remark}
Our aim is then to prove that given $\epsilon>0$ and $x_0\in\W$, there are
arbitrarily small times $t_0$ such that $(\Ps_t)_{t\geq0}$ is $\W$--strong
Feller at $(t_0,x_0)$ up to the error $\epsilon$. This aim is accomplished
if we can prove the hypotheses of the following lemma. Now, consider the
a.\ s.\ semigroups $(\Ps_t)_{t\geq0}$ and $(P_t^\er{R})_{t\geq0}$ as in the statement
of Theorem \ref{t:coupling}.
\begin{lemma}\label{l:cut-off}
Given $x_0\in\W$, assume that for every $\ep>0$ there are $(t_0,R_0)$, with the
possibility to choose $t_0$ arbitrarily small, and $\delta_\ep>0$ such that:
\begin{enumerate}
\item for every $\varphi\in B_b(\H)$, with $|\varphi|_\infty\leq1$, and for every
      $h\in\W$, with $|h|_\W<\delta_\ep$,
$$
|(\Ps_{t_0}^\er{R_0}\varphi)(x_0+h)-(\Ps_{t_0}\varphi)(x_0+h)|\le\ep;
$$
\item $(\Ps_t^\er{R_0})_{t\geq0}$ is $\W$--strong Feller at $(t_0,x_0)$ (\emph{without error}).
\end{enumerate}
Then $(\Ps_t)_{t\geq0}$ is $\W$--strong Feller at $(t,x_0)$ for all $t>0$.
\end{lemma}
The proof is obvious, from triangle inequality and Lemma \ref{l:preliminare}. We
can now proceed to the proof of Theorem \ref{t:coupling}.
\begin{lemma}\label{l:approssima}
Under the assumptions of Theorem \ref{t:coupling}, for
every $\rho>0$, $x\in\W$ with $|x|_\W\leq\rho$, $t\geq0$ and $\varphi\in B_b(\H)$,
$$
|(\Ps_t^\er{R_\rho}\varphi)(x)-(\Ps_t\varphi)(x)|
\leq 2|\varphi|_\infty\,P^\er{R_\rho}_x[\tau_\rho<t],
$$
where $R_\rho$ and $\tau_\rho$ are given by the assumptions.
\end{lemma}
\begin{proof}
We have
$$
(\Ps_t^\er{R_\rho}\varphi)(x)-(\Ps_t\varphi)(x)
=\E^{P^\er{R_\rho}_x}[\varphi(\xi_t)\uno_{\{\tau_\rho<t\}}]-\E^{P_x}[\varphi(\xi_t)\uno_{\{\tau_\rho<t\}}]
$$
hence
$$
|(\Ps_t^\er{R_\rho}\varphi)(x)-(\Ps_t\varphi)(x)|
\leq|\varphi|_\infty(P^\er{R_\rho}_x[\tau_\rho<t]+P_x[\tau_\rho<t]).
$$
Now, again from the assumptions,
$$
P_x[\tau_\rho<t]=1-P_x[\tau_\rho\geq t]
=1-P^\er{R_\rho}_x[\tau_\rho\geq t]=P^\er{R_\rho}_x[\tau_\rho<t]
$$
and this proves the lemma.
\end{proof}
\begin{proof}[Proof of Theorem \ref{t:coupling}]
We only need to prove condition (\textbf{1}) of Lemma \ref{l:cut-off},
since (\textbf{2}) is clearly verified. Given $x\in\W$, we set
$\rho=1+|x|_\W$ and consider $\tau_\rho$ and $R_\rho$ from the
statement of the theorem. Fix $\ep>0$, then there is $t_\ep>0$ such
that $P_{x+h}^\er{R_\rho}[\tau_\rho<t_\ep]<\frac{\ep}{4}$ for all
$h\in\W$ with $|h|_\W\le1$. Since $|x+h|_\W\le\rho$, Lemma
\ref{l:approssima} immediately implies (\textbf{1}) of Lemma
\ref{l:cut-off}.
\end{proof}
\subsection{The strong Feller property for the Navier-Stokes equations}
In this section we apply the abstract result of previous section to the
Markov selections of equations \eqref{e:nseq}. We shall prove that the
strong Feller property holds for Markov selections under the assumption
on the noise given below. Indeed, as we shall see in Proposition \ref{p:BEL},
Assumption \ref{a:noiseass2} below ensures that the solutions to a
regularised version of the Navier-Stokes equations have the strong
Feller property. In view of possible extensions to weaker assumptions
on the noise (e.\ g., noise less regular, more degenerate, etc.),
we remark that all it is needed here to apply the abstract results
of the previous section are exactly the conclusions of Proposition
\ref{p:BEL}, and Assumption \ref{a:noiseass2} below is only a way
to ensure this.
\subsubsection{Further assumptions}
In order to introduce the assumptions that we need,
we recall that $\Q$ denotes the covariance of the driving noise
that we have introduced at the beginning of Section \ref{s:martingale}.
\begin{assumption}[Noise is non-degenerate and regular]\label{a:noiseass2}
There are an isomorphism $\Q_0$
of $H$ and a number $\al_0>\frac16$ such that
$$
\Q^{\frac12}=A^{-\frac34-\al_0}\Q_0^{\frac12}.
$$
\end{assumption}
We first remark that the assumption above implies Assumption
\ref{a:noiseass1}. Indeed, since the embedding of
$H^{\sss{\frac32}+2\ep}(\T)$ into $L^2(\T)$ is Hilbert-Schmidt
in dimension three for every $\ep>0$, and
$A^{-\sss{\frac34}-\ep}$ is a bounded operator from
$H$ to $H^{\sss{\frac32}+2\ep}(\T)$, the operator
$A^{-\sss{\frac34}-\ep}$ is Hilbert-Schmidt in $H$,
for every $\ep>0$.

Moreover, $A^{-\sss{\frac34}-\ep}{\Q_0}^{\sss{\frac12}}W(t)$
is a Brownian motion in $H$, for every $\ep>0$ and
every isomorphism $\Q_0$ of $H$, where $W(t)$ is a cylindrical
Wiener process on $H$ (see Da Prato \& Zabczyk \cite{DPZa2}).
In conclusion, $A^{-\sss{\frac34}-\al_0}{\Q_0}^{\sss{\frac12}}W(t)$
is a Brownian motion in $D(A^\al)$ for every $\al_0>\al\geq0$.

The restriction $\al_0>\frac16$ is needed essentially
for technical reasons in the proof of Theorem \ref{t:weakstrong}
(see Section \ref{ss:aapprox}).
\subsubsection{The main continuity result}
Let us introduce the function $\theta:(0,\infty)\to(0,\infty)$ defined as
\begin{equation}\label{e:theta}
\theta(\al)=
\begin{cases}
\frac12+\frac\al2  &\qquad\text{if }0<\al\leq\frac12\\
\frac14+\al        &\qquad\text{if }\al>\frac12.
\end{cases}
\end{equation}

Given $\al_0$, we set
$$
\W=D(A^{\theta(\al_0)})
\qquad\text{and}\qquad
|x|_\W:=|A^{\theta(\al_0)}x|.
$$
The general result of the previous section is applied to the Navier-Stokes
equations using the above spaces. It is the content of the following theorem.
\begin{theorem}\label{t:maincont}
Assume that the covariance $\Q$ satisfies Assumption \ref{a:noiseass2}.
Let $(P_x)_{x\in H}$ be \emph{any} a.\ s.\ Markov process associated to the
stochastic Navier-Stokes equations \eqref{e:nseq} and let $(\Ps_t)_{t\geq0}$
be the operators on $B_b(H)$ defined as in \eqref{e:tran}. Then
$(\Ps_t)_{t\geq0}$ is $\W$--strong Feller.
\end{theorem}
In view of Theorem \ref{t:coupling}, we follow the approach of Flandoli
\& Maslowski \cite{FlMa} to construct $P^\er{R}_x$. We introduce an
equation which differs from the original one by a cut-off only, so that
with large probability they have the same trajectories on a small
deterministic time interval. We consider the equation
\begin{equation}\label{e:eqR}
du+\bigl[Au+B(u,u)\chi_R(|u|_\W^2)\bigr]\,dt=\Q^{\frac12}\,dW,
\end{equation}
where $R\ge1$ and $\chi_R:[0,\infty)\to[0,1]$ is a non-increasing smooth
function equal to $1$ over $[0,R+1]$, to $0$ over $[R+2,\infty)$, and with
derivative bounded by $1$.
\begin{theorem}[Weak-strong uniqueness]\label{t:weakstrong}
For every $x\in\W$, equation \eqref{e:eqR} has a unique martingale solution
$P^\er{R}_x$, with
$$
P^\er{R}_x[C([0,\infty);\W)]=1.
$$
Let $\tau_R:\Omns\to[0,\infty]$ be defined as
\begin{equation}\label{e:stoppingtime}
\tau_R(\om)=\inf\{\,t\geq0\,:\,|\om(t)|_\W^2\geq R\,\},
\end{equation}
and $\tau_R(\om)=\infty$ if this set is empty. If $x\in\W$ and $|x|_\W^2<R$,
then
\begin{equation}\label{e:errore}
\lim_{\ep\to0}P^\er{R}_{x+h}[\tau_R\ge\ep]=1,
\qquad\text{uniformly in }h\in\W,\ |h|_\W<1.
\end{equation}
Moreover, on $[0,\tau_R]$, the probability measure $P^\er{R}_x$ coincides with any
martingale solution $P_x$ of the original stochastic Navier-Stokes
equations \eqref{e:nseq}, namely
\begin{equation}\label{e:weakstrong}
\E^{P^\er{R}_x}[\varphi(\xi_t)\uno_{\{\tau_R\geq t\}}]
=\E^{P_x}[\varphi(\xi_t)\uno_{\{\tau_R\geq t\}}],
\end{equation}
for every $t\geq0$ and $\varphi\in B_b(H)$.
\end{theorem}
The proof of this theorem is given in Appendix \ref{ss:aapprox}, where we shall give
a slightly stronger form of \eqref{e:weakstrong}. In order to apply
Theorem \ref{t:coupling} we only need the following result.
\begin{proposition}\label{p:BEL}
For every $R>0$, the transition semi-group $(\Ps_t^\er{R})_{t\geq0}$ associated to
equation \eqref{e:eqR} is $\W$--strong Feller.
\end{proposition}
Compared to Flandoli \& Maslowski \cite{FlMa}, Da Prato \& Debussche \cite{DPDe},
Da Prato \& Zabczyk \cite{DPZa} and other references, the proof of this proposition
is classical, but since this equations does not fall in any of the cases considered
in the above references, we provide the main lines of proof in Appendix \ref{ss:aapprox}.
The proof of Theorem \ref{t:maincont} is therefore complete.
\section{Some consequences of strong Feller regularity}\label{s:consequences}
In this section we show that the regularity results of the
previous section provide some additional results on
Markov selections for the Navier-Stokes equations.

We shall prove that the strong Feller property allows
to improve the point-wise regularity of the solutions.
As a consequence, the set of exceptional times (see
Remark \ref{r:exceptime}) of the martingale solutions
starting from a regular initial condition is empty
and the energy inequalities of Definition \ref{d:asNSsol}
hold for all times. Moreover, for such regular initial
conditions, the Markov property holds for every time.
Finally, we shall prove a \emph{condition} for
global well-posedness.

Throughout this section, we shall work under Assumptions
\ref{a:noiseass2} on the covariance of the noise (which,
via Theorem \ref{t:maincont}, ensures that any Markov solution
has the $\W$--strong Feller property).
\subsection{A support theorem}\label{ss:support}
In this part of the section we show a support theorem for
Markov solutions. Actually, the result below is not a
consequence of the strong Feller property, but we will use
it in the rest of the section. As it concerns the proof,
we follow closely Flandoli \cite{Fla3}; we also use some
results on the approximate problem from Appendix
\ref{ss:control}.

We preliminarily introduce a notation. We say that a Borel
probability measure $\mu$ on $H$ is \emph{fully supported}
on $\W$ if $\mu(U)>0$ for every open set $U\subset\W$.
\begin{proposition}\label{p:support}
Let $(P_x)_{x\in H}$ be an a.\ s.\ Markov selection and assume
Assumption \ref{a:noiseass2} on the covariance. For every $x\in\W$
and every $T>0$, the image measure of $P_x$ at time $T$ is fully
supported on $\W$ (as defined above).
\end{proposition}
\begin{proof}
Fix $x\in\W$ and $T>0$. We need to show that for every $y\in\W$
and $\ep>0$, $P_x[|\xi_T-y|_\W<\ep]>0$. Let $\ov{y}\in\W$ such
that $A\ov{y}\in\W$ and $|y-\ov{y}|_\W<\frac{\ep}2$. Choose
$R>0$ such that $3|x|^2_\W<R$ and $3|y|_\W^2<R$.
Then by Theorem \ref{t:weakstrong},
\begin{align*}
P_x[|\xi_T-y|<\ep]
&\ge P_x[|\xi_T-\ov{y}|_\W<\frac\ep2]
 \ge P_x[|\xi_T-\ov{y}|_\W<\frac\ep2,\ \tau_R>T]\\
&=   P_x^\er{R}[|\xi_T-\ov{y}|_\W<\frac\ep2,\ \tau_R>T]
\end{align*}
By Lemma \ref{l:appcon}, there is a control $\ov{w}\in W^{s,p}([0,T];D(A^{\al_1}))$,
with $s$, $p$ and $\al_1$ chosen as in Lemma \ref{l:concon}, such that the solution
$\ov{u}$ to the control problem \eqref{e:ceq} corresponding to $\ov{w}$ satisfies
$$
\ov{u}(0)=x,
\qquad
\ov{u}(T)=\ov{y}
\quad\text{and}\quad
|\ov{u}(t)|_\W^2\le\frac23R.
$$
By Lemma \ref{l:concon}, there exists $\delta>0$ such that for all $w\in W^{s,p}([0,T];D(A^{\al_1}))$
with $|w-\ov{w}|_{W^{s,p}(0,T;D(A^{\al_1}))}<\delta$, we have
$$
|u(T,w)-\ov{y}|_\W<\frac\ep2
\qquad\text{and}\qquad
\sup_{t\in[0,T]}|u(t,w)|_\W^2<R,
$$
where each $u(\cdot,w)$ is the solution to the controlled problem
\eqref{e:ceq} corresponding to $w$ and starting at $x$. Hence
$$
P_x^\er{R}[|\xi_T-\ov{y}|_\W<\frac\ep2,\ \tau_R>T]
\ge P_x^\er{R}[|\eta_T-\ov{w}|_{W^{s,p}(0,T;D(A^{\al_1}))}<\delta]>0
$$
where $\eta_t=\xi_t-x+\int_0^t(A\xi_s+\chi_R(|\xi_s|^2_W)B(\xi_s,\xi_s))\,ds$
and the probability on $\eta$ above is positive since by Assumption \ref{a:noiseass2}
$\eta$ is a Brownian motion in $D(A^{\al_1})$, with full support
on $W^{s,p}(0,T;D(A^{\al_1}))$.
\end{proof}
\begin{remark}
We remark that if a probability $\mu$ is \emph{fully supported}
on $\W$, it means that its support \emph{contains} $\W$. As
we shall see in Theorem \ref{t:regular} below, the support
of the marginal of $P_x$ at every positive time is indeed $\W$.
\end{remark}
\subsection{Improved regularity}
This part of the section is devoted to the proof of the following theorem
and some additional useful results. For technical reasons (see Lemma
\ref{l:regbase}) we shall restrict to the case $\al_0>\frac12$, so that
$\th(\al_0)>\frac34$.
\begin{theorem}\label{t:regular}
Assume Assumption \ref{a:noiseass2} with $\al_0>\frac12$ and
let $(P_x)_{x\in H}$ be a Markov solution. For every $x\in\W$
and every $t>0$,
$$
P_x[\text{there is $\ep>0$ such that }\xi\in C([t-\ep,t+\ep];\W)]=1.
$$
\end{theorem}
We point out that in the above formula the radius $\ep$ of the
neighbourhood is random and depends on $\om\in\Omns$. In particular, this
theorem does not really improve our knowledge on the global regularity
of trajectories. Anyway, it turns out to be useful for Theorem
\ref{t:allmarkov} and the following result (an immediate consequence
of the improved regularity).
\begin{corollary}
Under the assumptions of the previous theorem, for every $x\in\W$
the set of exceptional times (see Remark \ref{r:exceptime}) of
$P_x$ is empty and for each $n\ge1$, the process $(E_t^n)_{t\ge0}$
(given in Definition \ref{d:asNSsol}) is a super-martingale under $P_x$.
\end{corollary}
\begin{proof}
The above theorem implies that, for a fixed $s$, $E^n_{s_k}\to E^n_s$
$P_x$--a.\ s., if $s_k\uparrow s$. As in the proof of Lemma \ref{l:exsuper},
each $(E^n_{s_k})_{k\in\N}$ is uniformly integrable and these facts
easily imply the conclusion.
\end{proof}

A standard technique to analyse regularity is to consider stationary solutions
and then disintegrate them. The stationary solutions have uniform average
bounds in regular topologies coming from stationarity and energy inequality;
sometimes they imply some additional a.\ s.\ regularity of paths at given time
$t$ (see Flandoli \& Romito \cite{FlRo2} for an example); then by disintegration
the same regularity is transferred to solutions with deterministic initial
conditions $u_0$, for a.\ e.\ $u_0$ with respect to the time-zero marginal
of the stationary solution. If irreducibility holds, the set of such $u_0$
is dense and the strong Feller property allows to extend the result to every
$u_0$.

The implementation of this program in our case is difficult due to a number
of issues. The first one is that, given an a.\ s.\ Markov process $(P_x)_{x\in H}$,
one needs to construct a stationary solution associated to the process (namely,
whose disintegration at every time produces elements of the family $(P_x)_{x\in H}$).
Without stationarity, the program above does not start.

We describe here a modification of this first step. We construct a solution
which is almost stationary and is associated to the given $(P_x)_{x\in H}$
by construction. It seems that it is not easy to prove that a suitable
limit of these almost stationary solutions is a true stationary solution,
but we conjecture that it should be true. Anyway, for the sequel of the
program, the solutions that we construct are sufficient. In particular,
they provide an example of non-stationary solution with uniform average
bounds in regular topologies, which is, as far as we know, a not entirely
trivial fact.

Given an a.\ s.\ Markov process $(P_x)_{x\in H}$, we say that a solution
$P$ of the martingale problem with initial condition a given probability
measure $\nu$ is \emph{associated to} $(P_x)_{x\in H}$ if
$$
P=\int_H P_x\,\nu(dx).
$$
Define, for $s<t$, the event
$$
\Omns^\W(s,t)=\{\xi\in C([s,t];\W)\}
$$
and, for every $t\ge0$,
\begin{align*}
\Omns^\W(t)
&=\{\text{there is $\ep>0$ such that }\xi\in C([t-\ep,t+\ep];\W)\}\\
&=\bigcup_{\ep>0}\Omns^\W(t-\ep,t+\ep).
\end{align*}
\begin{lemma}\label{l:quasistaz}
Under the assumptions of Theorem \ref{t:regular}, for every a.\ s.\ Markov
process $(P_x)_{x\in H}$ there is a measure $\wt{\nu}$, fully supported
on $\W$ (as defined in Section \ref{ss:support}) such that the solution
$\wt{P}$ with initial condition $\wt{\nu}$ associated to $(P_x)_{x\in H}$
has the property
$$
\wt{P}[\Omns^\W(t)]=1,
\qquad\text{for all $t>0$.}
$$
\end{lemma}
\begin{proof}
Let $\wt{\nu}=\int_0^1\Ps_s^*\delta_0\,ds$. Following Chow
and Khasminskii \cite{ChHa}, we have
\begin{equation}\label{e:boundonKB}
\E^{\wt{\nu}}[\,|\cdot|_V^2]
=\int_0^1\E^{P_0}|\xi_s|_V^2\,ds
\le\frac12\Tr\Q
\end{equation}
(see the definition of solution of the martingale problem). By Proposition
\ref{p:support}, $\wt{\nu}$ is fully supported on $\W$. Let $\wt{P}$ be
the solution to the martingale problem with initial distribution $\wt{\nu}$
associated to $(P_x)_{x\in H}$. The basic fact is the following bound.
For every $t>0$,
\begin{align*}
\E^{\wt{P}}[|\xi_t|^2_V]
&=\int_H\E^{P_x}[|\xi_t|^2_V]\,\wt{\nu}(dx)
=\int_H\E^{\Ps_t^*\delta_x}[|\cdot|_V^2]\,\wt{\nu}(dx)\\
&=\int_0^1\E^{\Ps_{t+s}^*\delta_0}[|\cdot|_V^2]\,ds
=\int_0^1\E^{P_0}[|\xi_{s+t}|_V^2]\,ds\\
&=\int_t^{1+t}\E^{P_0}[|\xi_r|_V^2]\,dr
\le \frac{1+t}{2}\Tr\Q.
\end{align*}
Fix $t>0$ and $\ep>0$ such that $t-\ep>0$, then by disintegrating
$\wt{P}$ at time $t-2\ep$, for every $R>0$ sufficiently large,
\begin{align*}
\wt{P}[\Omns^\W(t-\ep,t+\ep)]
&=  \int_H P_y[\Omns^\W(\ep,3\ep)]\,\wt{P}_{t-2\ep}(dy)\\
&\ge\int_{\{|y|_V^2\le R\}} P_y[\Omns^\W(\ep,3\ep)]\,\wt{P}_{t-2\ep}(dy)\\
&\ge\Bigl(\inf_{|y|_V^2\le R}P_y[\Omns^\W(\ep,3\ep)]\Bigr)\wt{P}[|\xi_{t-2\ep}|_V^2\le R]\\
&\ge\Bigl(\inf_{|y|_V^2\le R}P_y[\Omns^\W(\ep,3\ep)]\Bigr)(1-\frac{1+t}{2R}\Tr\Q)
\end{align*}
where $\wt{P}_{t-2\ep}$ is the marginal of $\wt{P}$ at time $t-2\ep$.

So, we only need to prove that, for every $R$ large enough, 
$$
\lim_{\ep\to0}\inf_{|y|_V^2\le R}P_y[\Omns^\W(\ep,3\ep)]=1.
$$
This property follows easily from Lemma \ref{l:regbase}
of Appendix \ref{ss:modified}, since the Lemma gives an
estimate of the time where the solutions are regular depending
only on $R$ and on the size of the noise.
\end{proof}
We can now prove Theorem \ref{t:regular}.
\begin{proof}[Proof of Theorem \ref{t:regular}]
Let $\wt{P}$ be the solution associated to $(P_x)_{x\in H}$
provided by the previous lemma. We know, by the same lemma, that
$\wt{P}[\Omns^\W(t)]=1$, hence by disintegration it follows that
$P_x[\Omns^\W(t)]=1$ for $\wt{\nu}$--a.\ e.\ $x$, where $\wt{\nu}$
is the marginal of $\wt{P}$ at time $0$. We use the strong Feller
property to show that the conclusion is true for all $x$. Indeed,
if $x_n\to x$ in $\W$, then
\begin{multline*}
P_{x_n}[\Omns^\W(t)]
=\int P_{\om(s)}[\Omns^\W(t-s)]\,P_{x_n}(d\om)
\longrightarrow\\\longrightarrow
 \int P_{\om(s)}[\Omns^\W(t-s)]\,P_x(d\om)
=P_x[\Omns^\W(t)],
\end{multline*}
and the conclusion follows.
\end{proof}
\subsection{The Markov property for all times}
The regularity result of the previous section allows
to prove that the Markov property holds for all
times. We assume again, as in the previous section,
that $\al_0>\frac12$. First we need the following lemma.
\begin{lemma}
Assume the hypotheses of Theorem \ref{t:regular} and let
$(P_x)_{x\in H}$ be an a.\ s.\ Markov process. Given
$\varphi\in C_b(H)$ and $r_n$, $r\ge0$, with $r_n\to r$,
define
$$
\psi_n(x)=\E^{P_x}[\varphi(\xi_{r_n})]
\qquad\text{and}\qquad
\psi(x)=\E^{P_x}[\varphi(\xi_r)]
$$
Then $(\psi_n)_{n\in\N}$ converges to $\psi$ in $C_b(\W)$.
\end{lemma}
\begin{proof}
Given $x\in\W$, we know, by Theorem \ref{t:regular} above,
that $\xi_{r_n}\to\xi_r$, $P_x$--a.\ s.\ in $\W$, and so
$\psi_n(x)\to\psi(x)$. Since $(\psi_n)_{n\in\N}$ are
uniformly bounded, we only need to prove that the
sequence $(\psi_n)_{n\in\N}$ is equi-continuous.
Ascoli-Arzel\`a theorem next concludes the proof.
But, equi-continuity is a direct consequence
of strong Feller property (Theorem \ref{t:maincont}).
\end{proof}
\begin{theorem}\label{t:allmarkov}
Under Assumption\ref{a:noiseass2}, with $\al_0>\frac12$,
if $(P_x)_{x\in H}$ is an a.\ s.\ Markov process,
then $(P_x)_{x\in\W}$ is a Markov process.

More precisely, for every $x\in\W$, every
$\varphi\in C_b(H)$ and every $0\le s\le t$,
\begin{equation}\label{e:allmarkov}
\E^{P_x}[\varphi(\xi_t)|\BNS_s]
=\E^{P_{\xi_s}}[\varphi(\xi_{t-s})],
\qquad P_x\text{-a.\ s.}
\end{equation}
\end{theorem}
\begin{proof}
Let $(P_x)_{x\in H}$ be an a.\ s.\ Markov process in $H$,
we know that it is $\W$--strong Feller. Given $x\in\W$,
let $T\subset(0,\infty)$ be such that for every
$\varphi\in C_b(H)$, every $s\not\in T$ and $t>s$,
equality \eqref{e:allmarkov} holds.

Fix $s\in(0,\infty)$ and $\varphi\in C_b(H)$, we
show that the above equality \eqref{e:allmarkov} holds
for $s$. Let $s_n\downarrow s$, with $s_n<t$ and
$s_n\not\in T$, and define, as in the previous lemma,
$$
\psi_n(x)=\E^{P_x}[\varphi(\xi_{t-s_n})]
\qquad\text{and}\qquad
\psi(x)=\E^{P_x}[\varphi(\xi_{t-s})].
$$
Since $\psi_n$ converges uniformly to $\psi$ and
$\xi_{s_n}\to\xi_s$, $P_x$--a.\ s.\ in $\W$ (by Theorem
\ref{t:regular}), we have that
$\psi_n(\xi_{s_n})\to\psi(\xi_{s})$, $P_x$--a.\ s.

Notice that $s_n\not\in T$, for all $n\in\N$, and so,
by \eqref{e:allmarkov} for $s_n$,
$$
\E^{P_x}[\varphi(\xi_t)\uno_A]
=\E^{P_x}[\psi_n(\xi_{s_n})\uno_A]
\longrightarrow\E^{P_x}[\psi(\xi_s)\uno_A]
=\E^{P_x}\bigl[\E^{P_{\xi_s}}[\varphi(\xi_{t-s})]\uno_A\bigr],
$$
for all $A\in\BNS_s$. In conclusion, $\psi(\xi_s)$
is the conditional expectation of $\varphi(\xi_t)$,
given $\BNS_s$, or, in different terms,
\eqref{e:allmarkov} holds for $s$.
\end{proof}
\subsection{A condition for well-posedness}\label{ss:wellposed}
In this last part of the section we show that the $\W$--strong
Feller property implies global well-posedness in $\W$ (that is,
for all initial conditions and all times) if there exists a
single initial condition for which the problem is well-posed
up to a \emph{deterministic} time. We give two versions of
this result, the first for path-wise uniqueness, the second
for uniqueness in law.
\begin{theorem}
Under Assumption \ref{a:noiseass2} on the covariance, assume
that there are $x_0\in\W$, $t_0>0$ and a solution $\wt{P}_{x_0}$
to the martingale problem starting at $x_0$, such that
$$
\widehat{P_{x_0}}[C([0,t_0];\W)]=1.
$$
Then, for every Markov selection $(P_x)_{x\in H}$,
$$
P_x[C([0,\infty);\W)]=1
$$
for every $x\in\W$. In particular, path-wise
uniqueness holds for every $x\in\W$.
\end{theorem}
\begin{proof}
We recall that, by Proposition \ref{p:support}, the marginal
of $P$ at each time is \emph{fully supported} on $\W$, for
every martingale solution $P$.

Let $(P_x)_{x\in H}$ a Markov selection. Then $P_{x_0}$ and
$\widehat{P_{x_0}}$ coincide on $[0,t_0]$. In particular,
$$
P_{x_0}[C([0,t_0];\W)]=1.
$$

Take any $\ep\in(0,t_0)$. By disintegration,
$$
P_x[C([0,t_0-\ep];\W)]=1
$$
for $P_{x_0,\ep}$--a.\ e.\ $x$, hence on a dense set
$D\subset\W$, where $P_{x_0,\ep}$ is the marginal at
time $\ep$ of $P_{x_0}$.

We show that $C([0,t_0);\W)$ is a full set for $P_x$,
for every $x\in\W$. Indeed, for all $x\in D$, we know
that $P_x[C([\frac1n,t_0-\ep];\W)]=1$,
for every $n\in\N$. By the strong Feller property, the
same property is true on the whole $\W$. In conclusion,
$C((0,t_0);\W)$ is a full set for $P_x$, for all $x\in\W$.
Finally, since there is a strictly positive random time
$\tau_x$ such that all trajectories starting from $x$ are
continuous with values in $\W$ (see Theorem \ref{t:weakstrong}),
we can deduce that $P_x[C([0,t_0);\W)]=1$.

We can conclude the proof. Given $x\in\W$, by the Markov
property we know that for every $s>0$,
$$
P_x[C([s,s+t_0);\W)]=1,
$$
hence $P_x[C([0,2t_0);\W)]=1$, and, by repeating this
argument, we obtain $P_x[C([0,\infty);\W)]=1$.

Since from any $x$ there is at least one regular solution
(see Lemma \ref{l:regbase}), we have path-wise uniqueness.
\end{proof}
We conclude this section by giving a condition for uniqueness
in law. Again, it holds only for regular conditions.
\begin{proposition}\label{p:lawuniq}
Under Assumption \ref{a:noiseass2}, consider two Markov selections
$(P_x^\er{1})_{x\in H}$ and $(P_x^\er{2})_{x\in H}$. If there are
$x_0\in\W$ and $t_0>0$ such that $P_{x_0}^\er{1}=P_{x_0}^\er{2}$
on $[0,t_0]$, then $P_1^\er{1}=P_x^\er{2}$ for all $x\in\W$.
\end{proposition}
\begin{proof}
For every $\ep\in(0,t_0)$, the disintegrations of $P_{x_0}^\er{1}$
and $P_{x_0}^\er{2}$ at time $\ep$ are the same, up to time $t_0$.
Since they are solutions over $[0,t_0-\ep]$ with initial conditions
chosen at random by the common law $\mu_\ep$ at time $\ep$ of both
$P_{x_0}^\er{1}$ and $P_{x_0}^\er{2}$, we can deduce that
$P_y^\er{1}$ and $P_y^\er{2}$ coincide over $[0,t_0-\ep]$ for
$\mu_\ep$--a.\ e. $y$. As in the proof of previous theorem,
there is a dense set of initial conditions $y$ in $\W$ such that
the two measure $(P_x^\er{i})_{i=1,2}$ coincide on the interval $[0,t_0-\ep]$.

By the strong Feller property, the values of the laws on the closure
of such dense set, hence on $\W$, are determined by the values on
the dense set. Hence $P_x^\er{1}=P_x^\er{2}$ for every $x\in\W$ on
the interval $[0,t_0-\ep]$.

By the Markov property, the values over a finite deterministic time
can be replicated up to every time. Thus the two selections coincide.
\end{proof}
Using the properties of regular conditional probability distributions
the following result follows easily from the above proposition.
\begin{corollary}\label{c:lawuniq}
Under Assumption \ref{a:noiseass2}, let $(P_x^\er{1})_{x\in H}$ and
$(P_x^\er{2})_{x\in H}$ be two Markov selections. If there is an initial
distribution on $\W$ such that the two solutions corresponding to the
given selections coincide, then $P_x^\er{1}=P_x^\er{2}$ for all $x\in\W$.
\end{corollary}
We remark that, as in Stroock \& Varadhan \cite[Theorem 12.2.4]{StVa},
uniqueness of the martingale problem follows from uniqueness of
the Markov selection (see also Remark \ref{r:singleton}), so the
previous result implies immediately the following criterion for
uniqueness in law.
\begin{theorem}\label{t:lawuniq}
Under Assumption \ref{a:noiseass2} on the covariance,
assume that there are $x_0\in H$ and $t_0>0$ such that there is only one
solution to the martingale problem with initial condition $x_0$
on the interval $[0,t_0]$. Then for every $x\in\W$ there is
only one solution to the martingale problem with initial
condition $x$.
\end{theorem}
\appendix
\section{Existence for the martingale problem}\label{s:exist}
This section has a twofold purpose. On one hand, we prove Theorem \ref{t:mainex}
on existence for solutions to the martingale problem.
On the other hand, the same line of demonstration works, almost flawlessly,
for Lemma \ref{l:compmeas} and, in this way, we lighten its proof.
\subsection{Proof of Theorem \ref{t:mainex}}\label{ss:exist}
We assume the hypotheses of Theorem \ref{t:mainex}, that is, we are given a probability
measure $\mu_0\in\Pr(H)$ with all finite moments. We aim to show that there exists a
solution, starting at $\mu_0$, to the martingale problem associated to the Navier-Stokes
equations \eqref{e:nseq}, as defined in Definition \ref{d:asNSsol}. The proof of this
result will be developed in independent general steps, suitable to be applied in Lemma
\ref{l:compmeas}. In the sequel we shall prove, for the sake of conciseness, only those
properties of the definition of solution to the martingale problem that present some
novelty with respect to the known literature. For all other results, we refer to
Flandoli \& Gatarek \cite{FlGa}.

Consider the Galerkin approximations of Navier-Stokes equations \eqref{e:nseq},
\begin{equation}\label{e:gnseq}
du_t^n+[A_n u_t^n+B_n(u_t^n)]\,dt=\sum_{i=1}^n\s_i e_i\,d\be_i(t),
\end{equation}
where $A_n$ and $B_n$ are the projection of operators $A$ and $B$ onto the linear
space $H_n=\Span[e_1,\dots,e_n]$ and $(\be_i)_{i\in\N}$ are independent one-dimensional
Brownian motions. Similarly, the process solution to \eqref{e:gnseq} has initial
distribution given by the projection onto $H_n$ of $\mu_0$.
Let $P^G_n$ be the law of $u_n$ on $\Omns$. Standard arguments (see for example
the proof of Theorem $3.1$ of Flandoli \& Gatarek \cite{FlGa}) show that
\eqref{e:gnseq} has a unique strong solution, for each $n\in\N$, and that
the sequence $(P^G_n)_{n\in\N}$ satisfies all assumptions of Lemma \ref{l:excomp}
below, once we take $\xi_1^n(t)=-\int_0^t(A_n\xi_s+B_n(\xi_s))\,ds$.
\begin{lemma}\label{l:excomp}
Given a sequence $(P_n)_{n\in\N}$ of probability measures on $\Omns$, assume that
there are $\al\in(0,\frac12)$, $p>\frac1\al$ and $\gamma\in(\frac32,2)$ such that
for all $T>0$,
\begin{enumerate}
\item $\E^{P_n}[\|\xi\|_{L^\infty(0,T;H)}+\|\xi\|_{L^2(0,T;V)}]\le C(T)$,
\item $\E^{P_n}[\|\xi_1^n\|_{W^{1,2}(0,T;D(A^{-\gamma}))}]\le C(T)$,
\item $\E^{P_n}[\|\xi-\xi_1^n\|_{W^{\al,p}(0,T;H)}]\le C(T)$,
\end{enumerate}
where $C(T)$ is independent of $n$ and $\xi_1^n$ is a suitable adapted process.

Then $(P_n)_{n\in\N}$ is tight in $\Omns\cap L^2_\loc([0,\infty);H)$.
\end{lemma}
Notice that, although it is not proved in Flandoli \& Gatarek \cite{FlGa},
it is easy to see that each $P_n^G$ fulfils the properties corresponding
to \map{1}, \dots, \map{5} associated with equation \eqref{e:gnseq}\footnote{We
do not write down such corresponding properties for the sake
of brevity. Truly, only \map{2} really needs to be patched, and its formulation
can be found in Lemma \ref{l:exmart}.}.

We set
\begin{equation}\label{e:Uspace}
\mathcal{U}=\Omns\cap L^2_\loc([0,\infty);H),
\end{equation}
by Lemma \ref{l:excomp} above, $(P_n^G)_{n\in\N}$ is tight in $\mathcal{U}$,
hence there is a sub-sequence $(P_{k_n}^G)_{n\in\N}$ converging weakly to
some probability measure $P^G_\infty$ on $\mathcal{U}$. It is now sufficient
to show that $P^G_\infty$ is a solution to the martingale problem associated
to the Navier-Stokes equations (say, Definition \ref{d:asNSsol}). First,
\map{5} is obvious, given the choice of the initial condition for
equation \ref{e:gnseq}. Next, $P^G_\infty$ satisfies \map{2} by the lemma below.
\begin{lemma}\label{l:exmart}
The measure $P_\infty^G$ defined above satisfies property \map{2} of
Definition \ref{d:asNSsol}.
\end{lemma}
\begin{proof}
Let $\varphi\in\test$ and set
$$
\zeta(t)=t|\Q^\sss{\frac12}\varphi|_H^2
\qquad\text{and}\qquad
\zeta_n(t)=t|\Q_n^\sss{\frac12}\varphi|_H^2,
$$
where $\Q_n$ is the projection of $\Q$ onto $H_n$. We know that,
under $P^G_{k_n}$,
$$
M_t^{\varphi,k_n}=\langle\xi_t-\xi_0,\varphi\rangle
+\int_0^t\langle\xi_s,A_n\varphi\rangle\,ds
-\int_0^t\langle B_n(\xi_s,\varphi),\xi_s\rangle\,ds
$$
is a continuous martingale, with quadratic variation $\zeta_{k_n}$. By
the L\'evy martingale characterisation, $(M_t^{\varphi,k_n})_{t\ge0}$
is a Brownian motion. Now, let $0\le s<t$ and $A\in\BNS_s$, we know that
\begin{align*}
&\E^{P^G_{k_n}}[\uno_A(M_t^{\varphi,k_n}-M_s^{\varphi,k_n})]=0,\\
&\E^{P^G_{k_n}}[\uno_A\bigl((|M_t^{\varphi,k_n}|^2-\zeta_n(t))-(|M_s^{\varphi,k_n}|^2-\zeta_n(s))\bigr)]=0,\\
&\sup_{n\in\N}\E^{P^G_{k_n}}[|M_t^{\varphi,k_n}|^{2+\ep}]<\infty,
\qquad\text{for all }t,
\end{align*}
since $\zeta_n(t)\uparrow\zeta(t)$. Since
$M_t^{\varphi,k_n}\to M_t^\varphi$ for all $\om\in\Omns$ (see for example
the proof of Theorem 3.1 in Flandoli \& Gatarek \cite{FlGa}) and
$P^G_{k_n}\rightharpoonup P^G_\infty$ weakly in $\mathcal{U}$, we have
that
\begin{align*}
&\E^{P^G_\infty}[\uno_A(M_t^\varphi-M_s^\varphi)]=0,\\
&\E^{P^G_\infty}[\uno_A\bigl((|M_t^\varphi|^2-\zeta_n(t))-(|M_s^\varphi|^2-\zeta_n(s))\bigr)]=0,
\end{align*}
or, in different words, \map{2} holds for $P^G_\infty$.
\end{proof}
Finally, the proof of Theorem \ref{t:mainex} is complete once
we prove that $P^G_\infty$ fulfils properties \map{1}, \map{3}
and \map{4}. This is ensured by the following lemma.
\begin{lemma}\label{l:exsuper}
Let $P_n$, $P\in\Pr(\Omns)$ such that
\begin{enumerate}
\item $P_n\rightharpoonup P$ in $\mathcal{U}$,
\item properties $\map{1}$, $\map{3}$ and $\map{4}$ hold for $P_n$, for all $n$,
\item for each $m\ge1$ there is a constant $c_m$ (independent of $n$) such that
$$
\E^{\mu_n}|x|^{2m}_H\le c_m,
\qquad m\ge1,
$$
where $\mu_n$ is the marginal at time $0$ of $P_n$.
\end{enumerate}
Then properties \map{1}, \map{3} and \map{4} hold for $P$.
\end{lemma}
\begin{proof}
First, we notice that by \map{3} and \map{4} for $P_n$, for all $m\ge1$
and all $t\ge0$,
\begin{equation}\label{e:unifbounds}
\E^{P_n}|\xi_t|^{2m}_H\le C(m,t,\s^2,c_1,\dots,c_m)
\quad\text{and}\quad
\E^{P_n}\int_0^t|\xi_s|_V^2\,ds\le C(t,\s^2,c_1).
\end{equation}
For each $T>0$ and $\la>0$ set
$$
A_{T,\la}=\{\sup_{t\in[0,T]}\bigl(|\xi_t|_H^2+2\int_0^t|\xi_s|^2_V\,ds\bigr)>\la\},
$$
this set is open in $\Omns$ by semi-continuity of the above norms in the topology
of $\Omns$, so that using Corollary \ref{c:doob} on $E_t^1$, we have,
$$
\la P[A_{T,\la}]
\le \liminf_{n\to\infty}\la P_n[A_{T,\la}]
\le C(T,\s^2,c_1),
$$
and so \map{1} holds for $P$.

We next prove \map{3}. Since $P_n\rightharpoonup P$ in $\mathcal{U}$, by
Skorokhod's theorem we know that there are a probability space
$(\Sigma,\mathscr{G},{\abP})$ and random variables $X_n$, $X$ on $\Sigma$ with
values in $\mathcal{U}$ such that $X_n\to X$ in $\mathcal{U}$ ${\abP}$--a.\ s.,
and $X_n$, $X$ have laws $P_n$, $P$ respectively.

First, $E_t^1$ is $P$-integrable since, by \eqref{e:unifbounds}, for every
$t>0$,
$$
\E^P|\xi_t|_H^2
  = \E^{\abP}|X(t)|_H^2
\le \liminf_{n\to\infty}\E^{\abP}|X_n(t)|_H^2
\le C(t,\s^2,c_1)
$$
and similarly
$$
\E^P\![\int_0^t|\xi_s|^2_V\,ds]
  = \E^{\abP}\![\int_0^t|X(s)|^2_V]
\le \liminf_{n\to\infty}\,\E^{\abP}[\int_0^t|X_n(s)|^2_V]
\le C(t,\s^2,c_1).
$$
Then, in order to prove the a.\ s.\ super-martingale property, we just need
to show the following claim.
\begin{center}
\begin{minipage}{0.9\linewidth}
Given $t>0$, there is a set $T_t\subset(0,t)$ of null
Lebesgue measure such that for all $s\not\in T_t$ and all positive bounded
continuous functions $\phi$ on $C([0,s];D(A)')$,
\begin{equation}\label{e:exsuper}
\E^{P}[\phi E^1_t]\le\E^P[\phi E_s^1].
\end{equation}
\end{minipage}
\end{center}
Indeed, if the above claim is true, we can first deduce from it that
$\E^P[E_t^1-E_s^1|\BNS_s]\le0$ for all such $s$. Then, if we set
$T=\bigcup_{t\in D} T_t$, where $D\subset[0,\infty)$ is a countable dense set,
it is easy to see that $E^1$ is an a.\ s.\ super-martingale having $T$ as its set
of exceptional times. In fact, if $t>s>0$ and $s\not\in T$, there is a sequence
$t_n\uparrow t$ in $D$ with $t_n\ge s$. Since $s\not\in T_{t_n}$, \eqref{e:exsuper}
holds for $t_n$ and so by Fatou's lemma (we use actually the version of Fatou's lemma
given at \eqref{e:newfatou}) inequality \eqref{e:exsuper} holds for $t$ as well.

We prove now the above claim. Fix $t>0$ and a positive bounded continuous
function $\phi$ on $C([0,s];D(A)')$. Since $X_n\to X$ ${\abP}$--a.\ s.\ in
$\Omns$, by lower semi-continuity of the time integral of the $V$-norm, we know
that for all $s\in[0,t]$,
$$
\phi(X)\int_s^t|X(r)|_V^2\,dr
\le\liminf_{n\to\infty}\phi(X_n)\int_s^t|X_n(r)|_V^2\,dr,
\qquad {\abP}-\text{a.\ s.}
$$
and so, by Fatou's lemma,
\begin{equation}\label{e:ingred1}
\E^{\abP}[\phi(X)\int_s^t|X(r)|_V^2\,dr]
\le\liminf_{n\to\infty}\E^{\abP}[\phi(X_n)\int_s^t|X_n(r)|_V^2\,dr].
\end{equation}
Next, since $X_n\to X$ ${\abP}$--a.\ s.\ in $L^2(0,t;H)$, and since $(X_n)_{n\in\N}$
is uniformly integrable in $L^2(\Omns,L^2(0,t;H))$ (thanks to \eqref{e:unifbounds}),
it follows that
$$
\int_0^t\E^{\abP}|X_n(r)-X(r)|_H^2\,dr
=\E^{\abP}\int_0^t|X_n(r)-X(r)|_H^2\,dr\longrightarrow0.
$$
Hence there are a sub-sequence $(X'_n)_{n\in\N}$ and a set $T_t\subset[0,t]$
(both independent of $\phi$), with null Lebesgue measure, such that for
all $s\not\in T_t$,
\begin{equation}\label{e:asH}
\E^{\abP}|X'_n(s)-X(s)|_H^2\longrightarrow0,
\qquad\text{as }n\to\infty.
\end{equation}
Since $X'_n\to X$ ${\abP}$--a.\ s.\ in $\Omns$, by proceeding as for \eqref{e:ingred1}, we
get
\begin{equation}\label{e:ingred21}
\E^{\abP}[\phi(X)|X(t)|_H^2]
\le\liminf_{n\to\infty}\E^{\abP}[\phi(X'_n)|X'_n(t)|_H^2].
\end{equation}
Now, let $s\not\in T_t$, then
\begin{multline*}
\E^{\abP}[\phi(X'_n)|X'_n(s)|_H^2]-\E^{\abP}[\phi(X)|X(s)|_H^2]=\\
=\E^{\abP}[\phi(X'_n)(|X'_n(s)|_H^2-|X(s)|_H^2)]+\E^{\abP}[(\phi(X'_n)-\phi(X))|X(s)|_H^2],
\end{multline*}
the first term on the right-hand side converges to zero by \eqref{e:asH}, since
$\phi$ is bounded, while the second term converges to zero by Lebesgue theorem,
since $\phi$ is continuous and $X_n'\to X$ ${\abP}$--a.\ s.\ in $\Omns$. In conclusion
\begin{equation}\label{e:ingred22}
\E^{\abP}[\phi(X'_n)|X'_n(s)|_H^2]\longrightarrow\E^{\abP}[\phi(X)|X(s)|_H^2],
\qquad n\to\infty,
\end{equation}
and \eqref{e:ingred1}, \eqref{e:ingred21} and \eqref{e:ingred22} together
ensure \eqref{e:exsuper}. We finally remark that $T_t$ can be chosen in such
a way that $0\not\in T_t$, since $E^1_0=0$.

The proof of \map{4} is entirely similar to the proof of \map{3} given above,
we only need to prove that for all $s\le t$,
$$
\E^{\abP}[\int_s^t|X_n-X|^{2m-2}_H\,ds]\longrightarrow0,
\qquad n\to\infty,
$$
so that we can conclude that
$$
\E^{\abP}[\phi(X_n)\int_s^t|X_n(r)|_H^{2m-2}\,dr]
\longrightarrow
\E^{\abP}[\phi(X)\int_s^t|X(r)|_H^{2m-2}\,dr],
$$
the former is a consequence of the convergence in $L^2(\Omns,L^2(0,t;H))$
and the bounds \eqref{e:unifbounds}.
\end{proof}
\section{Some results on disintegration and reconstruction}\label{s:tecnic}
In this section we give a few technical results concerning
a.\ s.\ super-martingales, mostly used in the proofs of
Theorem \ref{t:mainmarkov}. We start by showing that
a stopped a.\ s.\ super-martingale is again
an a.\ s.\ super-martingale.
\begin{proposition}\label{p:stoppata}
Given an a.\ s.\ super-martingale $(\theta_t,\B_t,P)_{t\ge0}$,
if $\tau$ is a stopping time with respect to $(\B_t)_{t\ge0}$, then
the stopped process $\th^\tau_t=\th_{t\wedge\tau}$ is again an
a.\ s.\ super-martingale with respect to the same filtration and
having the same exceptional set $T_\th$ of $\th$.
\end{proposition}
\begin{proof}
The result is classical for super-martingales (see for example
Proposition 1.3 in Revuz \& Yor \cite{ReYo}. In order to prove
the same for a.\ s.\ super-martingales, we simply observe that
if $(t_n)_{n\in\N}$ is an increasing sequence of times in
$T_\th^c$, then $(\th_{t_n})_{n\in\N}$ is a super-martingale.
So for example, to prove the a.\ s.\ super-martingale property
for the stopped process, fix $s\not\in T_\th$ and $t>s$, pick
an increasing sequence containing both $s$ and $t$ and use the
above remark to deduce that for all $A\in\B_s$,
$$
\E[\th_t^\tau\uno_A]\le\E[\th_s^\tau\uno_A].
$$
One can prove similarly that $\th_t^\tau$ is $\B_t$-measurable
and that $\E[|\th_t^\tau|]<\infty$.
\end{proof}
The next proposition shows that the martingale property is conserved under
disintegration and reconstruction. We give the proof for completeness to
clarify the details related to a.\ s.\ super-martingales, even though the
statement is similar to Theorem $1.2.10$ of Stroock \& Varadhan \cite{StVa}.
\begin{proposition}\label{p:condmart}
Given $P\in\Pr(\Om)$, two continuous adapted processes $\theta$,
$\zeta:[0,\infty)\times\Om\to\R$ and $t_0\ge0$, the following conditions are
equivalent:
\begin{itemize}
\item[\textit{(i)}] $(\theta_t,\B_t,P)_{t\ge t_0}$ is a $P$-square integrable
martingale with quadratic variation $(\zeta_t)_{t\ge t_0}$;
\item[\textit{(ii)}] there is a $P$-null set $N\in\B_{t_0}$ such that for all
$\om\not\in N$, the process $(\theta_t,\B_t,P|_{\B_t}^\om)_{t\ge t_0}$ is a
$P|_{\B_t}^\om$-square integrable martingale with quadratic variation
$(\zeta_t)_{t\ge t_0}$ and $\E^P[\E^{P|_{\B_t}^\cdot}[\zeta_t]]<\infty$.
\end{itemize}
\end{proposition}
\begin{proof}
Assume \textit{(i)}. First, we prove that if $t_0\le t_1<t_2$, then there is
a $P$-null set $N_{t_1t_2}\in\B_{t_0}$ such that for all $\om\not\in N_{t_1t_2}$,
\begin{equation}\label{e:*}
\E^{P|_{\B_{t_0}}^\om}[\theta_{t_2}|\B_{t_1}]=\theta_{t_1},
\qquad P|_{\B_{t_0}}^\om-a.\ s.
\end{equation}
Indeed, let $A\in\B_{t_1}$, then for each $B\in\B_{t_0}$ we have that $A\cap B\in\B_{t_1}$
and
$$
\E^P[\uno_B\E^{P|_{\B_{t_0}}^\cdot}[\theta_{t_2}\uno_A]]
=\E^P[\theta_{t_2}\uno_{A\cap B}]
=\E^P[\theta_{t_1}\uno_{A\cap B}]
=\E^P[\uno_B\E^{P|_{\B_{t_0}}^\cdot}[\theta_{t_1}\uno_A]]
$$
so that $\E^{P|_{\B_{t_0}}^\om}[\theta_{t_2}\uno_A]=\E^{P|_{\B_{t_0}}^\om}[\theta_{t_1}\uno_A]$
out of a $P$-null set in $\B_{t_0}$. Since $\B_{t_0}$ is countably generated,
the $P$-null set can be chosen independently of $A$.

Next, let $D$ be a dense set in $[t_0,\infty)$, then by the previous argument
we can find a $P$-null set $N\in\B_{t_0}$ such that \eqref{e:*} is true for
$\om\not\in N$ and $t_1$, $t_2\in D$. By Lemma $1.2.9$ of Stroock \& Varadhan
\cite{StVa}, \eqref{e:*} is true for all $t\ge t_0$.

One can proceed similarly to prove that $\theta_t$ is $P|_{\B_t}^\om$-square
integrable with quadratic variation $(\zeta_t)_{t\ge t_0}$, since $\theta_t^2$
is a sub-martingale and $\theta_t^2-\zeta_t$ is a martingale. Finally,
$\E^P[\E^{P|_{\B_t}^\cdot}[\zeta_t]]=\E^P[\zeta_t]$.

Vice versa, assume \textit{(ii)}. Since $\om\mapsto P|_{\B_{t_0}}^\om[\zeta_t]$
is $P$-integrable, $\theta_t$ is $P$-square integrable and it is easy to see
that, by disintegration on $\B_{t_0}$, $\theta_\cdot$ is a martingale with
quadratic variation $\zeta_\cdot$.
\end{proof}
\begin{proposition}[Doob's maximal inequality]
Let $\theta:[0,\infty)\times\Omega\to\R$ be an adapted left lower semi-continuous
process and assume that $(\theta_t,\B_t,P)_{t\ge0}$ is an a.\ s.\ super-martingale.
Let $T_\theta$ be the set of exceptional times of $\theta$ and let
$[a,b]\subset[0,\infty)$, with $a$, $b\not\in T_\theta$. Then for each
$\la>0$,
\begin{equation}\label{e:doobsup}
\la P[\sup_{t\in[a,b]}\theta_t\ge\la]
\le \E^P[\theta_a]+\E^P[\theta_b^-].
\end{equation}
Moreover, under the assumption of uniform integrability of $(\theta_t^-)_{t\in[a,b]}$,
\end{proposition}
\begin{proof}
We prove the first inequality, under the assumption that $b\not\in T_\theta$.
It is easy to see that the inequality holds if we replace
$\sup_{t\in[a,b]}\theta_t$ with $\sup_{t\in D}\theta_t$, where $D$ is a
countable set of non-exceptional times containing both $a$
and $b$ (see for example Theorem $9.4.1$ of Chung \cite{Chu}). If the
countable set is dense in $[a,b]$, by the left lower semi-continuity
of $\theta$ it follows that $\sup_{t\in[a,b]}\theta_t=\sup_{t\in D}\theta_t$
and the inequality is proved. 

If $b\in T_\theta$, one can use a sequence $b_n\uparrow b$, with
$b_n\not\in T_\theta$. The left-hand side of \eqref{e:doobsup} converges,
as $n\to\infty$ by monotone convergence. Since $\theta$ is left lower
semi-continuous, it follows that $\theta^-$ is left upper semi-continuous
so that, by virtue of uniform integrability, an extended Fatou's lemma
(see for example Theorem $7.5.2$ of Ash \cite{Ash}) gives that
$\limsup_{n\to\infty}\E^P[\theta_{b_n}^-]\le\E^P[\theta_b^-]$.
\end{proof}
We state an easy consequence of the above inequality, since we use the maximal
inequality in this form.
\begin{corollary}\label{c:doob}
Under the assumptions of the previous proposition, assume also that
$\theta_t=\al_t-\be_t$, where $\al$ and $\be$ are positive and $\be$ is
non-decreasing. Then for each $a\not\in T_\theta$ and $b\in[0,\infty)$,
with $a<b$,
$$
\la P[\sup_{t\in[a,b]}\al_t\ge\la]
\le 2(\E^P\theta_a+\E^P\theta_b^-+\E^P\be_b),
\qquad\la>0.
$$
\end{corollary}
\begin{proof}
Notice that
\begin{align*}
P[\sup\al_t\ge\la]
&\le P[\sup\theta_t+\be_b\ge\la]\\
&=   P[\sup\theta_t+\be_b\ge\la,\be_b\le\frac{\la}2]
    +P[\sup\theta_t+\be_b\ge\la,\be_b\ge\frac{\la}2]\\
&\le P[\sup\theta_t\ge\frac\la2]
    +P[\be_b\ge\frac\la2]\\
&\le \frac2\la(\E^P[\theta_a]+\E^P[\theta_b^-])+\frac2\la\E^P[\be_b].
\end{align*}
Finally, $b$ can be also in $T_\theta$ since uniform integrability is given
by the estimate $\theta_t^-\le\be_b$, for each $t\le b$.
\end{proof}
We aim to show that the a.\ s.\ super-martingale property is conserved under
disintegration and reconstruction. Prior to this, it is necessary to state a
fairly simple extension to Fatou's lemma, which fits our needs. Assume
that the process $\theta$ of the above proposition is given as
$\theta_t=\al_t-\beta_t$, where $\al$, $\be:[0,\infty)\times\Om\to[0,\infty)$
and $\be$ is non-decreasing. If $t_n\uparrow t$, then
\begin{equation}\label{e:newfatou}
\E^P[\theta_t]\le\liminf_{n\to\infty}\E^P[\theta_{t_n}].
\end{equation}
Indeed, by monotone convergence, $\E^P[\be_{t_n}]\to\E^P[\sup\be_{t_n}]\le\E^P[\be_t]$.
Using the left lower semi-continuity of $\theta$, one gets
$$
\liminf_{n\to\infty}\al_{t_n}
=  \liminf_{n\to\infty}\theta_{t_n}+\sup_{n\in\N}\be_{t_n}
\ge\theta_t+\sup_{n\in\N}\be_{t_n},
$$
so that by Fatou's lemma,
\begin{align*}
\E^P[\theta_t]
&\le\E^P[\liminf_{n\to\infty}\al_{t_n}]-\E^P[\sup_{n\in\N}\be_{t_n}]\\
&\le\liminf_{n\to\infty}\E^P[\al_{t_n}]-\lim_{n\to\infty}\E^P[\be_{t_n}]\\
&=\liminf_{n\to\infty}\E^P[\theta_{t_n}].
\end{align*}
\begin{proposition}\label{p:condsuper}
Let $\al$, $\be:[0,\infty)\times\Omega\to\R_+$ be two adapted processes such
that $\be$ is non-decreasing and
$$
\theta=\al-\be
$$
is left lower semi-continuous. Given $P\in\Pr(\Om)$ and $t_0\ge0$, the following
conditions are equivalent.
\begin{itemize}
\item[\textit{(i)}] $(\theta_t,\B_t,P)_{t\ge t_0}$ is an a.\ s.\ super-martingale
and for all $t\ge t_0$,
$$
\E^P[\al_t+\be_t]<\infty,
$$
\item[\textit{(ii)}] there is a $P$-null set $N\in\B_{t_0}$ such that for all
$\om\not\in N$ the process $(\theta_t,\B_t,P|_{\B_{t_0}}^\om)_{t\ge t_0}$
is an a.\ s.\ super-martingale and for all $t\ge t_0$,
$$
\E^{P|_{\B_{t_0}}^\om}[\al_t+\be_t]<\infty
\qquad\text{and}\qquad
\E^P[\E^{P|_{\B_{t_0}}^\cdot}[\al_t+\be_t]]<\infty.
$$
\end{itemize}
\end{proposition}
\begin{proof}
We preliminarily show that for each $s$ and $t$, with $t_0\le s\le t$, there
is a $P$-null set $N_{s,t}\in\B_{t_0}$ such that for all $\om\not\in N_{s,t}$,
\begin{equation}\label{e:dis-claim}
\E^{P|_{\B_{t_0}}^\om}[\theta_s-\theta_t|\B_s]
=\E^P[\theta_s-\theta_t|\B_s],
\qquad P|_{\B_{t_0}}^\om-\text{a.\ s.}
\end{equation}
For simplicity, we denote by $Z_\om$ the r.v.\ on the left-hand side, and by
$Z$ the r.v.\ on the right-hand side. Both $Z$ and $Z_\om$ are $\B_s$-measurable.
For each $A\in\B_s$ and $B\in\B_{t_0}$, we have
$$
\E^P[\uno_{A\cap B}(\theta_s-\theta_t)]
=\E^P[\uno_{A\cap B}Z]
=\E^P[\uno_B\E^{P|_{\B_{t_0}}^\om}[\uno_AZ]]
$$
and, at the same time,
$$
\E^P[\uno_{A\cap B}(\theta_s-\theta_t)]
=\E^P[\uno_B\E^{P|_{\B_{t_0}}^\om}[\uno_A(\theta_s-\theta_t)]]
=\E^P[\uno_B\E^{P|_{\B_{t_0}}^\om}[\uno_AZ_\om]],
$$
hence there is a $P$-null set $N_{s,t,A}\in\B_{t_0}$ such that for
all $\om\not\in N_{s,t,A}$,
$$
\E^{P|_{\B_{t_0}}^\om}[\uno_AZ]
=\E^{P|_{\B_{t_0}}^\om}[\uno_AZ_\om]
$$
Since $\B_s$ is countably generated, we can find a single $P$-null
set $N_{s,t}\in\B_{t_0}$ such that the above equality holds for all
$\om\not\in N_{s,t}$ and all $A$ in a countable set of generators.
In conclusion, $Z_\om=Z$, $P|_{\B_{t_0}}^\om$--a.\ s., for all
$\om\not\in N_{s,t}$ and the claim \eqref{e:dis-claim} is proved.

Assume \textit{(i)}. Since $\E^P[\E^{P|_{\B_{t_0}}^\om}[\be_t]]=\E^P[\be_t]<\infty$,
for each $t>t_0$ there is a $P$-null set $N_{1,t}\in\B_{t_0}$ such that
$\E^{P|_{\B_{t_0}}^\om}[\be_t]<\infty$ for all $\om\not\in N_{1,t}$.
Choose a sequence $t_n\uparrow\infty$ and let $N_1=\bigcup_{n\in\N}N_{1,t}$.
Since $\be$ is non-decreasing, we have that $\E^{P|_{\B_{t_0}}^\om}[\be_t]<\infty$
for all $\om\not\in N_1$. In particular, for  all $\om\not\in N_1$,
$\be\in L^1_\loc([t_0,\infty);L^1(\Om,P|_{\B_{t_0}}^\om))$.

Now, consider again a sequence $t_n\uparrow\infty$, and notice that
for each $n\in\N$,
$$
\E^P[\int_{t_0}^{t_n}\E^{P|_{\B_{t_0}}^\om}[\al_s]\,ds]
=   \int_{t_0}^{t_n}\E^P[\al_s]
\le \int_{t_0}^{t_n}\E^P[\theta_s+\be_s]
\le (t_n-t_0)\E^P[\theta_{t_0}+\be_{\overline{t}}].
$$
From this it follows, as above, that there is a $P$-null set $N_2\in\B_{t_0}$ such that
$\al\in L^1_\loc([t_0,\infty);L^1(\Om,P|_{\B_{t_0}}^\om))$ for $\om\not\in N_2$.
In particular, $\E^{P|_{\B_{t_0}}^\om}[\al_t]<\infty$ for a.\ e.\ $t\ge t_0$,
for each $\om\not\in N_2$. At this stage of the proof the null set of times
does depend on $\omega$, but what we know is enough for the next computations.

Fix $t>t_0$, then for every Borel set $T\subset [t_0,t]$, $A\in\B_{t_0}$ and
$B\in\B_t$, we have
$$
\E^P\Bigl[\int_{t_0}^t\uno_T\uno_{A\cap B}\E^P[\theta_s-\theta_t|\B_s]\,ds\Bigr]\ge0.
$$
By using the generalised Fatou's lemma \eqref{e:newfatou} and Fubini theorem, it follows that
$$
\E^P\Bigl[\uno_A\E^{P|_{\B_{t_0}}^\om}\bigl[\int_{t_0}^t\uno_T\uno_B\E^{P|_{\B_{t_0}}^\om}[\theta_s-\theta_t|\B_s]\,ds\bigr]\Bigr]\ge0.
$$
so that there is a $P$-null set $N_{t,T,B}\in\B_{t_0}$ such that for all
$\om\not\in N_{t,T,B}$,
$$
\int_{t_0}^t\uno_T\E^{P|_{\B_{t_0}}^\om}\bigl[\uno_B\E^{P|_{\B_{t_0}}^\om}[\theta_s-\theta_t|\B_s]\bigr]\,ds\ge0.
$$
Since both the $\s$-algebra of Borel sets of $[t_0,t]$ and $\B_t$ are countably
generated, there is a single $P$-null set $N_t\in\B_{t_0}$ such that the above
inequality holds for all $T$ and $B$ running over a set of generators of their
respective $\s$-algebras. Fix $\om\not\in N_t$, then there is a set
$T_t(\om)\subset[t_0,t)$ with null Lebesgue measure such that for all
$s\not\in T_t(\om)$,
$$
\E^{P|_{\B_{t_0}}^\om}[\theta_s-\theta_t|\B_s]\ge0
\qquad P|_{\B_{t_0}}^\om-\text{a.\ s.}
$$
Let $D\subset(t_0,\infty)$ be a countable dense set and let $N_3=\bigcup_{t\in D}N_t$.
Fix $\om\not\in N_3$ and set $T_\theta(\om)=\bigcup_{t\in D}T_t(\om)$. We show
that $T_\theta(\om)$ is the set of exceptional times of
$(\theta_t,\B_t,P|_{\B_{t_0}}^\om)_{t\ge t_0}$. Let $t>s>t_0$, with $s\not\in T_\theta(\om)$,
and take a sequence $t_n\uparrow t$ of points of $D$. By the simple extension
to Fatou's lemma \eqref{e:newfatou} given above, it follows that for each $A\in\B_s$,
$$
\E^{P|_{\B_{t_0}}^\om}[\uno_A\theta_t]
\le\liminf_{n\to\infty}\E^{P|_{\B_{t_0}}^\om}[\uno_A\theta_{t_n}]
\le\E^{P|_{\B_{t_0}}^\om}[\uno_A\theta_s]
$$
and so $\E^{P|_{\B_{t_0}}^\om}[\theta_t|\B_s]\le\theta_s$. Finally, we show that
we can choose $N_3$ in such a way that $t_0\not\in T_\theta(\om)$. Indeed, one can
proceed as above with $s=t_0$ and without integrating with respect to time, getting
an additional $P$-null set out of which everything is fine.

We finally set $N=N_1\cup N_2\cup N_3\in\B_{t_0}$ and such $P$-null set is the one
needed for the validity of \textit{(ii)}. In order to conclude the proof of \textit{(ii)},
we only need to show that $\E^{P|_{\B_{t_0}}^\om}[\al_t]<\infty$ holds for all $t\ge t_0$
and all $\om\not\in N$. Fix $\om\not\in N$, then for every $t\not\in T_\theta$ such that
$\E^{P|_{\B_{t_0}}^\om}[\al_t]<\infty$ (by the choice of $N_2$ above, this is true
for a.\ e.\ $t$), we have
$$
\E^{P|_{\B_{t_0}}^\om}[\al_t]
=  \E^{P|_{\B_{t_0}}^\om}[\theta_t]+\E^{P|_{\B_{t_0}}^\om}[\be_t]
\le\E^{P|_{\B_{t_0}}^\om}[\theta_{t_0}]+\E^{P|_{\B_{t_0}}^\om}[\be_t],
$$
and so the semi-continuity argument concludes the proof.

Assume next that statement \textit{(ii)} holds. We show first that for all
$T>t_0$, $\al$, $\be\in L^1(t_0,T;L^1(\Om,P))$. Indeed,
$\E^P\be_t\le\E^P\be_T$ for all $t\in[t_0,T]$, and for all
$t\not\in T_\theta$,
\begin{align*}
\E^P[\al_t]
&=   \E^P[\theta_t+\be_t]
 \le \E^P[\E^{P|_{\B_{t_0}}^\cdot}[\theta_t]]+\E^P[\be_T]\\
&\le \E^P[\E^{P|_{\B_{t_0}}^\cdot}[\theta_{t_0}]]+\E^P[\be_T]
 =   \E^P[\theta_{t_0}]+\E^P[\be_T].
\end{align*}

For a fixed a time $t>t_0$ we show that there is a set
$T_{\theta,t}\subset[t_0,t]$ with zero Lebesgue measure, such that
$\E^P[\theta_s-\theta_t|\B_s]\ge0$, $P$--a.\ s., for all $s\not\in T_{\theta,t}$.
For each $A\in\B_t$, $\om\not\in N$ and $s\not\in T_\theta(\om)$, where $T_\theta(\om)$
is the set of exceptional times of $(\theta_t,\B_t,P|_{\B_{t_0}}^\om)$, we know that
$$
\E^{P|_{\B_{t_0}}^\om}[\uno_A\E^{P|_{\B_{t_0}}^\om}[\theta_s-\theta_t|\B_s]]\ge 0,
$$
so that for all Borel set $T\subset[t_0,t]$ we have that
$$
\E^P[\int_{t_0}^t\uno_T\E^{P|_{\B_{t_0}}^\om}[\uno_A\E^{P|_{\B_{t_0}}^\om}[\theta_s-\theta_t|\B_s]]\,ds]\ge0
$$
All such integrals are finite thanks to the integrability properties of $\al$ and $\be$.
By using \eqref{e:newfatou} and Fubini theorem, we get that
$$
\int_{t_0}^t\uno_T\E^P[\uno_A\E^P[\theta_s-\theta_t|\B_s]]\,ds\ge0.
$$
Hence there is a set $T_{\theta,t,A}\subset[t_0,t)$ of null Lebesgue measure such that
for $s\not\in T_{\theta,t,A}$, $\E^P[\uno_A\E^P[\theta_s-\theta_t|\B_s]]\ge0$. Since
$\B_t$ is countably generated, it is possible to find, as before, a set
$T_{\theta,t}$ of null Lebesgue measure such that $\E^P[\theta_s-\theta_t|\B_s]\ge0$
holds $P$--a.\ s.\ for all $s\not\in T_{\theta,t}$.

Now, as in the first part of this proof, we can show that the set $T_{\theta,t}$ can
be found to be independent of $t$, by using a countable dense set of times and
lower semi-continuity. Similarly, $t_0\not\in T_\theta$.
\end{proof}
\section{Analysis of some equations related to Navier-Stokes}\label{s:related}
This section is divided into three parts. In the first part we prove
a regularity result for the auxiliary equation \ref{e:eqxv} below,
that we use in Section \ref{s:consequences} to show the regularity
improvement. In the second part we prove a couple of results
stated in Section \ref{s:regular} on equation \eqref{e:eqR}.
Finally, in the third part we show some controllability
properties for the same equation \eqref{e:eqR} used in
Section \ref{s:consequences}.

In this section we will be quite sloppy in the management of
constants in the various inequalities. The value of constants will
change from line to line, but we use the same symbol. We shall
only state the quantities they depend upon.
\subsection{Regularity boost for an auxiliary equation}\label{ss:modified}
Given $z\in C([0,T];H)$, with $z(0)=0$, and $u_0\in H$, consider
the equation
\begin{equation}\label{e:eqxv}
\begin{cases}
\frac{dv}{dt}+Av+B(v+z,v+z)=0,\\
v(0)=u_0.
\end{cases}
\end{equation}
We say that $v\in L^\infty(0,T;H)\cap C([0,T];D(A)')$ is a solution to
the above problem if
\begin{align*}
 \ps{v(t)}{\varphi}_H
+\int_0^t\ps{v(s)}{A\varphi}_H\,ds
-\int_0^t\ps{B(v+z,\varphi)}{v+z}_H\,ds
=\ps{u_0}{\varphi}_H
\end{align*}
for every $\varphi\in D(A)$ with bounded gradient. Notice that
all terms are well defined under the given regularities.
The following result shows that nice bounds of the solution
to the above equation in $V$ allow to improve its regularity up
to the one allowed by the data.
\begin{lemma}\label{l:regbase}
Let $u_0\in V$ and $z\in C([0,T];D(A^\gamma))$ be given,
with $\gamma>\frac34$. Then there exists $T_*>0$, depending
only on $|u_0|_V$ and on the norm of $z$ in $C([0,T];D(A^\gamma))$,
such that the equation \eqref{e:eqxv} has a solution
\begin{equation}\label{e:regforv}
v_*\in C([0,T_*];V)\cap L^2(0,T_*;D(A)).
\end{equation}
Such a solution is unique in the class $L^\infty(0,T;H)\cap C([0,T];D(A)')$.
Moreover,
\begin{equation}\label{e:addreg}
v_*\in C((0,T_*];D(A^\gamma)).
\end{equation}
\end{lemma}
\begin{proof}
We do not give all the rigorous details (they require a careful use
of Galerkin approximations for the existence part and finite dimensional
projection for the uniqueness) but only the formal estimates behind the
results. 

\underline{\textsl{Step 1}.} The basic computation for existence is
an a-priori estimate in $V$,
\begin{align*}
\frac{d}{dt}|v|_V^2+2|Av|_H^2
&  = -2\ps{Av}{B(v+z,v+z)}_H\\
&\le 2|A^{\frac34+\ep}v|_H\,|A^{\frac14-\ep}B(v+z,v+z)|_H\\
&\le C_B |A^{\frac34+\ep}v|_H\, |A^{\frac34-\frac\ep2}(v+z)|_H^2\\
&\le |Av|_H^2+C(B,\ep)(|v|_V+|A^{\frac34}z|_H)^{p_\ep},
\end{align*}
(for a suitable exponent $p_\ep>0$) where we have used Lemma
\ref{l:Breg} with $\al=\sss{\frac12}-\ep$ and interpolation
inequalities, and $\ep>0$ is small enough. The proof of
the existence statement then easily follows.

\underline{\textsl{Step 2}.} As to uniqueness, if
$v_1\in L^\infty(0,T_*;H)\cap C([0,T_*];D(A)')$ is another
solution to \eqref{e:eqxv} and $w=v_*-v_1$, then,
using the embedding $|v|_{L^3}\le|A^{\sss{\frac14}}v|_H$,
\begin{align*}
\frac12\frac{d}{dt}|w|_H^2+|w|_V^2
&=  -\ps{w}{B(w,v_*+z)+B(v_1+z,w)}_H
 =   \ps{v_*+z}{B(w,w)}_H\\
&\le |w|_V\, |w|_{L^3(\T)}|v_*+z|_{L^6(\T)}
 \le |w|_V^{\frac32}|w|_H^{\frac12}|v_*+z|_V\\
&\le |w|_V^2+|v_*+z|_V^4|w|_H^2,
\end{align*}
hence by Gronwall's lemma, $w\equiv0$.

\underline{\textsl{Step 3}.} In this step, we prove \eqref{e:addreg}
under the assumption $\gamma=\frac{n}2$, where $n\ge2$ is an integer.
We shall use the following claim:
\begin{center}
\framebox{\begin{minipage}{0.8\linewidth}
if for $t_0\in(0,T_*)$ and for $m\in\N$, with $2\le m\le n$ we have
$v_*(t_0)\in D(A^{\frac{m}2})$, then
$$
v_*\in C([t_0,T_*];D(A^{\frac{m}2}))\cap L^2([t_0,T_*];D(A^{\frac{m+1}2})).
$$
\end{minipage}}
\end{center}
Since $v_*\in L^2(0,T_*;D(A))$, it follows that $V_*(t)\in D(A)$ for
a.\ e.\ $t\in(0,T_*]$. The above claim yields
$$
v_*\in C((0,T_*];D(A))\cap L^2((0,T_*];D(A^{\frac32})),
$$
and so $v_*(t)\in D(A^{\frac32})$ for a.\ e.\ $t\in(0,T_*]$. It is
then sufficient to iterate this argument a finite number of times,
to deduce \eqref{e:addreg} for these values of $\gamma$.

Let us prove the framed claim above.
Again we just show the basic estimate, which follows easily
from inequality \eqref{e:Btemam}:
\begin{align*}
\frac{d}{dt}|A^{\frac{m}2}v_*|_H^2+2|A^{\frac{m+1}2}v_*|_H^2
&  =-\ps{A^mv_*}{B(v_*+z,v_*+z)}_H\\
&\le |A^{\frac{m+1}2}v_*|^2_H+C_m(1+|v_*|_V+|A^{\frac{m}2}z|_H)^{p_m},
\end{align*}
the proof of the claim is then easy and omitted.

\underline{\textsl{Step 4}.} The transition from $\gamma=\frac{n}2$
to any $\gamma>\frac34$ follows from a variation of the framed claim
given above. Assume that for given $\kappa\ge\frac12$, $\be\in[0,\frac12)$,
with $\kappa+\be>\frac34$, and $t_0\in(0,T_*]$ we have
$$
v_*\in C([t_0,T_*];D(A^\kappa))\cap L^2([t_0,T_*];D(A^{\kappa+\frac12}))
$$
and $z\in C([0,T_*];D(A^{\kappa+\beta}))$, then
$$
v_*\in C([t_0,T_*];D(A^{\kappa+\be}))\cap L^2([t_0,T_*];D(A^{\kappa+\be+\frac12})).
$$
This claim follows from the estimate
\begin{multline*}
\ps{A^{\kappa+\be}v}{A^{\kappa+\be}B(v+z,v+z)}_H
\le\frac12|A^{\kappa+\be+\frac12}v|_H^2+\wt{C}|A^{\frac{m}2+\be}z|_H^4+\\
\qquad+\wt{C}\Bigl(|A^{\kappa+\frac12}v|_H^2|A^{\kappa}v|_H^{\frac{4\be}{1-2\be}}+|A^{\kappa+\be}z|_H^2\Bigr)|A^{\kappa+\be}v|_H^2,
\end{multline*}
whose proof follows. First, by using Lemma \ref{l:Breg} with $\al=\kappa+\frac14$
and interpolation inequalities, we get
\begin{align*}
\ps{A^{\kappa+\be}v}{A^{\kappa+\be}B(v,v)}_H
&=   \ps{A^{\kappa+2\be}v}{A^{\kappa}B(v,v)}_H
 \le \wt{C} |A^{\kappa+2\be}v|_H|A^{\kappa+\frac12}v|^2_H\\
&\le \wt{C} |A^{\kappa+2\be}v|_H|A^{\kappa+\frac12}v|^{1-2\be}_H|A^{\kappa+\frac12}v|^{1+2\be}_H\\
&\le \wt{C} |A^{\kappa+\frac12}v|_H^{1-2\be}|A^{\kappa}v|_H^{2\be}|A^{\kappa+\be}v|_H^{1-2\be}
         |A^{\kappa+\frac12+\be}v|_H^{1+2\be}\\
&\le \frac18|A^{\kappa+\frac12+\be}v|_H^2+\wt{C}|A^{\kappa+\frac12}v|_H^2|A^{\kappa}v|_H^{\frac{4\be}{1-2\be}}|A^{\kappa+\be}v|_H^2.
\end{align*}
Similarly, by using Lemma \ref{l:Breg} with $\al=\kappa+\beta-\frac14$,
\begin{align*}
\ps{A^{\kappa+\be}v}{A^{\kappa+\be}B(v,z)}_H
&=   \ps{A^{\kappa+\be+\frac12}v}{A^{\kappa+\be-\frac12}B(v,z)}_H\\
&\le |A^{\kappa+\be+\frac12}v|_H|A^{\kappa+\be-\frac12}B(v,z)|_H\\
&\le C_B |A^{\kappa+\be+\frac12}v|_H\,|A^{\kappa+\be}v|_H\,|A^{\kappa+\be}z|_H\\
&\le \frac18|A^{\kappa+\be+\frac12}v|_H^2+C_B|A^{\kappa+\be}v|_H^2|A^{\kappa+\be}z|_H^2
\end{align*}
and in the same way,
$$
\ps{A^{\kappa+\be}v}{A^{\kappa+\be}B(z,z)}_H
\le \frac18|A^{\kappa+\be+\frac12}v|_H^2+C_B|A^{\kappa+\be}z|_H^4.
$$
All together, the above estimates yield the estimate of the non-linearity.
\end{proof}
\subsection{A few results on the regularised problem}\label{ss:aapprox}
In this appendix we give details of the proofs of Theorem \ref{t:weakstrong}
and Proposition \ref{p:BEL}. We denote by $u^\er{R}_x$ the solution to
equations \eqref{e:eqR} (possibly dropping the $x$ sub-script when there is
no ambiguity), and by $P_x^\er{R}$ its law on $\Omns$. Moreover, we write
$u^\er{R}=v^\er{R}+z$, where $v$ is the solution to the auxiliary problem
\begin{equation}\label{e:auxNS}
\dot v^\er{R}+Av^\er{R}+B(v^\er{R}+z,v^\er{R}+z)\chi_{R}(|v^\er{R}+z|_\W^2)=0,
\end{equation}
where for all $\om\in\Omns$, $z(\cdot,\om)$ solves the Stokes problem
\begin{equation}\label{e:stokes}
z(t)+\int_0^tAz(s)\,ds=\om(t).
\end{equation}
\begin{proof}[Proof of Theorem \ref{t:weakstrong}]
By Assumption \ref{a:noiseass2} the trajectories of the noise belong to 
$$
\Omns^*=\bigcap_{\beta\in(0,\frac12),\ \al\in[0,\al_0)} C^\be([0,\infty);D(A^\al)).
$$
with probability one, hence analyticity of the semi-group generated
by $A$ implies that, corresponding to each $\om\in\Omns^*$,
$z\in C([0,\infty);D(A^{\sss{\frac12}+\al_0-\ep}))$ for every $\ep>0$
(see Flandoli \cite{Fla}), and in particular $z\in C([0,\infty);\W)$,
since $\th(\al_0)<\al_0+\frac12$.

Given $\om\in\Omns^*$, one can prove that equation \eqref{e:auxNS} has
a unique global weak solution in the space $C([0,\infty);\W)$. The
proof of this claim can be carried on by means of standard arguments
such as Galerkin approximations (see for instance Flandoli \&
Gatarek \cite{FlGa}) and here is omitted for the sake of brevity.
Anyway, the crucial estimates \eqref{e:crucial1}, \eqref{e:crucial2}
and \eqref{e:crucial3} are given below.

Next, we prove \eqref{e:errore}. In order to do so, it is
sufficient to show that $P^\er{R}_x[\tau_R<\ep]\le C(\ep,R)$
with $C(\ep,R)\downarrow0$ as $\ep\downarrow0$, for all
$x\in\W$, with $|x|^2_\W\le\frac{R}{4}$,
and for $\ep$ small enough (depending only on $R$). Fix then
$\ep>0$ small enough, then, for $t<\tau_R(\om)$,
$$
\frac{d}{dt}|v^\er{R}|_\W^2+2|A^{\frac12}v^\er{R}|_\W^2
\le|\langle v^\er{R},B(v^\er{R}+z,v^\er{R}+z)\rangle_\W|.
$$
Consider first the case \framebox{$\alpha_0>$\text{\tiny$\frac12$}}.
Let $\Theta_{\ep,R}=\sup_{[0,\ep]}|z|_\W$ and
assume that $\Theta_{\ep,R}^2\le\frac{R}{4}$. Then, by Lemma \ref{l:Breg},
\begin{align}\label{e:crucial1}
|\langle v^\er{R},B(v^\er{R}+z,v^\er{R}+z)\rangle_\W|
&\le C|A^{\frac12}v^\er{R}|_\W|v^\er{R}+z|_\W^2\notag\\
&\le 2|A^{\frac12}v^\er{R}|_\W^2+C(|v^\er{R}|_\W^2+\Theta_{\ep,R}^2)^2,
\end{align}
and, if we set $\varphi(t)=|v^\er{R}(t)|_\W^2+\Theta_{\ep,R}^2$, we
get $\dot\varphi\le C\varphi^2$. This implies, together with the bounds
on $x$ and $\Theta_{\ep,R}$, that $\varphi(t)\le R(2-C\ep R)^{-1}$
and in conclusion
$$
|u^\er{R}(t)|^2
\le 2(|v^\er{R}(t)|_\W^2+|z(t)|_\W^2)
\le 2\varphi(t)
\le\frac{2R}{2-C\ep R}
\le R+1
$$
for $\ep\le\frac{2}{C R (R+1)}$. In particular, since this
holds for all $t\le\ep$, it follows that $\tau_R\ge\ep$.
Hence
$$
P_x^\er{R}[\tau_R<\ep]
\le P_x^\er{R}[\sup_{[0,\ep]}|z|_\W^2>\frac{R}{4}],
$$
and, since the last probability above is independent of $x$ (it depends only
on the law of the Stokes problem) and converges to $0$ as $\ep\downarrow0$,
the claim is proved.

In the special case \framebox{$\alpha_0=$\text{\tiny$\frac12$}}, we set
$\Theta_{\ep,R}=\sup_{[0,\ep]}|A^\oep z|_\W$, where
$\oep\in(0,\frac14)$ can be chosen arbitrarily small, and \eqref{e:crucial1}
is replaced by
\begin{align}\label{e:crucial2}
&|\langle v^\er{R},B(v^\er{R}+z,v^\er{R}+z)\rangle_\W|=\notag\\
&\qquad\qquad=   |\langle A^{\frac{5}{4}-\oep}v^\er{R},A^{\frac{1}{4}+\oep}B(v^\er{R}+z,v^\er{R}+z)\rangle_H|\notag\\
&\qquad\qquad\le C|A^{\frac{5}{4}-\oep}v^\er{R}|_H|A^{\frac{3}{4}+\oep}(v^\er{R}+z)|_H^2\notag\\
&\qquad\qquad\le C|v^\er{R}|_\W^{2\oep}|A^{\frac{1}{2}}v^\er{R}|_\W^{1-2\oep}\bigl(|v^\er{R}|_\W^{1-2\oep}|A^{\frac{1}{2}}v^\er{R}|_\W^{2\oep}+\Theta_{\ep,R}\bigr)^2\\
&\qquad\qquad\le |A^{\frac{1}{2}}v^\er{R}|_\W^2+C\bigl(|v^\er{R}|_\W^{4\frac{1-\oep}{1-2\oep}}+|v^\er{R}|_\W^{\frac{4\oep}{1+2\oep}}\Theta_{\ep,R}^{\frac{4}{1+2\oep}}\bigr)\notag\\
&\qquad\qquad\le |A^{\frac{1}{2}}v^\er{R}|_\W^2+C\bigl(|v^\er{R}|_\W^2+\Theta_{\ep,R}^{2(1-2\oep)}\bigr)^{2\frac{1-\oep}{1-2\oep}}\notag\\
&\qquad\qquad\le |A^{\frac{1}{2}}v^\er{R}|_\W^2+C\bigl(|v^\er{R}|_\W^2+\Theta_{\ep,R}^2+2\oep\bigr)^{2\frac{1-\oep}{1-2\oep}}\notag
\end{align}
which follows from Lemma \ref{l:Breg}, interpolation inequalities
and Young's inequality. We set $\varphi(t)=|v^\er{R}(t)|_\W^2+\Theta_{\ep,R}^2+2\oep$,
and the proof proceed as in the previous case to get $|u^\er{R}(t)|_\W^2\le R+1$ for
$\ep$ small enough and depending only on $R$.

If finally \framebox{$\alpha_0\in($\text{\tiny$\frac{1}{6},\frac{1}{2}$}$)$},
we set $\Theta_{\ep,R}=\sup_{[0,\ep]}|A^{\theta(\gamma)}z|$, where
$\gamma=\frac{\alpha_0}{2}+\frac{1}{4}$ (by this choice, in particular we
have $\theta(\alpha_0)<\theta(\gamma)<\alpha_0+\frac12$) and again by Lemma
\ref{l:Breg}, interpolation and Young's inequalities,
\begin{align}\label{e:crucial3}
|\langle v^\er{R},B(v^\er{R}+z,v^\er{R}+z)\rangle_\W|
&\le |A^{\frac{1}{2}}v^\er{R}|_\W|A^{-\frac{1}{2}}B(v^\er{R}+z,v^\er{R}+z)|_\W\notag\\
&\le C|A^{\frac{1}{2}}v^\er{R}|_\W|A^{\theta(\gamma)}(v^\er{R}+z)|^2\notag\\
&\le C|A^{\frac{1}{2}}v^\er{R}|_\W\bigl(|A^{\frac{1}{2}}v^\er{R}|_\W^{\frac{1-2\alpha_0}{4}}\!|v^\er{R}|_\W^{\frac{3+2\alpha_0}{4}}\!\!\!\!\!\!\!\!+\!\!\Theta_{\ep,R}\bigr)^2\notag\\
&\le |A^{\frac{1}{2}}v^\er{R}|_\W^2+C(|v^\er{R}|_\W^{2\frac{3+2\alpha_0}{1+2\alpha_0}}+\Theta_{\ep,R}^4)\\
&\le |A^{\frac{1}{2}}v^\er{R}|_\W^2+C\bigl(|v^\er{R}|_\W^2+\Theta_{\ep,R}^{4\frac{1+2\alpha_0}{3+2\alpha_0}}\bigr)^{\frac{3+2\alpha_0}{1+2\alpha_0}}\notag\\
&\le |A^{\frac{1}{2}}v^\er{R}|_\W^2+C\bigl(|v^\er{R}|_\W^2+\Theta_{\ep,R}^2+\frac{1}{3})^{1+a},\notag
\end{align}
where $a=2(1+2\alpha_0)^{-1}$. If $\varphi(t)=|v^\er{R}(t)|^2_\W+\Theta_{\ep,R}^2+\frac{1}{3}$
we get the same conclusions as in the previous cases.


Finally, in order to prove \eqref{e:weakstrong}, we show a
stronger property, namely that $\tau_R(u^\er{R}_x)=\tau_R(\xi)$
and that $\xi$ coincides with $u^\er{R}_x$ on $[0,\tau_R]$,
$P_x$-a.\ s.. Indeed, set $w=u^\er{R}_x-\xi$ and consider
the process $\widetilde{E}_t$ defined below in \eqref{e:Etilde}.
By Lemma \ref{l:Etilde}, $\widetilde{E}$ is an a.\ s.\ super-martingale.
Since by Lemma \ref{l:tauisstopping} $\tau_R$ is a stopping time,
Proposition \ref{p:stoppata} ensures that the stopped process
$(\widetilde{E}_t^{\tau_R})_{t\ge0}$ is an a.\ s.\ super-martingale too,
and in particular
$\E^{P_x}[\widetilde{E}_t^{\tau_R}]\le\E^{P_x}[\widetilde{E}_0^{\tau_R}]=0$.

Now, for every $s\le t\wedge\tau_R$ we have $|u^\er{R}(s)|_\W^2\le R$
and so $\chi_R(\|u^\er{R}(s)|_\W)=1$. Moreover, using the
properties of Navier-Stokes non-linearity,
$$
B(u^\er{R},u^\er{R})-B(\xi,\xi)
=B(w,u^\er{R})+B(\xi,w)
$$
and
$$
\langle w,B(u^\er{R},u^\er{R})-B(\xi,\xi)\rangle
=\langle w,B(w,u^\er{R})+B(\xi,w)\rangle
=\langle w,B(w,u^\er{R})\rangle
$$
The last term in the equality above can be estimated using
Lemma {\S}I.2.1 of Temam \cite{Tem2}. Indeed, if $\al<1$,
\begin{align*}
\langle w,B(w,u^\er{R})\rangle
&\le C|A^{\frac58-\frac{\th(\al)}{2}}w|^2_H |u|_{\mathcal{W}}\\
&\le C(R)|A^{\frac58-\frac{\th(\al)}{2}}w|^2_H\\
&\le C(R)|w|_V^{\frac52-2\th(\al)}|w|_H^{2\th(\al)-\frac12}\\
&\le |w|_V^2+C(R,\al)|w|_H^2,
\end{align*}
by interpolation, since $\sss{\frac58-\frac{\th(\al)}{2}<\frac12}$ for all $\al>0$.
If $\al>1$, the same inequality holds by simply replacing
$|A^{\sss{\frac58-\frac{\th(\al)}{2}}}w|^2_H$ with $|w|_H^2$.
The special case $\al=1$ (corresponding to the Sobolev critical
exponent) can be handled by simply writing the same inequality
for a slightly smaller $\al=1-\ep$.

Next, consider $t>0$, then we know that $\E^{P_x}[\widetilde{E}_t^{\tau_R}]\le0$, and so,
using the previous inequality,
$$
\E^{P_x}[|w(t\wedge\tau_R)|_H^2]
\le C_{\al,R}\E^{P_x}\Bigl[\int_0^{t\wedge\tau_R}|w(s)|_H^2\,ds\Bigr]
\le C_{\alpha,R}\int_0^t\E^{P_x}[|w(s\wedge\tau_R)|_H^2].
$$
By Gronwall's lemma we finally deduce that $\E[|w(t\wedge\tau_R)|_H^2]=0$,
for all $t\ge0$. We next prove that 
\begin{equation}\label{e:claim}
\E^{P_x}[|w(\tau_R(u^\er{R}))|_H^2]=0.
\end{equation}
On $\{\tau_R(u^\er{R})=+\infty\}$ the above formula
is obvious, so we need to prove that the above
expectation is zero on $\{\tau_R(u^\er{R})<\infty\}$.
If $\om\in\{\tau_R(u^\er{R})<\infty\}$ and $t\uparrow\infty$,
then $t\wedge\tau_R(u^\er{R}(\om))\uparrow\tau_R(u^\er{R}(\om))$
and so, by semi-continuity (we recall that $w$ is continuous in
time with respect to the weak topology of $H$)
$$
|w((\tau_R(u^\er{R}(\om)))|_H^2\le\liminf_{t\uparrow\infty}|w(t\wedge\tau_R(u^\er{R}(\om)))|_H^2.
$$
Fatou's lemma then implies \eqref{e:claim}. Finally, from formula
\eqref{e:claim} and the above considerations we deduce that
$\tau(u_x^\er{R})=\tau_R(\xi)$, $P_x$-a.\ s.\ and that
$u_x^\er{R}=\xi$, $P_x$-a.\ s.\ on $[0,\tau_R]$.
\end{proof}
The previous proof used the following two lemmas. The first shows
that the \emph{blow-up} time for the regular solution is a stopping
time with respect to the natural filtration, the second prove
an a.\ s.\ super-martingale property for the difference between
the weak and the strong solution.
\begin{lemma}\label{l:tauisstopping}
Let $P_x$ be a martingale solutions to Navier-Stokes equations \eqref{e:nseq}.
Under the assumptions of Theorem \ref{t:weakstrong}, consider the process $u^\er{R}_x$
and the random time $\tau_R$ defined in \eqref{e:stoppingtime}. Then $\tau_R(u_x^\er{R})$
is a stopping time with respect to the filtration $(\BNS_t)_{t\ge0}$.
\end{lemma}
\begin{proof}
Let $\tau_R'$ be the time $\tau_R$ on $\Omns'=C([0,\infty);\W)$.
By Proposition $4.5$ of Revuz \& Yor \cite{ReYo}, $\tau_R'$ is
a stopping time with respect to the canonical filtration of
$\Omns'$. On the other hand,
$\tau_R'(u_x^\er{R}(\om))=\tau_R(u_x^\er{R}(\om))$
for $P_x$-a.\ e.\ $\om\in\Omns$, and therefore
$$
\{\om\in\Omns:\tau_R(u_x^\er{R}(\om))\le t\}
=\{\om\in\Omns:u_x^\er{R}(\om)\in\{\om'\in\Omns':\tau_R'(\om')\le t\}\}.
$$
The set on right-hand side is finally $\BNS_t$-measurable since
$u_x^\er{R}$ is progressively measurable.
\end{proof}
\begin{lemma}\label{l:Etilde}
Let $P_x$ be a martingale solutions to Navier-Stokes equations \eqref{e:nseq}.
Under the assumptions of Theorem \ref{t:weakstrong}, set $w=u^\er{R}_x-\xi$.
and define
\begin{equation}\label{e:Etilde}
\widetilde{E}_t
=|w(t)|_H^2
 +2\int_0^t|w(s)|_V^2\,ds
 +2\int_0^t\langle\chi_R(\|u^\er{R}\|_{\W}^2)B(u^\er{R},u^\er{R})-B(\xi,\xi),w\rangle.
\end{equation}
Then the process $(\widetilde{E}_t,\BNS_t,P_x)_{t\ge0}$ is an a.\ s.\ super-martingale.
\end{lemma}
\begin{proof}
Since $\xi_t$ already satisfies an energy inequality (property \map{3} of Definition
\ref{d:asNSsol}), if we write $|w|_H^2=|u^\er{R}|_H^2-2\langle u^\er{R},\xi\rangle_H+2|\xi|_H^2$,
then we only need to show energy equalities (in terms of super-martingales) for
$|u^\er{R}|_H^2$ and $\langle u^\er{R},\xi\rangle_H+2|\xi|_H^2$. We actually show
that for every $s\not\in T_{P_x}$ (where $T_{P_x}$ is the set of exceptional times
of $P_x$) and every $t>s$,
\begin{align*}
\langle u^\er{R}(t),\xi(t)\rangle
&=\langle u^\er{R}(s),\xi(s)\rangle
 -\int_s^t\langle Au^\er{R}+\chi_R(\|u^\er{R}\|_{\W}^2)B(u^\er{R},u^\er{R}),\xi\rangle\,dr+\\
&\quad-\int_s^t\langle A\xi+B(\xi,\xi),u^\er{R}\rangle\,dr
      +\int_s^t\langle\mathcal{Q}^{\frac12}\,dW(r),\xi\rangle+\\
&\quad+\int_s^t\langle\mathcal{Q}^{\frac12}\,dW(r),u^\er{R}\rangle
      +\sigma^2 t
\end{align*}
and
$$
|u^\er{R}(t)|_H^2
+2\int_s^t|u^\er{R}|_V^2\,dr
=|u^\er{R}(s)|_H^2+\int_s^t\langle u^\er{R},\mathcal{Q}^{\frac12}\,dW\rangle+\sigma^2(t-s).
$$
hold $P_x$-a.\ s. (in particular, all the quantities in the two formulae above
are well defined and finite).

The proof is in the spirit of the results of Serrin \cite{Ser} on weak-strong
uniqueness: the process $u^\er{R}$ is regular enough (more precisely,
$\sup|u^\er{R}|_\W$ has exponential moments) so that the estimates can
be obtained with standard arguments. We prove the first of the two equalities
above, the other's being entirely similar. Let $\pi_n$ be the projection
on the first $n$ Fourier modes and let $\xi_n(t)=\pi_n\xi(t)$ and
$u_n^\er{R}(t)=\pi_n u^\er{R}(t)$, then
\begin{align*}
u_n^\er{R}(t)&=u_n^\er{R}(0)-\int_0^t\pi_n(Au^\er{R}+\chi_R(|u^\er{R}|_\W)B(u^\er{R},u^\er{R}))\,ds+\pi_n\Q^{\frac12}W(t),\\
\xi_n(t)&=\xi_n(0)-\int_0^t\pi_n(A\xi+B(\xi,\xi))\,ds+\pi_n\Q^{\frac12}W(t),
\end{align*}
so that, by standard arguments,
\begin{align*}
&\langle u_n^\er{R}(t),\xi_n(t)\rangle-\langle u_n^\er{R}(0),\xi_n(0)\rangle=\\
&=-\int_0^t\langle\pi_n(Au+\chi_R(|u^\er{R}|_\W)B(u^\er{R},u^\er{R})),\xi_n\rangle
  -\int_0^t\langle u_n^\er{R},\pi_n(A\xi+B(\xi,\xi))\rangle\\
&\quad +\int_0^t\langle\pi_n\Q^{\frac12}\,dW(s),\xi_n\rangle
       +\int_0^t\langle u_n^\er{R},\pi_n\Q^{\frac12}\,dW(s)\rangle
       +t\Tr(\pi_n\Q\pi_n).
\end{align*}
The terms on the left-hand side above converge to the corresponding term
since $u^\er{R}$ and $\xi$ are continuous in time with values in $H$ with
respect to respectively the strong and the weak topology. The limit in
the It\^{o} integrals and the correcting term are equally easy.
As for the two Lebesgue integrals, the scalar products converge
by Lebesgue theorem, since for both weak and strong solutions, a.\ s.\ in time,
$Au^\er{R}+\chi_R(|u^\er{R}|_\W)B(u^\er{R},u^\er{R})$ and $A\xi+B(\xi,\xi)$
are in $V'$ and $u^\er{R}$, so that the integrands converge
a.\ s.\ in time, and since the following bounds
\begin{align*}
|\langle\pi_n(Au^\er{R}+\chi_R(|u^\er{R}|_\W)B(u^\er{R},u^\er{R})),\xi_n\rangle_H|
&\le C|\xi|_V(|u^\er{R}|_V+|u^\er{R}|_{L^4}^2)\\
&\le C|\xi|_V\sup_{[0,T]}|u^\er{R}|_V^2
\end{align*}
and
\begin{align*}
|\langle u_n^\er{R},\pi_n(A\xi+B(\xi,\xi))\rangle_H|
&\le C|u^\er{R}|_V(|\xi|_V+|\xi|_{L^4}^2)\\
&\le C(|\xi|_V+|\xi|_H^{\frac12}|\xi|_V^{\frac32})\sup_{[0,T]}|u^\er{R}|_V,
\end{align*}
hold by standard inequalities (the terms on the right-hand side of both
equations are integrable $P$-a.\ s.). The proof is complete. 
\end{proof}
\begin{proof}[Proof of Proposition \ref{p:BEL}]
We aim to show that the transition semi-group
$(\Ps_t^\er{R})_{t\ge0}$ is $\W$--strong Feller. As in
the proof of the previous lemma, we shall provide
formal estimates, that can be made rigorous
only at the level of Galerkin approximations.
Let $(\Sigma, \F, (\F_t)_{t\ge0},\abP)$ be
a filtered probability space, $(W_t)_{t\ge0}$
the cylindrical Wiener process on $H$ and,
for every $x\in\W$, let $u_x^\er{R}$ be the
solution to equations \eqref{e:eqR}. By
the Bismut, Elworthy \& Li formula,
$$
D_y(\Ps_t^\er{R}\psi)(x)
=\frac1{t}\E^{\abP}[\psi(u_x^\er{R}(t))\int_0^t\ps{\Q^{-\frac12}D_y u_x^\er{R}(s)}{dW(s)}]
$$
and thus, for $|\psi|_\infty\leq1$, by Burkholder, Davis \& Gundy inequality,
\begin{multline*}
|(\Ps_t^\er{R}\psi)(x_0+h)-(\Ps_t^\er{R}\psi)(x_0)|\le\\
\le\frac{C}{t}\sup_{\eta\in[0,1]}
\E^{\abP}[(\int_0^t|Q_0^{-\frac12}A^{\frac34+\al_0}D_h u_{x_0+\eta h}^\er{R}(s)|_H^2\,ds)^{\frac12}].
\end{multline*}
The proposition is proved once we can prove that the right-hand side
of the above inequality converges to $0$ as $|h|_\W\to0$.

Assume first that $\al_0\neq\frac12$. Fix $x\in\W$ and
write $u=u_x^\er{R}$. The term $Du$ solves the following equation
\begin{align*}
\frac{d}{dt}(A^{\frac14+\al_0}Du)+A(A^{\frac14+\al_0}Du)
&=\chi_R(|u|_\W^2)A^{\frac14+\al_0}(B(Du,u)+B(u,Du))\\
&\quad+2\chi_R'(|u|_\W^2){\ps{u}{Du}}_\W\,A^{\frac14+\al_0}B(u,u),
\end{align*}
hence
\begin{align}\label{e:stimaDu}
&\frac{d}{dt}|A^{\al_0+\frac14}Du|_H^2+2|A^{\al_0+\frac34}Du|_H^2=\notag\\
&\qquad\quad =2\chi_R(|u|_\W^2)\ps{A^{\al_0+\frac14}Du}{A^{\al_0+\frac14}(B(Du,u)+B(u,Du))}+\\
&\qquad\qquad+4\chi_R'(|u|_\W^2){\ps{u}{Du}}_\W\ps{A^{\al_0+\frac14}Du}{A^{\al_0+\frac14}B(u,u)}.\notag
\end{align}
We use Lemma \ref{l:Breg} to estimate the terms on the right-hand side.
\begin{align*}
&2\chi_R(|u|_\W^2)\ps{A^{\al_0+\frac14}Du}{A^{\al_0+\frac14}(B(Du,u)+B(u,Du))}=\\
&\qquad=  2\chi_R(|u|_\W^2)\ps{A^{\al_0+\frac34}Du}{A^{\al_0-\frac14}(B(Du,u)+B(u,Du))}\\
&\qquad\le C(B)\chi_R(|u|_\W^2)|A^{\al_0+\frac34}Du|_H\,|Du|_\W\,|u|_\W\\
&\qquad\le\frac12|A^{\al_0+\frac34}Du|_H^2+C(B,R)|Du|^2_\W,
\end{align*}
and, similarly,
\begin{align*}
&4\chi_R'(|u|_\W^2){\ps{u}{Du}}_\W\ps{A^{\al_0+\frac14}Du}{A^{\al_0+\frac14}B(u,u)}\le\\
&\qquad\le C(B)\chi_R'(|u|_\W^2)|u|_\W^3|Du|_\W\,|A^{\al_0+\frac34}Du|_H\\
&\qquad\le\frac14|A^{\al_0+\frac34}Du|_H^2+C(B,R)|Du|^2_\W.
\end{align*}
We plug the above estimates into \eqref{e:stimaDu} and we use
an interpolation argument to get
$$
\frac{d}{dt}|A^{\al_0+\frac14}Du|_H^2+|A^{\al_0+\frac34}Du|_H^2
\le C(B,R,\al_0)|A^{\al_0+\frac14}Du|_H^2,
$$
since, by definition, $\th(\al_0)\in[\al_0+\frac14,\al_0+\frac34)$. From
Gronwall's inequality we finally obtain
\begin{equation}\label{e:pippo}
\E^{\abP}\int_0^T|A^{\al_0+\frac34}Du(s)|_H^2\,ds\le C(B,R,\al_0,T)|A^{\al_0+\frac14}h|_H^2
\end{equation}
and the proof for $\al_0\neq\frac12$ is completed.

The case $\al_0=\frac12$ is slightly more complicated, since for
this value of $\al_0$ Lemma \ref{l:Breg} is weaker. As in estimate
\eqref{e:crucial2}, we introduce a small penalisation
exponent $\oep$ and, by proceeding as in the case $\al_0\neq\frac12$
we get
\begin{align*}
&2\chi_R(|u|_\W^2)\ps{A^{\frac34}Du}{A^{\frac34}(B(Du,u)+B(u,Du))}\le\\
&\qquad\qquad\le C(B)\chi_R|A^{\frac54-\oep}Du|\cdot|A^{\frac34+\oep}Du|\cdot|A^{\frac34+\oep}u|\\
&\qquad\qquad\le\frac12|A^{\frac54}Du|^2+C(B,\oep)\chi_R|A^{\frac34+\oep}u|^2|A^{\frac34}Du|^2
\end{align*}
and
\begin{align*}
&4\chi_R'(|u|_\W^2){\ps{u}{Du}}_\W\ps{A^{\al_0+\frac14}Du}{A^{\al_0+\frac14}B(u,u)}\le\\
&\qquad\qquad\le C(B,R)\chi_R|A^\frac34Du|\cdot|A^{\frac54-\oep}Du|\cdot|A^{\frac34-\oep}u|^2\\
&\qquad\qquad\le\frac12|A^\frac54Du|^2+C(B,R,\oep)|A^{\frac34+\oep}u|^{\frac4{1+2\oep}}|A^\frac34Du|^2.
\end{align*}
Now, since by interpolation $\chi_R |A^{\frac34+\oep}u|\le C(R)|A^{1-\ep}u|^{e(\oep,\ep)}$, where
$\ep$ is suitably chosen and $e(\oep,\ep)<<1$ (when $\oep<<1$ and $\ep<<1$), we finally get
$$
\frac{d}{dt}|A^{\frac34}Du|_H^2+|A^{\frac54}Du|_H^2
\le C(B,R,\oep,\ep)(1+|A^{1-\ep}u|^{e(\oep,\ep)})|A^{\frac34}Du|_H^2,
$$
and \eqref{e:pippo} follows for this case by Gronwall's lemma, provided that
\begin{equation}\label{e:fernique}
\E^{\abP}\Bigl[\exp\Bigl(C\int_0^T|A^{1-\ep}u|^e\,ds\Bigr)\Bigr]
\end{equation}
is finite. If we split $u=v+z$ as in the previous proof, the quantity \eqref{e:fernique}
above for $z$ is finite by Fernique's theorem (see also Proposition 2.16, Da Prato \& Zabczyk
\cite{DPZa2}). As it concerns the same quantity for $v$, we need a few additional computations.
Indeed, it is sufficient to find a small $\be>0$ such that $\E^\abP[F(T)]$ is finite for the
quantity $F(t)=\exp(\be|A^{\frac34}v(t)|^2+\be\int_0^T|A^{\frac54}v|_H^2\,ds)$. Now,
$$
\dot F(t)
=F(t)[-\be|A^\frac54 v|^2+2\be\chi_R\langle A^{\frac54-\oep}v,A^{\frac14+\oep}B(u,u)\rangle]
$$
and, by using again Lemma \ref{l:Breg} and interpolation
\begin{align*}
\chi_R\langle A^{\frac54-\oep}v,A^{\frac14+\oep}B(u,u)\rangle
&\le C(B)\chi_R|A^{\frac54-\oep}v|\cdot|A^{\frac34+\oep}u|^2\\
&\le C(B,R,\oep,\ep)|A^\frac54 v|(\Theta^\frac12+|A^{1-\ep}v^{\frac{8\oep}{1-4\ep}})\\
&\le\frac12|A^\frac54 v|^2+C(B,R,\oep,\ep)(1+\Theta),
\end{align*}
where $\Theta=\sup_{[0,T]}|A^{1-\ep}z|^{\frac{16\oep}{1-4\ep}}$. In conclusion
we get $\dot F\le \beta\,C(B,R,\oep,\ep)(1+\Theta)\,F(t)$ and the proof is complete.
\end{proof}
\subsection{Controllability results for the regularised problem}\label{ss:control}
In this part we show the auxiliary results used in the proof of
Proposition \ref{p:support}. In the following lemmas we follow
closely Lemma $2.1$ and Lemma $3.1$ respectively, of Flandoli
\cite{Fla3}.
\begin{lemma}[Approximate controllability]\label{l:appcon}
Let $R>0$, $T>0$ and let $x\in\W$ and $y\in\W$, with $Ay\in\W$, such that
$$
|x|^2_\W\le\frac{R}2
\qquad\text{and}\qquad
|y|^2_\W\le\frac{R}2,
$$
then there exist $w\in\text{Lip}([0,T];\W)$ and
$$
u\in C([0,T];\W)\cap L^2([0,T];D(A^{\th(\al_0)+\sss{\frac12}})),
$$
such that $u$ solves the equation
\begin{equation}\label{e:ceq}
\partial_t u+Au+\chi_R(|u|_\W^2)B(u,u)=\partial_t w,
\end{equation}
with $u(0)=x$ and $u(T)=y$, and
\begin{equation}\label{e:cbd}
\sup_{t\in[0,T]}|u(t)|_\W^2\le R.
\end{equation}
\end{lemma}
\begin{proof}
Consider first $w=0$, then by an inequality similar to \eqref{e:crucial1}--\eqref{e:crucial3},
we get
$$
\frac{d}{dt}|u|_\W^2+|A^{\frac12}u|_\W^2\le C(B,\al_0)(R+2)^e|u|_\W^2,
$$
for some exponent $e$ and some constant $C(B,\al_0)$ depending on $\al_0$ and
the non-linearity, so that by Gronwall's lemma,
$$
|u(t)|_\W^2+\int_0^t|A^{\frac12}u|_\W^2\,ds\le\frac{R}2\e^{tC(B,\al_0,R)}.
$$
Hence, $u\in D(A^{\th(\al_0)+\sss{\frac12}})$ almost everywhere and, by starting again
the equation on one of this regular points (notice that the equation has unique solution),
we can find a small $T_*\in(0,\frac{T}2)$ such that
\begin{align*}
&|u(t)|_\W^2\le R,\qquad\text{for all }t\le T_*,\\
&Au(T_*)\in\W.
\end{align*}
Define $u$ to be the solution above for $t\in[0,T_*]$ and set for $t\in[T_*,T]$,
$$
u(t)=\frac{T-t}{T-T_*}u(T_*)+\frac{t-T_*}{T-T_*}y.
$$
First, we obviously have \eqref{e:cbd}. Next, if we set
$$
\eta=\partial_t u+Au+\chi_R(|u|_\W^2)B(u,u)
$$
and $w$ to be $0$ for $t\le T_*$ and $w(t)=\int^t_{T_*}\eta_s\,ds$ for $t\in[T_*,T]$,
we also have \eqref{e:ceq}. We only have to prove that $\eta\in L^\infty(0,T;\W)$.
The first two terms of $\eta$ are obvious, for the non-linear term we observe that
by Lemma \ref{l:Breg} it follows that
$$
|B(u_1,u_2)|_\W\le C|u_1|_{D(A^{\th(\al_0)+\sss{\frac12}})}|u_2|_{D(A^{\th(\al_0)+\sss{\frac12}})},
$$
for any $u_1$, $u_2\in D(A^{\th(\al_0)+\sss{\frac12}})$.
\end{proof}
Let $s\in(0,\frac12)$ and $p>1$ such that $s-\sss{\frac1p}>0$, under this assumption
one can see (see Flandoli \cite{Fla3}) that for every $\al_1<\al_0$ the map
$$
w\mapsto z(\cdot,w):
W^{s,p}([0,T];D(A^{\al_1}))
\longrightarrow
C([0,T];D(A^{\al_1+s-\frac1p-\ep}))
$$
is continuous, for all $\ep>0$, where $z$ is the solution to the
Stokes problem we have introduced in the proof of Theorem \ref{t:weakstrong}.
In particular, it is possible to find, for any value of $\al_0$,
values $\al_1\in(0,\al_0)$, $s$ and $p$ such that the above map is
continuous from $W^{s,p}(0,T;D(A^{\al_1}))$ with values in
$C([0,T];D(A^{\oth}))$, where $\oth$ corresponds to the space of
regularity where $z$ is evaluated in inequalities \eqref{e:crucial1},
\eqref{e:crucial2} and \eqref{e:crucial3}, hence $\oth=\th(\al_0)$
for $\al_0>\frac12$, $\oth=\th(\al_0)+\oep$ for $\al_0=\frac12$ and
$\oth=\th(\frac{\al_0}2+\frac14)$ for $\al_0\in(\frac16,\frac12)$.
\begin{lemma}[Continuity along the controllers]\label{l:concon}
Let $s$, $p$ and $\al_1$ be chosen as above, and let $w_n\to w$
in $W^{s,p}([0,T];D(A^{\al_1}))$. Let $u$ be the solution to
equation \eqref{e:ceq} corresponding to $w$ and some initial
condition $x$ and let
$$
\tau=\inf\{t\ge0\,:\,|u(t)|_\W^2\ge R\}
$$
(and $\tau=+\infty$ if the set is empty). Define similarly, for each $n\in\N$, $u_n$
and $\tau_n$ corresponding to $w_n$ and the same initial condition $x$.

If $\tau>T$, then $\tau_n>T$ for $n$ large enough and
$$
u_n\longrightarrow u
\qquad\text{in }C([0,T];\W).
$$
\end{lemma}
\begin{proof}
Set $v_n=u_n-z_n$ for each $n\in\N$, and $v=u-z$, where $z_n$, $z$
are the solutions to the Stokes problem corresponding to $w_n$, $w$ respectively
(see the proof of Theorem \ref{t:weakstrong}). Since $w_n\to w$ in $W^{s,p}(0,T;D(A^{\al_1}))$,
this gives a common lower bound for $(\tau_n)_{n\in\N}$ and $\tau$. For every time
smaller than this lower bound, we can estimate,
\begin{align*}
\frac{d}{dt}|v-v_n|_\W^2+2|A^{\frac12}(v_n-v)|^2_\W
&=2\ps{v_n-v}{B(u,u)-B(u_n,u_n)}_\W\\
&=2\langle{v_n-v},B(u_n,v-v_n)+B(u_n,z-z_n)+\\
&\quad+B(v-v_n,u)+B(z-z_n,u)\rangle_\W.
\end{align*}
In order to estimate the term on the right-hand side above, we consider three
cases, depending on the value of $\al_0$. If $\al_0>\frac12$, we can estimate
the right-hand side of the above formula as in inequality \eqref{e:crucial1}.
For instance,
\begin{align*}
\langle v_n-v,B(u_n,v-v_n)\rangle_\W
&\le C|A^{\frac12}(v-v_n)|_\W|u_n|_\W|v-v_n|_\W\\
&\le \frac18|A^{\frac12}(v-v_n)|_\W^2+C|u_n|_\W^2|v-v_n|_\W^2,
\end{align*}
and the other pieces can be handled similarly. So, by Gronwall's lemma,
$$
|v_n(t)-v(t)|_\W^2
\le\Theta_n\exp(\int_0^t(|u_n|_\W^2+|u|_\W^2)\,ds)\int_0^t(|u_n|_\W^2+|u|_\W^2)\,ds
$$
where $\Theta_n=\sup_{[0,T]}|z-z_n|_\W^2$. Now, since $\tau>T$, if
$S=\sup_{t\in[0,T]}|u(t)|^2_\W$, then $S<R$ and we can find a $\delta>0$
(depending only on $R$ and $S$) and $n_0\in\N$ such that $\Theta_n<\delta$ and $|v_n-v|_\W<\delta$
for all $n\ge n_0$, and so
$$
|u_n(t)|_\W\le|v_n(t)-v(t)|_\W+\Theta_n^\frac12+|u(t)|_\W\le2\sqrt\delta+\sqrt{S}\le\sqrt{R-\delta}.
$$
We can conclude that $u_n\to u$ in $C([0,T];\W)$ and $\tau_n>T$ for all $n\ge n_0$.

For the other values of $\al_0$, we proceed similarly, using inequalities \eqref{e:crucial2}
and \eqref{e:crucial3}. If $\al_0=\frac12$, then $\Theta_n=\sup_{[0,T]}|A^{\frac34+\oep}(z-z_n)|^2$
and one gets
$$
\frac{d}{dt}|v_n(t)-v(t)|_\W^2
\le C(1+|A^{\frac34+\oep}u_n|^2+|A^{\frac34+\oep}u_n|^2)(|v-v_n|_\W^2+\Theta_n).
$$
Finally, if $\al_0\in(\frac16,\frac12)$, then we set $\Theta_n=\sup_{[0,T]}|A^{\oth}(z_n-z)|^2$,
with $\oth=\th(\frac{\al_0}2+\frac14)$, and we get
$$
\frac{d}{dt}|v_n(t)-v(t)|_\W^2
\le C(1+|u_n|_\W^2|A^{\frac12}u_n|_\W^a+|u|_\W^2|A^{\frac12}u|_\W^a)(|v-v_n|_\W^2+\Theta_n),
$$
with $a=\frac{2-4\al_0}{3+2\al_0}$.
\end{proof}
\section{Estimates on the non-linearity}\label{s:nonlinear}
We show two continuity result on the Navier-Stokes non-linearity.
The first one derives from Lemma 4.1 of Temam \cite{Tem2}. The
second is a proof of continuity of the bi-linear term in spaces
of powers of the Stokes operator.
\begin{lemma}
Let $m\ge2$ be an integer and let $v\in D(A^{\sss{\frac{m+1}2}})$ and
$z\in D(A^{\sss{\frac{m}2}})$. Then
\begin{equation}\label{e:Btemam}
|\ps{A^mv}{B(v+z,v+z)}_H|\le\frac12|A^{\frac{m+1}2}v|^2_H+C_m(1+|v|_V+|A^{\frac{m}2}z|_H)^{p_m},
\end{equation}
where $C_m$, $p_m>0$ are positive constants depending only on $m$.
\end{lemma}
\begin{proof}
The result is a variant of Lemma 4.1 of Temam \cite{Tem2}, to whom we shall
rely heavily. Consider first the term $\ps{A^mv}{B(w,v)}_H$, with $w$ equal to $v$
or $z$, it is given by terms of the form (integration by parts is used to produce
the second term)
$$
\int_\T w_i(D_{x_i}v_j)(D^{2\al}v_j)\,dx
=\int_\T(D^\al v_j)(D^\al w_i D_{x_i}v_j)\,dx,
$$
where $\al=(\al_1,\al_2,\al_3)$ is a multi-index with $\al_1+\al_2+\al_3=m$
and, as usual, $D^\al=D^{\al_1}_{x_1}D^{\al_2}_{x_2}D^{\al_3}_{x_3}$.
By expanding the derivatives of the product, we obtain terms of the form
$$
\int_\T w_i(D_{x_i}D^\al v_j)(D^\al v_j)\,dx
\qquad\text{and}\qquad
\int_\T(\delta^k w_i)(\delta^{m-k+1}v_j)(\delta^m v_j)\,dx,
$$
where $\delta^k$ denotes a generic differential operator of order $k$ and
$k=1,\dots,m$. The terms on the left are equal to zero, due to the
divergence-free constraint, while the terms on the right can be estimated
as in Lemma 4.1 of Temam \cite{Tem2} (H\"older inequality, plus the Sobolev
embeddings $L^3\subset H^{\sss{\frac12}}(\T)$ and $L^6\subset H^1(\T)$, plus
interpolation),
$$
\int_\T(\delta^k w_i)(\delta^{m-k+1}v_j)(\delta^m v_j)\,dx
\le |v|_V^{\frac{k}{m}}|A^{\frac{m+1}2}v|_H^{\frac{2m-k}{m}}|A^{\frac{2k+1}{4}}w|_H.
$$
If $w=v$, interpolation yields the bound
$|v|_V^{1+\sss{\frac{1}{2m}}}|A^\sss{\frac{m+1}2}v|_H^{2-\sss{\frac1{2m}}}$
and so \eqref{e:Btemam} by Young's inequality. If $w=z$, the estimate above
again leads to \eqref{e:Btemam}, since $\frac{2k+1}{4}\le\frac{m}{2}$,
unless $k=m$. For $k=m$, we use the same estimates as above, in a different
order,
\begin{align*}
\int_\T(\delta^m z)(\delta v_j)(\delta^m v_j)\,dx
&\le |A^{\frac{m}{2}}z|_H|\delta v_j|_{L^6}|\delta^m v_j|_{L^3}
 \le |A^{\frac{m}{2}}z|_H|Av|_H|A^{\frac{2m+1}{4}}v|_H\\
&\le |v|_V^{1-\frac1{2m}}|A^{\frac{m+1}2}|_H^{1+\frac1{2m}}|A^{\frac{m}{2}}z|_H.
\end{align*}
Consider next the terms $\ps{A^mv}{B(w,z)}_H$, with $w$ equal to $v$ or $z$.
By integration by parts,
\begin{align*}
\ps{A^mv}{B(w,z)}_H
&=-\ps{z}{B(w,A^mv)}
 =-\sum\int_\T z_j w_i D^{2\al}D_{x_i}v_j\\
&=\sum\int_\T(\delta^{m+1}v_j)(\delta^kw_i)(\delta^{m-k}z_j),
\end{align*}
where, as above, the sum is over multi-indexes $\al$ with $|\al|=m$
and $k=0,\dots,m$. In order to estimate the generic integral
$\int(\delta^{m+1}v)(\delta^kw)(\delta^{m-k}z)$, we use H\"older
inequality with exponents $2$, $3$ and $6$ respectively for $k>0$,
and exponents $2$, $+\infty$ and $2$ for $k=0$. In the case $k=0$
we use the embedding $D(A)\subset C(\T)$.

If $w=v$, we use interpolation inequalities and Young's inequality,
as above, to obtain \eqref{e:Btemam}, while for $w=z$ we only need
Young's, but for the case $k=m$, where again we use the embedding
$D(A)\subset C(\T)$,
$$
\int(\delta^{m+1}v)(\delta^mw)(\delta^kz)
\le |A^{\frac{m+1}2}v|_H\,|A^{\frac{m}2}z|_H\,|Az|_H
\le |A^{\frac{m+1}2}v|_H\,|A^{\frac{m}2}z|_H^2,
$$
and again Young's inequality.
\end{proof}
\begin{lemma}\label{l:Breg}
For each $\al>0$, with $\al\neq\frac12$, the bi-linear operator $B$ maps
$D(A^{\sss{\th(\al)}})\times D(A^{\sss{\th(\al)}})$ continuously to
$D(A^{\al-\sss{\frac14}})$, where $\th(\al)$ is the function defined
in \eqref{e:theta}.

If $\al=\frac12$, $B$ maps $D(A^{\sss{\frac34}})\times D(A^{\sss{\frac34}})$ continuously
to $D(A^{\frac14-\ep})$, for every $\ep>0$.
\end{lemma}
\begin{proof}
If $\al\le\frac14$, the lemma follows by Sobolev embeddings, as in the proof
of Lemma 2.1 of Temam \cite[Part I]{Tem2}. Assume that $\al>\frac14$ and
fix $u$, $v\in D(A^{\sss{\th(\al)}})$. If $(u_\ka)_{\ka\in\Z^3}$ and
$(v_\ka)_{\ka\in\Z^3}$ are the Fourier coefficients of $u$ and $v$,
the norm of $B(u,v)$ in $D(A^{\al-\sss{\frac14}})$ in terms of Fourier coefficients is given by
$$
\|B(u,v)\|_{\al-\frac14}^2
=  \sum_{\ka\in\Z^3_*}|\ka|^{4\al-1}\Bigl|\sum_{\el+\ma=\ka}(u_\el\cdot\ma)(v_\ma-\frac{v_\ma\cdot\ka}{|\ka|^2}\ka)\Bigr|^2.
$$
We split the inner sum in three terms, corresponding to the following three subsets,
$$
\begin{array}{rcl}
E_\ka&=&\{(\el,\ma)\in\Z^3_*\times\Z^3_*\,:\,\el+\ma=\ka,\ |\el|\ge\frac{|\ka|}2,\ |\ma|\ge\frac{|\ka|}2\,\}\\
F_\ka&=&\{\ma\in\Z^3_*\,:\,|\ma|<\frac12|\ka|\,\}\\
G_\ka&=&\{\el\in\Z^3_*\,:\,|\el|<\frac12|\ka|\,\},
\end{array}
$$
and we shall estimate the three terms separately. For simplicity, we write
$U_\ka=|\ka|^{2\th(\al)}|u_\ka|$ and $V_\ka=|\ka|^{2\th(\al)}|v_\ka|$, for all $\ka\in\Z^3$,
in such a way that $\sum_\ka|U_\ka|^2=\|u\|_{\th(\al)}^2$.

We start by the term in $E_\ka$. By using H\"older's inequality and Cauchy Schwartz,
\begin{align*}
&\sum_{E_\ka}|\ma|\cdot|u_\el|\cdot|v_\ma|
\le C_\al\sum_{E_\ka}(\frac1{|\ma|^{4\th-1}}+\frac1{|\el|^{4\th-1}})|U_\el|\cdot|V_\ma|\\
&\qquad\le C_\al\|u\|_\th\Bigl(\sum_{E_\ka}\frac1{|\ma|^{8\th-2}}|V_\ma|^2\Bigr)^\frac12
    +C_\al\|v\|_\th\Bigl(\sum_{E_\ka}\frac1{|\el|^{8\th-2}}|U_\el|^2\Bigr)^\frac12,
\end{align*}
and, since the two terms are similar, we handle just the first one. By summing in $\ka$,
and exchanging the sums,
$$
\sum_{\ka\in\Z^3_*}|\ka|^{4\al-1}\sum_{|\ma|\ge\frac{|\ka|}2}\frac1{|\ma|^{8\th-2}}|V_\ma|^2
\le C_0\sum_{\ma\in\Z^3_*}\frac1{|\ma|^{8\th-4\al-4}}|V_\ma|^2
\le C_0\|v\|_{\th},
$$
since there are at most $C_0|\ma|^3$ points $\ka$ of $\Z^3_*$ such that $|\ka|\le2|\ma|$
(where $C_0$ is a universal constant) and, for any value of $\al$, $8\th(\al)-4\al-4\ge0$.

Next, we estimate the term in $F_\ka$. We have
\begin{align*}
&\sum_{\ka\in\Z^3_*}|\ka|^{4\al-1}\bigl|\sum_{\ma\in F_\ka}|\ma|\cdot|v_\ma|\cdot|u_{\ka-\ma}|\bigr|^2=\\
&\qquad=\sum_{\ka\in\Z^3_*}\sum_{\ma_1,\ma_2\in F_\ka}\frac{|\ka|^{4\al-1}}{(|\ma_1|\cdot|\ma_2|)^{2\th-1}}|V_{\ma_1}|\cdot|V_{\ma_2}|\cdot|u_{\ka-\ma_1}|\cdot|u_{\ka-\ma_2}|\\
\intertext{and, by exchanging the sums,}
&\qquad=\sum_{\ma_1,\ma_2\in\Z^3_*}\frac{|V_{\ma_1}|\cdot|V_{\ma_2}|}{(|\ma_1|\cdot|\ma_2|)^{2\th-1}}\sum_{|\ka|>2|\ma_1|\vee|\ma_2|}|\ka|^{4\al-1}|u_{\ka-\ma_1}|\cdot|u_{\ka-\ma_2}|\\
&\qquad\le\Bigl[\sum_{\ma\in\Z^3_*}\frac{|V_\ma|}{|\ma|^{2\th-1}}\bigl(\sum_{|\ka|>2|\ma|}|\ka|^{4\al-1}|u_{\ka-\ma}|^2\bigr)^{\frac12}\Bigr]^2\\
&\qquad\le\|v\|_\th^2\sum_{\ma\in\Z^3_*}\frac1{|\ma|^{4\th-2}}\sum_{|\ka|>2|\ma|}|\ka|^{4\al-1}|u_{\ka-\ma}|^2\\
&\qquad=\|v\|_\th^2\sum_{\el\in\Z^3_*}\frac{|U_\el|^2}{|\el|^{4\th}}\sum_{|\el+\ma|>2|\ma|}\frac{|\el+\ma|^{4\al-1}}{|\ma|^{4\th-2}}.
\end{align*}
It is elementary to see that $\frac1{|\el|^{4\th}}\sum_{|\ma|<|\el|}\frac{|\el+\ma|^{4\al-1}}{|\ma|^{4\th-2}}$
is bounded by a constant depending only on $\al$.

Prior to the estimate of the term in $G_\ka$, we notice that for the terms in $E_\ka$
and $F_\ka$ we only used the fact that $\al>\frac14$. The term in $G_\ka$ is the
most delicate, since it is the only estimate where we need to make special assumptions
for the case $\al=\frac12$. We have
\begin{align*}
&\sum_{\ka\in\Z^3_*}|\ka|^{4\al-1}\bigl|\sum_{\el\in G_\ka}|\ka-\el|\cdot|v_{\ka-\el}|\cdot|u_\el|\bigr|^2\\
&\qquad=\sum_{\ka\in\Z^3_*}\sum_{\el_1,\el_2\in G_\ka}\frac{|\ka|^{4\al-1}}{(|\ka-\el_1|\cdot|\ka-\el_2|)^{2\th-1}}|V_{\ka-\el_1}|\cdot|V_{\ka-\el_2}|\cdot|u_{\el_1}|\cdot|u_{\el_2}|\\
\intertext{and, by exchanging the sums,}
&\qquad=\sum_{\el_1,\el_2\in\Z^3_*}|u_{\el_1}|\cdot|u_{\el_2}|\sum_{|\ka|>2|\el_1|\vee|\el_2|}\frac{|\ka|^{4\al-1}}{(|\ka-\el_1|\cdot|\ka-\el_2|)^{2\th-1}}|V_{\ka-\el_1}|\cdot|V_{\ka-\el_2}|\\
&\qquad\le\sum_{\el_1,\el_2\in\Z^3_*}\!|u_{\el_1}|\!\cdot\!|u_{\el_2}|\Bigl|\sum_{|\ka|>2|\el_1|}\frac{|\ka|^{4\al-1}|V_{\ka-\el_1}|^2}{|\ka-\el_1|^{4\th-2}}\Bigr|^{\frac12}\Bigl|\sum_{|\ka|>2|\el_2|}\frac{|\ka|^{4\al-1}|V_{\ka-\el_2}|^2}{|\ka-\el_2|^{4\th-2}}\Bigr|^{\frac12}\\
&\qquad=\Bigl[\sum_{\el\in\Z^3_*}\frac{|U_\el|}{|\el|^{2\th}}\bigl(\sum_{|\ka|>2|\el|}\frac{|\ka|^{4\al-1}}{|\ka-\el|^{4\th-2}}|V_{\ka-\el}|^2\bigr)^{\frac12}\Bigr]^2\\
&\qquad\le\|u\|_\th^2\sum_{\el\in\Z^3_*}\frac1{|\el|^{4\th}}\sum_{|\ka|>2|\el|}\frac{|\ka|^{4\al-1}}{|\ka-\el|^{4\th-2}}|V_{\ka-\el}|^2\\
&\qquad\le\|u\|_\th^2\sum_{\el\in\Z^3_*}\frac1{|\el|^{4\th}}\sum_{|\el+\ma|>2|\el|}\frac{|\el+\ma|^{4\al-1}}{|\ma|^{4\th-2}}|V_\ma|^2\\
\intertext{by exchanging the sums again and, since the set of all $\el$ such that $|\el+\ma|>2|\el|$ is contained in the set of all $|\el|$
such that $|\el|\le|\ma|$, we have}
&\qquad\le C_\al\|u\|_\th^2\sum_{\ma\in\Z^3_*}\frac{|V_\ma|^2}{|\ma|^{4\th-4\al-1}}\sum_{|\el|<|\ma|}\frac1{|\el|^{4\th}}.
\end{align*}
It is elementary to verify that $|\ma|^{-(4\th-4\al-1)}\sum_{1\le|\el|<|\ma|}|\el|^{-4\th}$ is bounded by
a constant depending only on $\al$ if $\al>\sss{\frac14}$ and $\al\neq{\sss{\frac12}}$. The same is true in the
special case $\al=\sss{\frac12}$ if we evaluate the norm of $B(u,v)$ in $D(A^{\sss{\frac14}-\ep})$.
\end{proof}
\bibliographystyle{amsplain}

\end{document}